\documentclass[11pt]{article}
\usepackage[a4paper,hmargin=26mm,vmargin=25mm]{geometry}
\usepackage[T1]{fontenc}
\usepackage[utf8]{inputenc}
\usepackage{lmodern,amsmath,amssymb,amsthm,mathtools,microtype}
\usepackage{enumitem,needspace}
\usepackage{xcolor}
\usepackage[colorlinks=true,linkcolor=blue!45!black,citecolor=blue!45!black,urlcolor=blue!45!black]{hyperref}
\hypersetup{pdftitle={Sharp bounds and asymptotic profiles for Lane–Emden systems at the Serrin threshold},
 pdfauthor={Shiguang Ma and Shengyang Zang},
 pdfsubject={Source masses, sharp joint growth, and critical system asymptotics}}
\allowdisplaybreaks[2]
\numberwithin{equation}{section}
\newtheorem{theorem}{Theorem}[section]
\newtheorem{lemma}[theorem]{Lemma}
\newtheorem{proposition}[theorem]{Proposition}
\newtheorem{corollary}[theorem]{Corollary}
\theoremstyle{remark}
\newtheorem{remark}{Remark}[section]
\newcommand{\R}{\mathbb R}
\newcommand{\Sph}{\mathbb S^{n-1}}
\newcommand{\dd}{\,d}
\newcommand{\Aop}{\mathcal A}
\newcommand{\av}{\mathop{\mathrm{av}}}
\newcommand{\CQ}{\ensuremath{C_{\mathcal Q}}}

\newcommand{\loc}{\mathrm{loc}}
\title{Sharp bounds and asymptotic profiles for Lane–Emden systems at the Serrin threshold}
\author{Shiguang Ma and Shengyang Zang}
\date{}
\begin{document}
\maketitle
\begin{abstract}
We classify isolated singularities of $p$-Laplace and $k$-Hessian
Lane--Emden systems on their Serrin critical curves without radial
symmetry. Our approach reconstructs a radial profile from the source
measures through a truncated Wolff potential, providing a nonlinear
counterpart of spherical averaging. For upper inequalities at equal
powers, we determine the sharp joint logarithmic lower bound on the
zero-mass singular branch. A radial pair attains this bound while both
normalized components oscillate between zero and infinity. For exact
systems at equal powers, we obtain a common second-order profile with
a shared phase and a sharp two-term logarithmic expansion of the
component difference. At unequal powers, we establish
power--logarithmic asymptotics with explicit leading constants. These
asymptotics hold uniformly in direction.
\end{abstract}

\section{Introduction}\label{sec:introduction}

For the scalar equation $-\Delta w=w^{n/(n-2)}$ in a punctured ball,
$n\ge3$, Aviles~\cite[Theorem~A]{Aviles1987} proved that a positive
solution has either a removable singularity or a logarithmically
corrected fundamental-solution profile. At this Serrin exponent, a
singularity can therefore persist with zero fundamental-solution
coefficient. For nonlinear operators, spherical averaging generally
does not preserve the equation, and coupled components need not
satisfy a closed equation for their sum.

We consider positive solutions of the $p$-Laplace system
\[
 -\Delta_pu=v^{q_1},\qquad -\Delta_pv=u^{q_2},
 \qquad 1<p<n,
\]
and negative $k$-admissible solutions of the Hessian system
\[
 \sigma_k(D^2u)=(-v)^{q_1},\qquad
 \sigma_k(D^2v)=(-u)^{q_2},
 \qquad 1\le k<n/2,
\]
in $B_1\setminus\{0\}$. We classify their isolated singularities for equal
and unequal powers on the respective Serrin critical curves. At the
symmetric point $q_1=q_2=q$, we also study the upper inequalities
obtained by allowing each source to range from zero to the indicated power.

A central observation of this paper is that a truncated Wolff potential~\cite{HW,MH}
provides a nonlinear counterpart of the spherical mean through
reconstruction from the actual source measures. We prescribe the flux
of a radial profile across each sphere to be the total source mass
enclosed. For the Laplacian, this recovers the spherical mean of the
sum up to an additive constant. For the $p$-Laplace and $k$-Hessian
operators, inverting the corresponding radial flux law gives the
Wolff reconstruction. To our knowledge, the use of this reconstruction
as a nonlinear spherical mean to obtain sharp componentwise asymptotics
for critical systems is new.
The total mass then sets a common scale on which the individual masses
describe the components, even when one mass fraction tends to zero.
This yields sharp joint bounds and
cluster-point constraints for upper inequalities, precise synchronization
for exact symmetric systems, and distinct power--logarithmic scales on
the nonsymmetric Serrin curve.

\subsection{The operators and the critical scale}

We use positive functions in both cases, writing $(w_1,w_2)=(u,v)$
for the $p$-Laplacian and $(w_1,w_2)=(-u,-v)$ for the Hessian operator. The parameters associated with the two operators are
\begin{equation}\label{parameters}
(a,\kappa,d_0)=
\begin{cases}
\bigl(p-1,\frac{n-p}{p-1},1\bigr),
& \mathcal Q[w]=-\Delta_p w,\\[1em]
\bigl(k,\frac{n-2k}{k},\frac1k\binom{n-1}{k-1}\bigr),
& \mathcal Q[w]=\sigma_k(D^2(-w)).
\end{cases}
\end{equation}
Here $k$ is an integer, $a$ is the degree of homogeneity in the dependent
variable, $\kappa$ is the fundamental-solution exponent, and $d_0$
normalizes the radial flux. The scalar Serrin exponent is
$q=q_S=n/\kappa$, and the spatial scaling exponent $\ell=n-a\kappa$
equals $p$ for $-\Delta_p$ and $2k$ for $\sigma_k$.
We also set
\[
 m=\frac qa>1,\qquad \beta=\frac1{q-a},\qquad
 L(r)=\log\frac1r,\qquad \omega=|\Sph|.
\]
Thus $\kappa q=n$ and $q\beta-1=a\beta$. In the $p$-Laplace case,
\begin{equation}\label{critical-identities}
 \kappa(q-a)=p,\qquad a(\kappa+1)=n-1,
 \qquad \beta=\frac\kappa p.
\end{equation}
For the Hessian parameters, $\kappa(q-a)=2k$. The borderline identity
$\kappa q=n$ is responsible for the logarithmic correction.

With this notation, both upper systems take the form
\begin{equation}\label{upper-system}
 0\le\mathcal Q[w_1]\le w_2^q,\qquad
 0\le\mathcal Q[w_2]\le w_1^q
 \quad\text{in }B_1\setminus\{0\}.
\end{equation}
For $-\Delta_p$, we consider solutions in
$W^{1,p}_{\loc}\cap L^\infty_{\loc}$ on the punctured ball, with the
inequalities understood distributionally, and use their
$C^{1,\alpha}_{\loc}$ representatives. For $\sigma_k$, we assume
$-w_i\in C^2\cap\Phi^k$ there, where $\Phi^k$ is the class of
$k$-convex functions. In this regularity class, admissibility means
$\sigma_j(D^2(-w_i))\ge0$ for $1\le j\le k$.
Unless stated otherwise, solutions on punctured domains are understood
in these classes.
The Riesz or Hessian measures of the canonical extensions are denoted by
$\mu_i$.

Throughout the paper, a singularity is called removable if the relevant
signed function extends as a locally bounded weak solution in the
$p$-Laplace case, or as a locally bounded $k$-convex solution in the
Hessian case. The extended Hessian equations and inequalities are
understood in the sense of Hessian measures. When both components are
bounded, the extensions are $C^{1,\alpha}$ in the first case and
locally H\"older continuous in the second.

\subsection{Background}\label{sec:background}

For semilinear equations with critical Sobolev growth, Caffarelli,
Gidas and Spruck~\cite[Theorem~1.1]{CGS} proved asymptotic radiality;
see Badiale and Bardi~\cite{BB} for related nonlinear results.
At the scalar Serrin exponent, Chan and
DelaTorre~\cite{CD2020} constructed radial solutions with higher-order
logarithmic corrections.

Taliaferro~\cite{T01,T06} treated scalar semilinear inequalities.
For the $p$-Laplace inequality $0\le-\Delta_pw\le\tau w^q$ at the Serrin
exponent, our earlier work~\cite{MZ} classified isolated singularities,
established asymptotic radiality, and proved on the zero-mass singular branch
\begin{equation}\label{scalar-lower-background}
 \liminf_{x\to0}|x|^\kappa L(|x|)^\beta w(x)
 \ge\tau^{-\beta}(a\beta\kappa^a)^\beta.
\end{equation}
Our mass argument also recovers the exact scalar logarithmic limits
of \cite[Theorem~1]{MZII} and \cite[Theorem~1.3]{MZKH} for classical
$p$-Laplace and $k$-Hessian solutions, respectively
(Corollary~\ref{hessian-scalar}).

We use Trudinger and Wang's Hessian measure theory~\cite{TW2,TW3}
and the potential estimates of Kilpel\"ainen and Mal\'y~\cite{KM94},
Labutin~\cite{Labutin}, and Phuc and Verbitsky~\cite{PV}.
For upper inequalities, we apply these estimates to the actual source
measures, since the power bounds need not be integrable.
For cross-coupled systems with nonnegative measure
data, Bidaut-V\'eron, Nguyen and V\'eron~\cite{BNV} proved existence
and potential estimates; their Theorems~C and~E exclude positive Dirac
data at the symmetric Serrin powers.

For Laplace upper systems, Ghergu, Taliaferro and
Verbitsky~\cite[Theorem~3.1]{GTV} proved fundamental-scale bounds;
applying \cite[Theorem~3.1]{T06} to the sum gives asymptotic radiality
and a zero-mass lower bound. For exact symmetric systems,
Bidaut-V\'eron~\cite[Proposition~2.1(v) and Theorem~1.2(i)]{BV2000}
and Liu and Liang~\cite[Theorems~1.1 and~1.4]{LL} obtained mass
exclusion and further estimates. Their precise classifications
\cite[Theorem~1.3]{BV2000}, \cite[Theorem~1.3]{LL}
exclude the Serrin boundary.

On the Serrin curve, Garc\'ia-Huidobro and
Yarur~\cite[Theorems~3.1--3.2 and Remark~3.1]{GHY1997} proved the
common leading term for radial Laplace systems at equal powers, and
two-sided power--logarithmic bounds at unequal powers, identifying the
constants if either normalized component converges. We obtain higher-order
asymptotics at equal powers and convergence to both constants at unequal
powers, uniformly in direction without radial symmetry.

\subsection{Source masses and precise asymptotics}

Let $F(r)=\sum_i\mu_i(B_r)$, including any mass at the origin. The nonincreasing
radial profile with enclosed mass $F$ satisfies
$\omega d_0r^{a\kappa}(-rz')^a=F(r)$, hence
\begin{equation}\label{intro-source-reconstruction}
 z(r)=K_0+\int_r^R
 \left(\frac{\sum_i\mu_i(B_s)}{\omega d_0s^{a\kappa}}\right)^{1/a}
 \frac{\dd s}{s}.
\end{equation}
Here $R>0$ is fixed and small and $K_0=z(R)>0$ is sufficiently large.
For $a=1$, Gauss--Green gives
$z(r)-z(R)=\bar S(r)-\bar S(R)$, where $\bar S$ is the spherical mean
of $\sum_iw_i$ (Remark~\ref{linear-consistency}).
The integral is $(\omega d_0)^{-1/a}$ times the Wolff
potential~\cite{HW,MH} centered at $0$, restricted to $[r,R]$,
with parameters $(1,p)$ or $(2k/(k+1),k+1)$, respectively
(see \cite[p.~900]{PV}). Differentiating the flux identity gives
\[
 \mathcal Q[z](r)=\av_{\Sph}\sum_i\mathcal Q[w_i](r\theta)
 \quad\text{for a.e. }r,
\]
where $\av_{\Sph}$ is normalized spherical averaging.

Suppose the total source is bounded above by $C(\sum_iw_i)^q$, all point
masses vanish, and the sum is unbounded near $0$. With
$P_i(t)=\mu_i(B_{e^{-t}})/(\omega d_0)$ and $J=\sum_iP_i$,
Proposition~\ref{mass-approx} gives
\begin{equation}\label{intro-mass-approximation}
 e^{-\kappa t}w_i(e^{-t}\theta)
 =\kappa^{-1}P_i(t)^{1/a}+o(J(t)^{1/a}),
\end{equation}
uniformly in direction, even when $P_i/J\to0$. Also,
$z(r)\sim(\sum_iw_i(r\theta)^a)^{1/a}$ uniformly in $\theta$.
This extends to any fixed finite number of components.
Unlike the spherical-extrema and bounded-variation phase-plane
methods of \cite{MZ,MZII,MZKH}, source masses retain individual
contributions and give sharp joint bounds without convergence
of the relative masses.

\subsection{Main results}

The sharp joint constant for~\eqref{upper-system} is determined by a
one-variable minimization. For $0\le s\le1$, set
\begin{gather}
 I_m(s)=\int_1^\infty
 (y^{m+1}-1+s^{m+1})^{-m/(m+1)}\dd y,\label{Idef}\\[4pt]
 \Phi_{a,\kappa}(s)=\kappa^{q\beta-1}
 (1+s^{1/a})I_m(s)^\beta,\label{Phi-def}\\[4pt]
 C_0=(d_0a\beta\kappa^a)^\beta,\qquad
 \CQ=d_0^\beta\min_{0\le s\le1}\Phi_{a,\kappa}(s).\label{optimalconstant}
\end{gather}
Here $C_0$ is the scalar constant and $\CQ$ is the optimal joint constant.
The function $I_m$ is finite and continuous on $[0,1]$.
When $a\ge1$ (all Hessian cases and $p\ge2$), the minimum is attained
uniquely at $s=0$, and $\CQ<2C_0$. Along sequences attaining the sharp
bound, the smaller-to-larger source-mass ratio then tends to zero
(Remark~\ref{mass-ratio-refinement}).
For $1<p<2$, there is a unique $p_*(n)\in(1,2)$ satisfying
\begin{equation}\label{subquadratic-threshold}
 (n-p_*)(4n-p_*)=(p_*-1)(2n-p_*)^2.
\end{equation}
The unique minimizer is $s=1$ for $1<p\le p_*(n)$ and
$s_*\in(0,1)$ for $p_*(n)<p<2$; correspondingly,
$\CQ=2C_0$ in the first range and $\CQ<2C_0$ in the second
(Lemma~\ref{minimum-subquadratic}).
For the Laplacian, the strict lower comparison $C_0<\CQ$ follows from
$I_q(0)>\int_1^\infty y^{-q}\dd y=1/(q-1)$.
The point masses determine the fundamental coefficients
\begin{equation}\label{fundamental-coefficients}
 b_i=\frac1\kappa
 \left(\frac{\mu_i(\{0\})}{\omega d_0}\right)^{1/a}\ge0.
\end{equation}
A positive function $W$ is asymptotically radial if
\begin{equation}\label{joint-asymptotic-radiality}
 \frac{\max_{|x|=r}W(x)}{\min_{|x|=r}W(x)}\to1.
\end{equation}

\begin{theorem}
\label{mainjoint}\label{hessian-joint}
Let $(w_1,w_2)$ solve \eqref{upper-system} for either operator.
Then $S=w_1+w_2$ is asymptotically radial and, uniformly in $\theta$,
\begin{equation}\label{component-fundamental-limits}
 r^\kappa(w_1,w_2)(r\theta)\to(b_1,b_2).
\end{equation}
Exactly one of the following holds:
\begin{enumerate}[label=(\roman*),itemsep=2pt]
\item both singularities are removable;
\item $b_1+b_2>0$;
\item $b_1=b_2=0$, $w_1(x),w_2(x)\to\infty$ as $x\to0$, and
\begin{gather}
 \liminf_{x\to0}|x|^\kappa L(|x|)^\beta S(x)\ge \CQ,
 \label{joint-lower}\\[3pt]
 \liminf_{x\to0}|x|^\kappa L(|x|)^\beta
 \max\{w_1(x),w_2(x)\}\ge C_0,
 \label{joint-max-lower}\\[3pt]
 \limsup_{x\to0}|x|^\kappa L(|x|)^\beta
 \min\{w_1(x),w_2(x)\}\ge C_0.
 \label{joint-min-upper}
\end{gather}
\end{enumerate}
All three bounds in alternative~(iii) are sharp.
\end{theorem}

Corollary~\ref{finite-cluster-region} determines the exact finite cluster
region of the normalized components in alternative~(iii) through the same
blow-up-time functional that gives \CQ. By Proposition~\ref{boundary-arc},
a solution attaining a sharp boundary state also attains every state on
the boundary arc joining it to $(C_0,C_0)$.

Every nonzero $(b_1,b_2)\in[0,\infty)^2$, including $b_1=0<b_2$,
is realized by $w_i=1+b_i r^{-\kappa}$, with
$\mu_i=\omega d_0(\kappa b_i)^a\delta_0$.
On branch~(iii), the next theorem gives joint attainment by a radial pair
whose individual normalized lower limits both vanish.

\begin{theorem}\label{individual-failure}
For either operator, there is a smooth radial solution of
\eqref{upper-system} on $B_{e^{-1}}\setminus\{0\}$ such that
$w_1,w_2\to\infty$ and $r^\kappa(w_1+w_2)\to0$ as $r\to0$.
Moreover, for $i=1,2$,
\begin{gather}
 \liminf_{r\to0}r^\kappa L(r)^\beta w_i(r)=0,\qquad
 \limsup_{r\to0}r^\kappa L(r)^\beta w_i(r)=\infty,
 \label{individual-oscillation}\\[3pt]
 \liminf_{r\to0}r^\kappa L(r)^\beta(w_1(r)+w_2(r))=\CQ.
 \label{simultaneous-joint-attainment}
\end{gather}
\end{theorem}

On the singular branch, exact symmetric coupling rules out this
imbalance and forces the two components to synchronize. Consider
\begin{equation}\label{exact-symmetric}
 \mathcal Q[w_1]=w_2^q,\qquad \mathcal Q[w_2]=w_1^q
 \quad\text{in }B_1\setminus\{0\},
\end{equation}
and set $\gamma=q\beta$, so that $\gamma=n/p$ for $-\Delta_p$
and $\gamma=n/(2k)$ for $\sigma_k$.

\begin{theorem}
\label{mainexact}\label{hessian-exact}
A solution of \eqref{exact-symmetric} either has two removable
singularities or satisfies, uniformly in $\theta$,
\begin{equation}\label{exact-common-limit}
 r^\kappa L(r)^\beta(w_1,w_2)(r\theta)\to C_0(1,1).
\end{equation}
More precisely, there is a common constant $b_*\in\R$ such that,
for $i=1,2$,
\begin{equation}\label{exact-common-phase}
 r^\kappa L(r)^\beta w_i(r\theta)
 =C_0\biggl[1-\frac{\beta\gamma}{\kappa}\frac{\log L(r)}{L(r)}
 +\frac{b_*}{L(r)}
 +O\left(\frac{(\log L(r))^2}{L(r)^2}\right)\biggr].
\end{equation}
Moreover, for some $c\in\R$,
\begin{equation}\label{nonradial-second-order-limit}
 w_1(r\theta)-w_2(r\theta)
 =cL(r)^{-\gamma}
 -\frac{c\gamma^2}{\kappa}L(r)^{-\gamma-1}\log L(r)
 +O\bigl(L(r)^{-\gamma-1}\bigr).
\end{equation}
All these asymptotics are uniform in $\theta$. If $c=0$, then
$w_1(r\theta)-w_2(r\theta)=O(r^\delta)$ for every $0<\delta<1$,
with the implied constant allowed to depend on $\delta$.
\end{theorem}

Proposition~\ref{nonzero-radial-coefficient} constructs radial
nondiagonal solutions and proves that, in the radial class, $c=0$ exactly
when $w_1=w_2$. Thus the logarithmic correction in
\eqref{nonradial-second-order-limit} cannot be absorbed into
$O(L^{-\gamma-1})$. For $c\ne0$, the difference extends continuously
by zero but has no positive H\"older exponent there. Individual derivative
limits are given in \eqref{exact-derivative-limits}.

\begin{remark}\label{three-dimensional-example}
For $(n,p)=(3,2)$, one has $q=3$, $\kappa=1$, $\beta=1/2$ and
$\gamma=3/2$. The joint and scalar constants satisfy
\[
 \CQ=\left[\frac14 B(1/2,1/4)\right]^{1/2}
 \simeq1.145<\sqrt2=2C_0,
\]
where $B$ is the Euler beta function; see \eqref{beta-constant}.
On the zero-mass singular branch, the least normalized lower limit of
the sum under the upper inequalities is therefore strictly smaller than
its limit $\sqrt2$ for the exact symmetric system. Exact symmetric
coupling also gives equal leading coefficients and a difference of
order $L(r)^{-3/2}$ when $c\ne0$.
\end{remark}

\medskip

For $q_1,q_2>0$ with $q_1q_2>a^2$, the critical boundary is determined
by the coupled scaling of the exact system
\begin{equation}\label{nonsymmetric-system}
 \mathcal Q[w_1]=w_2^{q_1},\qquad
 \mathcal Q[w_2]=w_1^{q_2}
 \quad\text{in }B_1\setminus\{0\}
\end{equation}
with homogeneities
\[
 \alpha_1=\frac{\ell(a+q_1)}{q_1q_2-a^2},\qquad
 \alpha_2=\frac{\ell(a+q_2)}{q_1q_2-a^2}.
\]
The Serrin boundary is $\max\{\alpha_1,\alpha_2\}=\kappa$; its diagonal
point corresponds to~\eqref{exact-symmetric}. Up to interchanging the
components, the remaining branch is $\alpha_1=\kappa>\alpha_2$.

\begin{theorem}
\label{main-nonsymmetric}
For either operator, let $(w_1,w_2)>0$ solve
\eqref{nonsymmetric-system} with
\begin{equation}\label{nonsymmetric-parameters}
 q_1>q_S,\qquad
 \eta=\frac{n}{q_1},\qquad
 q_2=\frac{\ell+a\eta}{\kappa},\qquad
 D=q_1q_2-a^2.
\end{equation}
Then either both singularities are removable or, uniformly in $\theta$,
\begin{equation}\label{nonsymmetric-limits}
 \begin{split}
 w_1(r\theta)&=A\,r^{-\kappa}L(r)^{-a/D}(1+o(1)),\\
 w_2(r\theta)&=B\,r^{-\eta}L(r)^{-q_2/D}(1+o(1)),
 \end{split}
\end{equation}
where
\begin{equation}\label{nonsymmetric-constants}
 \begin{gathered}
 C_1=\frac{d_0a^2\kappa^a}{D},\qquad
 C_2=d_0a\eta^a(\kappa-\eta),\\
 A=(C_1^aC_2^{q_1})^{1/D},\qquad
 B=(C_1^{q_2}C_2^a)^{1/D}.
 \end{gathered}
\end{equation}
\end{theorem}

Every parameter choice is realized by a smooth radial solution
(Section~\ref{sec:nonsymmetric}), including $q_2<1$ or $q_2<a$.
Corollary~\ref{nonsymmetric-masses-gradients} gives the corresponding
mass asymptotics. The diagonal endpoint $\eta=\kappa$ belongs to
Theorem~\ref{mainexact}. For $p=2$, $n\ge5$ and $q_2=1$, one has
$q_1=n/(n-4)$ and the system becomes
$\Delta^2w_2=w_2^{n/(n-4)}$, $w_1=-\Delta w_2>0$.
For this scalar biharmonic case, Soranzo~\cite[Theorem~5]{Soranzo1997}
proved radial logarithmic upper bounds, and Andrade and
do \'O~\cite[Theorem~2(b)]{AD2022} treated nonradial logarithmic
asymptotics. The constants in~\eqref{nonsymmetric-constants} are
derived directly from the coupled equations.

Corollary~\ref{mainvariable} extends the leading symmetric asymptotics
to bounded nonnegative coupling coefficients whose spherical averages
have positive logarithmic Ces\`aro limits. These averages determine the
leading constants, allowing persistent oscillation and no positive
pointwise lower bound.

\subsection{Proof strategy and organization}

Sections~\ref{sec:local}--\ref{sec:mass} derive fundamental-scale bounds
by local source finiteness, potential estimates and doubling. On the zero-mass
branch, radial comparison and slow flux variation give
\eqref{intro-mass-approximation}. We identify limits by annular $C^1$ compactness
and flux identities for $-\Delta_p$, or local $L^1$ convergence across
the origin and weak measure continuity for $\sigma_k$.

Sections~\ref{sec:sharp}--\ref{sec:construction} use the blow-up time of
$\dot X=Y^m$, $\dot Y=X^m$ for the sharp bound, minimizing ratios,
boundary-arc propagation and realizability. Comparison across the puncture
in Section~\ref{sec:coupling} bounds the difference for exact symmetric
systems. Section~\ref{sec:exact}
derives the common profile and synchronization via small-perturbation
Hessian regularity~\cite{Fan} and the spherical spectral gap.
Angular decay leaves exponentially small errors in the
averaged flux equations. Separating angular forcing from errors
proportional to the difference gives every decay exponent below $1$
when its leading coefficient vanishes. Monotone-weight cancellation
treats averaged coefficients.

Section~\ref{sec:nonsymmetric} reduces the fundamental-scale component
to a scalar critical inequality by doubling and concave-power comparison
(cf.~\cite[Theorem~1.2(i)]{BV2000}). Its mass asymptotics recover the
other component by radial source comparison; a dynamical formulation
following~\cite{BG} gives realization. Appendix~\ref{app:cylinder} proves the cylinder inverse
estimate.\label{sec:introduction-end}

\section{Potential estimates and fundamental coefficients}\label{sec:local}

For either operator, homogeneity of degree $a$ removes constant positive
coupling coefficients. If
$0\le\mathcal Q[w_1]\le\tau_1w_2^q$ and
$0\le\mathcal Q[w_2]\le\tau_2w_1^q$, set $E=q^2-a^2$ and
\begin{equation}\label{constant-coupling-normalization}
 \widetilde w_1=\tau_1^{a/E}\tau_2^{q/E}w_1,\qquad
 \widetilde w_2=\tau_1^{q/E}\tau_2^{a/E}w_2.
\end{equation}
Then $(\widetilde w_1,\widetilde w_2)$ satisfies \eqref{upper-system}, with
coupling coefficients $1$. A common
coefficient $\tau$ multiplies both normalized components by $\tau^\beta$,
so the sharp lower constant for $w_1+w_2$ on the zero-mass singular branch
is $\tau^{-\beta}\CQ$.

For a nonnegative measure $\nu$, write the Wolff potentials as
\begin{equation}\label{common-potential-definition}
 \mathcal P_{\mathcal Q}^{R}\nu(x)
 =\int_0^R
 \left(\frac{\nu(B_s(x))}{s^{a\kappa}}\right)^{1/a}\frac{\dd s}{s}.
\end{equation}
If $w\ge0$ and the relevant signed function is $p$-superharmonic or
$k$-convex on $B_{3R}(x)$, its source measure satisfies
\begin{equation}\label{common-potential}
 c\mathcal P_{\mathcal Q}^{R/8}\mu_w(x)\le w(x)
 \le C\inf_{B_R(x)}w+C\mathcal P_{\mathcal Q}^{2R}\mu_w(x).
\end{equation}
For $-\Delta_p$ this is the estimate of Kilpel\"ainen and
Mal\'y \cite[Theorem~1.6]{KM94}; for the Hessian operator it is
Labutin's estimate \cite{Labutin}, in the form of
\cite[Theorem~7.2]{PV}. Constants may depend on the operator parameters.

\subsection{Riesz measures and \texorpdfstring{$p$}{p}-Laplace regularity}

In the $p$-Laplace case, a nonnegative $p$-superharmonic
function on a punctured ball extends across $0$, since a point has zero
$p$-capacity when $p<n$; its canonical value is
$w(0)=\liminf_{x\to0}w(x)$ \cite[Theorem~7.35]{HKM}.
Nonnegativity supplies the required lower bound.
The extended Riesz measure is finite on compact subballs
\cite[Theorem~2.1]{KM92}. We use the generalized gradient of truncation
theory in its Riesz equation. On the punctured ball it agrees with the
ordinary weak gradient, and its flux is locally integrable even when the
function itself is not \cite[Theorem~1.15]{KM92}. For each component,
\begin{equation}\label{measure-decomposition}
 \mu_i=f_i\dd x+\mu_i(\{0\})\delta_0,\qquad
 f_i=-\Delta_pw_i\text{ away from }0,\qquad f_i\in L^1_{\loc}(B_1).
\end{equation}
Under the almost-everywhere upper inequalities, the density is bounded
on compact subsets of the punctured ball.

Bounded densities also give
\begin{equation}\label{gradient-input}
 \|\nabla w\|_{L^\infty(B_{r/2})}
 \le C\left(r^{-1}\|w\|_{L^\infty(B_r)}
       +(r\|f\|_{L^\infty(B_r)})^{1/a}\right),
 \qquad -\Delta_pw=f.
\end{equation}
The interior $C^{1,\alpha}$ estimate for bounded measurable
right-hand sides \cite[Theorem~1, p.~127]{Tolksdorf}, after scaling,
gives \eqref{gradient-input} and local $C^1$ compactness for families
with uniformly bounded solutions and densities. We write
$\av_{\Sph}g=\omega^{-1}\int_{\Sph}g\dd\theta$.

\subsection{Hessian measures and regularity}\label{sec:hessian}

In the $k$-Hessian case, write $h=-w\le0$,
$d=d_0$, and $F_k[h]=\sigma_k(D^2h)$. Hessian measures are locally
finite and weakly continuous under local
$L^1$ convergence in $\Phi^k$
\cite[Theorem~1.1 and Lemma~2.2]{TW2}. We also use uniqueness for
$\mu_k[h]=M\delta_0$ in $\R^n$, with $h\to0$ at
infinity, stated after the proof of
\cite[Theorem~4.5, equations~(4.14)--(4.15)]{TW3}:
\begin{equation}\label{hessian-fundamental-measure}
 h(x)=-\frac1\kappa\left(\frac{M}{\omega d}\right)^{1/k}
 |x|^{-\kappa}.
\end{equation}
For $M=0$, this follows from
\eqref{common-potential} by letting $R\to\infty$.

For $h\le0$ with
$h\in C^2(B_1\setminus\{0\})\cap\Phi^k(B_1\setminus\{0\})$,
\cite[Lemma~2.4]{MZKH} gives the canonical extension
$h(0)=\limsup_{x\to0}h(x)$ in $\Phi^k(B_1)$ and the decomposition
\begin{equation}\label{hessian-decomposition}
 \mu_k[h]=M\delta_0+f\dd x,\qquad f\in L^1_{\loc}(B_1),
 \qquad f=F_k[h]\text{ away from }0,
\end{equation}
where $M=\mu_k[h](\{0\})\ge0$.
Thus, for both operators, the actual densities are integrable near $0$.
This does not assert integrability of the upper bound $S^q$.

For compactness with bounded Hessian densities, suppose
$\mu_k[h]=f\dd x$ on $B_{7r}$. Set $H=\sup_{B_{7r}}h$ and apply
\eqref{common-potential} to the nonnegative function $H-h$ at
$x\in B_r$, with $R=2r$. Its signed $k$-convex function is $h-H$,
with the same Hessian measure $f\dd x$. Since $B_{2r}(x)$ contains
$B_r$, taking the supremum over $x\in B_r$ gives
\[
 H-\inf_{B_r}h\le C(H-\sup_{B_r}h)+Cr^2\|f\|_\infty^{1/k}.
\]
Let $O=\operatorname{osc}_{B_{7r}}h$, $o=\operatorname{osc}_{B_r}h$
and $A=H-\sup_{B_r}h\le O-o$. Enlarging $C$ so that $C>1$,
the inequality becomes $A+o\le CA+Cr^2\|f\|_\infty^{1/k}$.
Consequently, with $\eta=(C-1)/C\in(0,1)$,
\begin{equation}\label{hessian-oscillation}
 \operatorname{osc}_{B_r}h\le
 \eta\operatorname{osc}_{B_{7r}}h+Cr^2\|f\|_\infty^{1/k}.
\end{equation}
Iteration gives a uniform local H\"older modulus for locally bounded
$k$-convex representatives with bounded Hessian density.

\subsection{Elementary potential bounds}

\begin{lemma}\label{potential-elementary}
Let $\nu$ be a finite nonnegative measure of mass $M$ on $\R^n$.
Then
\begin{subequations}\label{hessian-potential-elementary}
\begin{equation}\label{weak-wolff}
 |\{x:\mathcal P_{\mathcal Q}^{\infty}\nu(x)>t\}|
 \le C M^{q/a}t^{-q}\qquad(t>0).
\end{equation}
If $0\le f\in L^1\cap L^\infty$, then
\begin{equation}\label{density-wolff}
 \|\mathcal P_{\mathcal Q}^{\infty}(f\dd x)\|_\infty
 \le C\|f\|_1^{\ell/(na)}\|f\|_\infty^{\kappa/n}.
\end{equation}
\end{subequations}
In particular, each canonically extended positive component belongs to
$L^{q,\infty}_{\loc}$.
\end{lemma}

\begin{proof}
Let $\mathcal M\nu(x)=\sup_{s>0}s^{-n}\nu(B_s(x))$. Using
$\nu(B_s(x))\le\min\{M,s^n\mathcal M\nu(x)\}$ in
\eqref{common-potential-definition}, and splitting at
$s=(M/\mathcal M\nu(x))^{1/n}$, gives
\[
 \mathcal P_{\mathcal Q}^{\infty}\nu(x)
 \le C M^{\ell/(na)}(\mathcal M\nu(x))^{\kappa/n}.
\]
The weak $(1,1)$ estimate for $\mathcal M$ yields
\eqref{weak-wolff}, since $q=n/\kappa$ and
$1+\ell/(a\kappa)=q/a$. For $\nu=f\dd x$, replacing
$\mathcal M\nu$ by $C\|f\|_\infty$ proves \eqref{density-wolff}.
The case $M=0$ is immediate.

For the local assertion, restrict the measure to a slightly larger compact
ball and extend it by zero. Cover the smaller ball by finitely many balls
where the upper estimate in \eqref{common-potential} applies. On each,
fix a point where the function is finite, and choose the potential radius
so that every infimum ball contains that point. Such points exist since
the signed representative is finite almost everywhere. The infimum terms
are uniformly bounded, and \eqref{weak-wolff} controls the potentials.
Taking the finite union gives the local weak $L^q$ bound.
\end{proof}

\subsection{Fundamental-solution bounds and coefficients}

\begin{lemma}
\label{upper-bound}\label{fundamental-vector}\label{hessian-fundamental}
Let $N\ge1$ be fixed and let $w_1,\ldots,w_N>0$ satisfy, in
$B_1\setminus\{0\}$,
\begin{equation}\label{upper-with-C}
 0\le f_i=\mathcal Q[w_i]\le CS^q,
 \qquad S=\sum_iw_i,\qquad 1\le i\le N.
\end{equation}
With $b_i$ defined by \eqref{fundamental-coefficients},
\begin{equation}\label{uniform-fundamental-limit}
 \sup_{\theta\in\Sph}|r^\kappa w_i(r\theta)-b_i|\to0
 \qquad(r\to0).
\end{equation}
In particular, $S(x)=O(|x|^{-\kappa})$, and each component with
$b_i>0$ is asymptotically radial.
\end{lemma}

\begin{proof}
\noindent\emph{The scale bound.}
The extensions have locally finite measures and integrable densities.
Lemma~\ref{potential-elementary} gives, on a fixed ball $B_{r_0}$ about $0$,
\begin{equation}\label{common-weak-bound}
 |\{x\in B_{r_0}:S(x)>t\}|\le A t^{-q},\qquad t>0.
\end{equation}
Here $\{S>t\}\subset\bigcup_i\{w_i>t/N\}$, so summing the
individual weak bounds requires no assumption that $q\ge1$.
If the upper bound failed, we could choose $z_j\to0$ such that
$2|z_j|<r_0$ and, for $M=S^{1/\kappa}$,
$|z_j|M(z_j)>4j$. Apply the doubling lemma
\cite[Lemma~5.1]{PQS} in
$D_j=B_{2|z_j|}\setminus\{0\}$. It gives $x_j\in D_j$ such that
\[
 M(x_j)\ge M(z_j),\qquad
 M(x_j)\operatorname{dist}(x_j,\partial D_j)>2j,
 \qquad
 M(x)\le2M(x_j)\quad\text{if }|x-x_j|\le j/M(x_j).
\]
The selection lemma applies since $M$ is positive and bounded on
every compact subset of $D_j$.
In particular $x_j\to0$. With
\[
 S_j=S(x_j),\qquad \rho_j=S_j^{-1/\kappa},\qquad
 \widetilde w_{i,j}(y)=S_j^{-1}w_i(x_j+\rho_jy),
\]
we have $\rho_j/|x_j|<1/(2j)$, the original selected balls avoid $0$,
and
\begin{equation}\label{common-doubling-normalization}
 \sum_i\widetilde w_{i,j}(0)=1,\qquad
 0\le\sum_i\widetilde w_{i,j}\le2^\kappa\quad\text{on }B_j.
\end{equation}
Homogeneity under scaling of the functions and spatial variables gives
\begin{equation}\label{scaled-bounds}
 g_{i,j}(y):=\mathcal Q[\widetilde w_{i,j}](y)
 =S_j^{-a}\rho_j^\ell f_i(x_j+\rho_jy)
 =\rho_j^n f_i(x_j+\rho_jy),\qquad 0\le g_{i,j}\le C_1.
\end{equation}
The last bound follows from \eqref{upper-with-C},
$\rho_j^nS_j^q=1$, and \eqref{common-doubling-normalization}.
For each fixed $R>0$,
\[
 \int_{B_{2R}}g_{i,j}(y)\dd y
 =\int_{x_j+\rho_jB_{2R}}f_i(x)\dd x\to0.
\]
These balls avoid $0$ and shrink to it, so integrability
of $f_i$ suffices; $S^q\in L^1$ is not required.
Moreover \eqref{common-weak-bound} is invariant under this scaling,
since $\rho_j^{-n}S_j^{-q}=1$. Thus, for each fixed $R$ and all large $j$,
\begin{equation}\label{scaled-weak-bound}
 |\{y\in B_R:\sum_i\widetilde w_{i,j}(y)>t\}|
 \le At^{-q},\qquad
 \sum_i\inf_{B_R}\widetilde w_{i,j}\le CR^{-\kappa}.
\end{equation}
Apply \eqref{common-potential} to each component at $0$, taking
$j>3R$, and estimate the potentials by \eqref{density-wolff}, with
the densities restricted to $B_{2R}$ and extended by zero. We obtain
\[
 1\le CR^{-\kappa}
   +C\sum_i\|g_{i,j}\|_{L^1(B_{2R})}^{\ell/(na)}.
\]
The constants are independent of $R$ and $j$. Letting $j\to\infty$
for fixed $R$, then $R\to\infty$, gives a contradiction.

\smallskip
\noindent\emph{Identification of the coefficients.}
For an arbitrary decreasing sequence $r_j\to0$, set
$v_{i,j}(y)=r_j^\kappa w_i(r_jy)$. The upper bound gives
$0\le v_{i,j}(y)\le C|y|^{-\kappa}$. On each compact annulus,
the scaled densities $r_j^nf_i(r_jy)$ are uniformly bounded, and
their integrals tend to zero.

For $\mathcal Q=-\Delta_p$, local $C^{1,\alpha}$ estimates give,
after extraction, $C^1$ convergence on compact annuli
to a nonnegative $p$-harmonic function $v_i$. Its whole-space bound
$v_i(y)\le C|y|^{-\kappa}$ and the classification
\cite[Theorem~5.3]{MZ} imply
$v_i(y)=A_i|y|^{-\kappa}$, with $A_i\ge0$; decay at infinity excludes
the additive constant. Choose $\rho\in(1,2)$ for which the
Gauss--Green identities hold simultaneously for all components and
all members of the countable scaled sequence. Then
\[
 \mu_i(B_{r_j\rho})
 =-\int_{\partial B_\rho}
 |\nabla v_{i,j}|^{p-2}\partial_\nu v_{i,j}\dd S
 \to\omega(\kappa A_i)^a.
\]
The left side tends to $\mu_i(\{0\})$. Hence $A_i=b_i$.
This argument uses annular $C^1$ convergence, without assuming local
$L^1$ convergence across $0$.

For $\mathcal Q[w]=F_k[-w]$, estimate \eqref{hessian-oscillation}
gives locally uniform convergence of $h_{i,j}=-v_{i,j}$ away from $0$ to
$h_i$. Since $\kappa<n$, the bound
$|h_{i,j}(y)|\le C|y|^{-\kappa}$ also gives local $L^1$ convergence
across $0$. The scaled Hessian measures satisfy
$\mu_k[h_{i,j}](E)=\mu_i(r_jE)$ and converge weakly to
$\mu_i(\{0\})\delta_0$. Weak continuity therefore identifies
this limit as $\mu_k[h_i]$. Since $h_i\to0$ at infinity,
\eqref{hessian-fundamental-measure} gives
$h_i=-b_i|y|^{-\kappa}$.

For both operators the limit is independent of the subsequence.
Since $r_j\to0$ was arbitrary, the coefficient limits are uniform.
The spherical ratio follows when $b_i>0$.
\end{proof}

\begin{lemma}\label{vanishing}
Under the hypotheses of Lemma~\ref{fundamental-vector} with $N=2$,
if $\mathcal Q[w_1]\ge c_1w_2^q$ and
$\mathcal Q[w_2]\ge c_2w_1^q$ for $c_1,c_2>0$, then
$\mu_1(\{0\})=\mu_2(\{0\})=0$.
\end{lemma}
\begin{proof}
A positive first mass would give $w_1(x)\ge c|x|^{-\kappa}$ by
Lemma~\ref{fundamental-vector}, and hence
$f_2\ge c_2w_1^q\ge c'|x|^{-n}$, contrary to $f_2\in L^1$ near $0$.
We exclude the second mass by interchanging the components.
\end{proof}

\section{Total source mass and asymptotic radiality}\label{sec:mass}

Fix either operator in~\eqref{parameters} and its associated parameters.
Let $w_i>0$ have source measures $\mu_i$, where $1\le i\le N$ for fixed finite $N$.
Write $S=\sum_iw_i$; constants may depend on $N$. We compare the branch
with vanishing fundamental-solution terms, without presupposing the joint
lower bound or component comparison. The next lemma absorbs the small
neighboring-annulus term.

\begin{lemma}\label{dyadic}
Suppose $h_j\ge0$, $h_j\le K_0b^j$ for some $b>1$, and, for all $j$ large,
\[
 h_j\le C_0+\varepsilon_j\max(h_{j-1},h_j,h_{j+1}),\qquad
 \varepsilon_j\to0.
\]
Then $\limsup_{j\to\infty}h_j\le C_0$.
\end{lemma}
\begin{proof}
We fix $\varepsilon>0$ with $\varepsilon b<1$ and choose $J$ so that
$\varepsilon_j\le\varepsilon$ for $j>J$. Iterating the inequality until the
index first reaches $J$, or until $k$ steps have been made, gives
\[
 h_j\le\frac{C_0}{1-\varepsilon}
       +\varepsilon^{j-J}h_J
       +\varepsilon^k\max_{J\le i\le j+k}h_i.
\]
Indeed, the a priori growth gives
\[
 \varepsilon^k\max_{J\le i\le j+k}h_i
 \le K_0b^j(\varepsilon b)^k,
\]
which tends to zero as $k\to\infty$ since $\varepsilon b<1$.
Then let $j\to\infty$, and finally $\varepsilon\to0$.
\end{proof}

\subsection{Comparison with a radial function}

For a radial profile $z$, define its flux by
\begin{equation}\label{radial-flux}
 J_z(r)=
 \begin{cases}
 r^{n-1}|z'(r)|^{p-2}(-z'(r)),&\mathcal Q=-\Delta_p,\\
 r^{n-k}(-z'(r))^k,&\mathcal Q[z]=\sigma_k(D^2(-z)).
 \end{cases}
\end{equation}
The radial identities give
\begin{equation}\label{common-radial-identity}
 \mathcal Q[z]=d_0r^{1-n}J_z'(r).
\end{equation}
This identity holds almost everywhere for $C^1$ functions with locally
absolutely continuous flux. The signed Hessian formula is
\begin{equation}\label{hessian-radial-identity}
 F_k[\phi]=d_0r^{1-n}\big(r^{n-k}(\phi')^k\big)'.
\end{equation}
\begin{lemma}\label{radial-admissibility}
A radial $C^1$ function $h$ with $h'\ge0$ is $k$-convex if
$r^{n-k}(h')^k$ is nondecreasing. If this flux is locally absolutely
continuous, \eqref{hessian-radial-identity} gives its Hessian measure
away from $0$.
\end{lemma}
\begin{proof}
For a smooth positive flux, monotonicity gives
$h''+(n-k)h'/(kr)\ge0$. For $j<k$, the $j$-th radial admissibility
bracket exceeds this by $n(k-j)h'/(jkr)>0$, so all
$\sigma_j(D^2h)$ are nonnegative. We approximate a general flux on compact
radial intervals by smooth positive nondecreasing fluxes and integrate
to recover $h$. Closure of $\Phi^k$ and weak continuity of Hessian
measures \cite{TW2} give both assertions, including intervals where
$h'=0$.
\end{proof}

With $F(r)=\sum_i\mu_i(B_r)$, the radial comparison function
\eqref{intro-source-reconstruction} takes the form
\begin{equation}\label{radial-z}
 z(r)=K_0+\int_r^R
 \left(\frac{F(s)}{\omega d_0s^{a\kappa}}\right)^{1/a}\frac{\dd s}{s},
\end{equation}
where $R>0$ is fixed and small, and $K_0>0$ is chosen sufficiently large.

\begin{lemma}\label{radial-comparison}
Assume
\begin{equation}\label{total-source}
 f_i=\mathcal Q[w_i]\ge0,\qquad
 \sum_i f_i\le AS^q,\qquad |x|^\kappa S(x)\to0.
\end{equation}
Then the reconstruction \eqref{radial-z} satisfies, for small $|x|=r$,
\begin{equation}\label{radial-comparison-bound}
 S(x)\le C_1z(r),\qquad z(r)\le K_1+C_2S(x),
\end{equation}
and
\begin{equation}\label{z-equation}
 0\le\mathcal Q[z](r)=\frac{F'(r)}{\omega r^{n-1}}
 \le C_3z(r)^q\quad\text{a.e.}
\end{equation}
\end{lemma}

Here $K_1<\infty$, and $C_1,C_2$ depend only on the operator parameters,
$N$, and $A$. The flux is locally absolutely continuous,
$r^\kappa z(r)\to0$, and $-z\in\Phi^k$ in the Hessian case.
The comparison and monotonicity of $z$ imply that $S$ is either
bounded near $0$ or tends uniformly to infinity.
\begin{proof}
A positive mass at $0$ would give $w_i(x)\ge c|x|^{-\kappa}$ by the
lower estimate in \eqref{common-potential}, contrary to
$|x|^\kappa S(x)\to0$. Thus all point masses vanish and $F(0)=0$.
Local finiteness gives $f_i\in L^1$ on every smaller full ball,
hence absolute continuity of $F$ there.
By \eqref{radial-z},
\[
 J_z(r)=\frac{F(r)}{\omega d_0}.
\]
Consequently \eqref{common-radial-identity} gives
\[
 \mathcal Q[z](r)=\frac{F'(r)}{\omega r^{n-1}}
 =\av_{\Sph}\sum_i f_i(r\theta)\quad\text{a.e.}
\]
Lemma~\ref{radial-admissibility} gives $-z\in\Phi^k$ in the
Hessian case, since $F$ is nonnegative and nondecreasing.

For $|x|=r$ and $s\ge r$, the inclusion $B_s(0)\subset B_{2s}(x)$,
the lower potential estimate, and
$(\sum_iX_i)^{1/a}\le C_{a,N}\sum_iX_i^{1/a}$ give
\[
 \int_r^{R/8}\left(\frac{F(s)}{s^{a\kappa}}\right)^{1/a}
 \frac{\dd s}{s}\le CS(x).
\]
The outer integral is finite and independent of $r$, proving the
second inequality in \eqref{radial-comparison-bound}. Multiplication
by $r^\kappa$ gives $r^\kappa z(r)\to0$.

For the opposite comparison, apply the upper estimate in
\eqref{common-potential} at $x$ with a fixed small radius. Bound each
infimum by its value at a fixed point common to the relevant balls,
and split the potential at $t=r/2$.
For $t\ge r/2$, the inclusion $B_t(x)\subset B_{3t}(0)$ bounds the
singular part by $Cz(r)$; outer scales contribute only a fixed constant.
For $t<r/2$, \eqref{total-source} gives
\[
 \mu_i(B_t(x))\le C t^n
 \left(\sup_{r/2\le|y|\le3r/2}S(y)\right)^q.
\]
Since $n-a\kappa=\ell$, integration of the small-scale potential yields
\begin{equation}\label{S-local-absorb}
 S(x)\le Cz(r)+Cr^{\ell/a}
 \left(\sup_{r/2\le|y|\le3r/2}S(y)\right)^{q/a},
\end{equation}
after enlarging $K_0$ to absorb fixed contributions.

For the doubling estimate on $z$, use monotonicity of $F$ for $2r<R$:
\[
 z(r/2)-z(r)\le cF(r)^{1/a}r^{-\kappa},\qquad
 z(r)-z(2r)\ge c'F(r)^{1/a}r^{-\kappa}.
\]
Thus $z(r/2)\le Dz(r)$ with a structural constant $D$.
Let $r_j=r_*2^{-j}$,
$M_j=\sup_{r_j/2\le|x|\le r_j}S(x)$,
$h_j=M_j/z(r_j)$, and $H_j=\max(M_{j-1},M_j,M_{j+1})$.
These three neighboring dyadic annuli cover the local annulus in
\eqref{S-local-absorb}, with comparable values of $z$.
After division by $z(r_j)$, the nonlinear term is bounded by
\[
 C\frac{r_j^{\ell/a}H_j^{q/a}}{z(r_j)}
 =C\frac{H_j}{z(r_j)}
       (r_j^\kappa H_j)^{(q-a)/a}.
\]
Here $\ell=\kappa(q-a)$ and $r_j^\kappa H_j\to0$.
Hence
\[
 h_j\le C+o(1)\max(h_{j-1},h_j,h_{j+1}).
\]
Also $h_j\le K2^{\kappa j}$ because $z\ge K_0$ and
$S=O(|x|^{-\kappa})$. Lemma~\ref{dyadic} proves the first comparison
in \eqref{radial-comparison-bound}. Spherical averaging of
\eqref{total-source} then proves \eqref{z-equation}.

Finally, $z$ is nonincreasing in $r$. If $S$ is unbounded
near $0$, the upper comparison gives $z(r)\to\infty$, and the
opposite comparison gives $S(x)\to\infty$ uniformly.
\end{proof}

\begin{remark}\label{linear-consistency}
When $a=1$, one has $\mathcal Q=-\Delta$, $d_0=1$, and
$\kappa=n-2$. Let $\bar S(r)=\av_{\Sph}S(r\theta)$. The
Gauss--Green identity gives
$F(r)=-\omega r^{n-1}\bar S'(r)$, so \eqref{radial-z} yields exactly
\[
 z(r)=\bar S(r)+K_0-\bar S(R).
\]
Thus source-mass reconstruction recovers the ordinary spherical average
up to the chosen additive constant. For $a\ne1$, it preserves the
averaged source through \eqref{z-equation}, rather than commuting the
nonlinear operator with averaging of the unknown.
\end{remark}

\subsection{The radial differential inequality}

\begin{lemma}\label{radial-ode}
Fix either operator in~\eqref{parameters} and its associated parameters.
Let $z>0$ be radial and $C^1$
away from $0$, with locally absolutely continuous flux
\eqref{radial-flux}; in the Hessian case assume $-z\in\Phi^k$.
If $0\le\mathcal Q[z]\le Dz^q$, $D>0$,
$z(r)\to\infty$, and $r^\kappa z(r)\to0$, then
\begin{equation}\label{radial-gradient}
 \frac{-rz'(r)}{z(r)}\to\kappa,\qquad
 \frac{z(r\lambda)}{z(r)}\to\lambda^{-\kappa},
\end{equation}
locally uniformly for $\lambda\in(0,\infty)$, and
\begin{equation}\label{radial-scalar-lower}
 \liminf_{r\to0}r^\kappa L(r)^\beta z(r)
 \ge D^{-\beta}C_0.
\end{equation}
\end{lemma}
\begin{proof}
First $-z'>0$ at all small radii. For $-\Delta_p$ this
follows from monotonicity of the signed flux: nonpositive flux
arbitrarily close to $0$ would contradict $z\to\infty$.
For the Hessian operator, apply the same argument to the nondecreasing
subharmonic flux $r^{n-1}(-z')$ of the $k$-convex function $-z$.
The flux \eqref{radial-flux} is therefore positive near $0$ and is
nondecreasing in $r$ by \eqref{common-radial-identity}.
We set $t=\log(1/r)$,
$Z(t)=e^{-\kappa t}z(e^{-t})$ and
$J(t)=J_z(e^{-t})$. In both cases,
$J=r^{a\kappa}(-rz')^a$ on this interval.
A positive limiting flux would give a positive lower limit of
$r^\kappa z(r)$ by integration, so $J(t)\to0$.
With $D_1=D/d_0$, \eqref{common-radial-identity} and $\kappa q=n$ give
\begin{equation}\label{radial-ODE}
 Z'=J^{1/a}-\kappa Z,\qquad 0\le-J'\le D_1 Z^q,
 \qquad Z,J\to0.
\end{equation}
The variation-of-constants formula and monotonicity of $J$ imply
\[
 Z(t)\ge\frac{1-e^{-\kappa(t-T)}}\kappa J(t)^{1/a}.
\]
Consequently $R(t)=J(t)/Z(t)^a$ has $\limsup R\le\kappa^a$. Almost everywhere,
\begin{equation}\label{Rineq}
 R'\ge-D_1 Z^{q-a}+aR(\kappa-R^{1/a}).
\end{equation}
Fix $0<\eta<\kappa^a$. For large $t$, \eqref{Rineq} gives a positive
lower bound for $R'$ whenever $R$ lies in a small closed interval about
$\eta$. Absolute continuity prevents a downward crossing of $\eta$.
If instead $R<\eta$ eventually, then
$Z'\le-\delta Z$ with $\delta=\kappa-\eta^{1/a}>0$.
Using $J(\infty)=0$ gives
\[
 J(t)\le D_1\int_t^\infty Z(s)^q\dd s
 \le\frac{D_1}{q\delta}Z(t)^q.
\]
It follows that
\[
 0\le(\log Z+\kappa t)'\le C Z^{q/a-1},
\]
whose right side is integrable. Thus $e^{\kappa t}Z(t)=z(e^{-t})$
stays bounded, a contradiction. Therefore $R\ge\eta$ eventually.
Letting $\eta\to\kappa^a$ from below proves $R\to\kappa^a$ and $Z'/Z\to0$,
which gives \eqref{radial-gradient}.

Finally, for fixed $\eta<\kappa^a$ and large $t$,
$-J'\le D_1\eta^{-q/a}J^m$. Integration gives
\[
 \liminf_{t\to\infty}t^\beta J(t)^{1/a}
 \ge[D_1\eta^{-q/a}(m-1)]^{-\beta}.
\]
We combine this with $Z/J^{1/a}\to1/\kappa$ and let
$\eta\to\kappa^a$ from below to obtain \eqref{radial-scalar-lower}, using
$D_1=D/d_0$.
The fixed-shift ratio $J(t+s)/J(t)\to1$ follows from
$Z'/Z\to0$ and $Z\sim\kappa^{-1}J^{1/a}$.
\end{proof}

\subsection{Asymptotic formulas in terms of the individual masses}

On the unbounded branch of Lemma~\ref{radial-comparison}, set
$t=L(r)$, $Z(t)=r^\kappa z(r)$, and
\[
 P_i(t)=\frac{\mu_i(B_{e^{-t}})}{\omega d_0},\qquad
 J(t)=\sum_iP_i(t)=\frac{F(r)}{\omega d_0}.
\]
Lemma~\ref{radial-ode} and its proof give
\begin{equation}\label{hessian-slow-mass}
 Z(t)\sim\kappa^{-1}J(t)^{1/a},\qquad
 \frac{J(t+s)}{J(t)}\to1,\qquad
 \frac{z(r\lambda)}{z(r)}\to\lambda^{-\kappa},
\end{equation}
for fixed $s\in\R$ and $\lambda>0$; the last convergence is locally
uniform in $\lambda\in(0,\infty)$.
We identify each component on the total-mass scale, even when its
mass share tends to zero.

\begin{proposition}\label{mass-approx}\label{hessian-mass}
Assume \eqref{total-source} and $S(x)\to\infty$ as $x\to0$.
Then, for $1\le i\le N$, uniformly in $\theta\in\Sph$,
\begin{equation}\label{mass-expansion}
 e^{-\kappa t}w_i(e^{-t}\theta)
 =\kappa^{-1}P_i(t)^{1/a}+o(J(t)^{1/a}).
\end{equation}
For $\mathcal Q=-\Delta_p$, the corresponding gradient formula is
\begin{equation}\label{mass-gradient-expansion}
 r^{\kappa+1}\nabla w_i(r\theta)
 =-P_i(L(r))^{1/a}\theta+o(J(L(r))^{1/a}),
 \qquad r\to0,
\end{equation}
again uniformly in direction.
\end{proposition}

In particular,
\begin{equation}\label{norm-z}
 \frac{(\sum_iw_i(r\theta)^a)^{1/a}}{z(r)}\to1,
 \qquad
 \frac{\max_{|x|=r}S(x)}{\min_{|x|=r}S(x)}\to1.
\end{equation}
In the $p$-Laplace case,
\begin{equation}\label{total-gradient-limit}
 \frac{r\nabla S(r\theta)}{S(r\theta)}\to-\kappa\theta
\end{equation}
uniformly in $\theta$. Hence $\partial_rS<0$ at all sufficiently small
radii and $|\nabla S|\sim\kappa S/r$. If $P_i/J$ is bounded away from
zero, \eqref{mass-expansion} also gives the individual relative formula
$r^\kappa w_i(r\theta)\sim\kappa^{-1}P_i(L(r))^{1/a}$.
\begin{proof}
Lemma~\ref{radial-comparison} gives $z\to\infty$ and
$cz(|x|)\le S(x)\le Cz(|x|)$ at small radii. For an arbitrary
decreasing sequence $r_j\to0$, set
\[
 W_{i,j}(y)=\frac{w_i(r_jy)}{z(r_j)}.
\]
The comparison and \eqref{radial-gradient} bound these functions on
every compact annulus. Homogeneity of the relevant operator gives
\[
 0\le\mathcal Q[W_{i,j}](y)
 \le Cr_j^\ell z(r_j)^{q-a}
 \left(\frac{S(r_jy)}{z(r_j)}\right)^q\to0,
\]
since $r_j^\ell z(r_j)^{q-a}=(r_j^\kappa z(r_j))^{q-a}\to0$.
We identify the limits and their masses for each operator.

For $-\Delta_p$, local regularity yields, after diagonal extraction,
$C^1$ convergence on compact annuli to nonnegative $p$-harmonic
functions $W_i$ on $\R^n\setminus\{0\}$. Classification
\cite[Theorem~5.3]{MZ} gives
$W_i=A_i|y|^{-\kappa}+B_i$, with $A_i,B_i\ge0$.
The same comparison constant applies at every fixed $|y|$, and
\eqref{radial-gradient} gives
$\sum_iW_i(y)\le C|y|^{-\kappa}$ on the whole punctured space.
Thus $B_i=0$. The flux identity
\[
 \mu_i(B_r)=-\int_{\partial B_r}
 |\nabla w_i|^{p-2}\partial_rw_i\dd S
\]
first follows for almost every radius by radial cutoff testing.
Both sides are continuous for $r>0$, since the source density is locally
bounded and the flux vector continuous away from $0$. The identity holds
for every radius, and $C^1$ convergence gives
\begin{equation}\label{flux-identification}
 \frac{\mu_i(B_{r_j})}
 {\omega d_0\,r_j^{a\kappa}z(r_j)^a}
 \to\kappa^aA_i^a.
\end{equation}
Here $d_0=1$ and $a\kappa=n-p$.

For the Hessian operator, apply \eqref{hessian-oscillation} to
$H_{i,j}=-W_{i,j}$. A diagonal subsequence converges locally uniformly
away from $0$ to a nonpositive $k$-convex function $H_i$. To retain the measure
at $0$, use monotonicity of $F$ in \eqref{radial-z}: for $0<\rho\le1$,
\[
 z(r\rho)\le z(r)+CF(r)^{1/a}r^{-\kappa}\rho^{-\kappa}.
\]
By \eqref{hessian-slow-mass},
$F(r)^{1/a}r^{-\kappa}/z(r)$ stays bounded. Therefore
$|H_{i,j}(y)|\le C|y|^{-\kappa}$ for $0<|y|\le1$.
Since $\kappa<n$ for the Hessian parameters, this gives local
$L^1$ convergence across $0$. Weak continuity identifies the limiting
Hessian measure, which vanishes away from $0$. Also,
$|H_i(y)|\le C|y|^{-\kappa}$ on the whole punctured space by the
comparison and \eqref{radial-gradient}. Hence the measure has the form
$M_i\delta_0$, and $H_i\to0$ at infinity. The uniqueness statement
\eqref{hessian-fundamental-measure} gives
$H_i=-A_i|y|^{-\kappa}$, with $A_i\ge0$.
The boundary of $B_1$ carries no limiting measure, so weak convergence
of the scaled measures yields \eqref{flux-identification} as well,
now with $a=k$, $a\kappa=n-2k$, and the stated $d_0$.

For either operator the definition of $z$ gives
\[
 \frac{F(r)}{\omega d_0\,r^{a\kappa}z(r)^a}
 =\left(\frac{-rz'(r)}{z(r)}\right)^a\to\kappa^a.
\]
Summing \eqref{flux-identification} therefore shows
$\sum_iA_i^a=1$, and then
\[
 \frac{P_i(t_j)}{J(t_j)}\to A_i^a,
 \qquad
 r_j^\kappa z(r_j)
 =\kappa^{-1}J(t_j)^{1/a}(1+o(1)),\qquad t_j=L(r_j).
\]
On $|y|=1$ the convergence to $A_i$ is uniform. Multiplying by
$r_j^\kappa z(r_j)$ gives \eqref{mass-expansion} along the extracted
subsequence. The coefficients $A_i$ may depend on the subsequence, but the leading
terms use the actual masses at each radius. Any sequence violating
uniform decay of the normalized remainder would have a subsequence
covered by this argument. Thus \eqref{mass-expansion} holds for all radii.
In the $p$-Laplace case the convergence is in $C^1$ on compact annuli,
so $\nabla W_{i,j}(\theta)\to-\kappa A_i\theta$ uniformly. Multiplying
by $r_j^\kappa z(r_j)$ and using the same mass identification proves
\eqref{mass-gradient-expansion} by the same subsequence argument.

We finally divide \eqref{mass-expansion} by
$r^\kappa z(r)\sim\kappa^{-1}J(t)^{1/a}$.
The continuous function $(x_i)\mapsto(\sum_i x_i^a)^{1/a}$ proves the
first limit in \eqref{norm-z}. The sum
$\sum_i(P_i/J)^{1/a}$ has a positive lower bound on the compact simplex
$\{P_i/J\ge0,\ \sum_iP_i/J=1\}$, so the uniform additive errors also
give the spherical max--min ratio in \eqref{norm-z}. Summing
\eqref{mass-gradient-expansion} and dividing by the corresponding sum
in \eqref{mass-expansion} proves \eqref{total-gradient-limit}.
\end{proof}

\section{The sharp joint lower bound}\label{sec:sharp}

We compare the mass inequalities with an increasing two-dimensional
flow. Its blow-up time determines the optimal joint constant.

\subsection{The associated two-dimensional system}
For $X,Y\ge0$, $(X,Y)\ne(0,0)$, consider the increasing system
\begin{equation}\label{reverse-flow}
 \dot X=Y^m,\qquad\dot Y=X^m.
\end{equation}
Denote its common finite blow-up time by $T_m(X,Y)$. It is symmetric and
homogeneous of degree $1-m$. If $X\ge Y$, direct separation using the
invariant $X^{m+1}-Y^{m+1}$ gives
\begin{equation}\label{Tformula}
 T_m(X,Y)=X^{1-m}I_m(Y/X).
\end{equation}
This also proves finiteness at the axes.

\begin{lemma}\label{time-properties}
On $[0,\infty)^2\setminus\{(0,0)\}$, the function $T_m$ is $C^1$,
decreasing in each variable, and satisfies
\begin{equation}\label{time-pde}
 (T_m)_XY^m+(T_m)_YX^m=-1,\qquad
 |\nabla T_m(X,Y)|\le C_m(X+Y)^{-m}.
\end{equation}
\end{lemma}
\begin{proof}
Comparison in \eqref{reverse-flow} gives monotonicity, strict in the
interior. The remaining lifetime decreases at speed $1$, giving the
first identity. Formula \eqref{Tformula} gives interior differentiability,
including the diagonal by symmetry. Differentiation under the
integral gives
\[
 (T_m)_Y(X,Y)=-mY^m\int_X^\infty
 (s^{m+1}-X^{m+1}+Y^{m+1})^{-(2m+1)/(m+1)}\dd s.
\]
We justify the axis limit by splitting the integral at $X+\delta$ and using
$s^{m+1}-X^{m+1}=(m+1)X^m(s-X)+O((s-X)^2)$ near $X$.  In the near part the
substitution $s-X=Y^{m+1}z$ gives
\[
 -mY^m\int_X^{X+\delta}
 (s^{m+1}-X^{m+1}+Y^{m+1})^{-(2m+1)/(m+1)}\dd s
 \to -X^{-m},
\]
whereas the far part tends to zero. Thus $(T_m)_Y(X,0)=-X^{-m}$;
Euler's identity gives $(T_m)_X(X,0)=(1-m)I_m(0)X^{-m}$.
This decomposition and the formula away from the axis prove continuity
of both derivatives. Compactness of $\{X+Y=1,\ X,Y\ge0\}$ and homogeneity
give the gradient estimate. The beta-integral
evaluation used in \eqref{beta-constant} follows by substituting
$z=y^{-(m+1)}$ in \eqref{Idef}.
\end{proof}

\subsection{Application to both source measures}

\begin{proof}[Proof of the bounds in Theorem~\ref{mainjoint}]
On the zero-mass singular branch, let $\mu_1,\mu_2$ be the source
measures and set
\[
 P(t)=\frac{\mu_1(B_{e^{-t}})}{\omega d_0},\qquad
 Q(t)=\frac{\mu_2(B_{e^{-t}})}{\omega d_0},\qquad J=P+Q,
 \quad S=w_1+w_2.
\]
Proposition~\ref{mass-approx} gives, for both operators,
\[
 e^{-\kappa t}(w_1,w_2)(e^{-t}\theta)
 =\kappa^{-1}(P^{1/a},Q^{1/a})+o(J^{1/a}).
\]
The errors are uniform in direction, including at limiting mass ratios
on the axes; the required source bound follows from
$w_1^q+w_2^q\le\max(1,2^{1-q})S^q$.
Because $\kappa q=n$, differentiation of the masses gives almost everywhere
\begin{equation}\label{mass-control}
 \begin{split}
 0\le-P'(t)&\le d_0^{-1}\kappa^{-q}Q(t)^m+o(J(t)^m),\\
 0\le-Q'(t)&\le d_0^{-1}\kappa^{-q}P(t)^m+o(J(t)^m).
 \end{split}
\end{equation}
Indeed $-P'=d_0^{-1}e^{-nt}\av_{\Sph}f_1(e^{-t}\theta)$ and
\eqref{mass-expansion} can be raised to the positive power $q$ with an
error $o(J^m)$, even when one component is small.

Define $\Theta(t)=d_0\kappa^q T_m(P(t),Q(t))$. The total mass $J(t)$ is
positive at all sufficiently large times. Indeed, if it vanished at some
such time, the total-mass formula for $z$ would make it constant
at smaller radii,
contrary to $z\to\infty$. Thus $T_m(P,Q)$ is defined, including at any point
where one individual mass is zero. Lemma~\ref{time-properties} supplies
$C^1$ regularity up to those axes, so the chain rule applies to the
absolutely continuous functions $P,Q$.
Write the two errors in \eqref{mass-control} as
$\varepsilon_i(t)J(t)^m$, with $\varepsilon_i(t)\to0$.  Since
$(T_m)_X,(T_m)_Y\le0$ and $P',Q'\le0$, the chain rule and
\eqref{time-pde} give
\[
 \Theta'(t)
 \le1+CJ(t)^m|\nabla T_m(P(t),Q(t))|
       (|\varepsilon_1(t)|+|\varepsilon_2(t)|)
 =1+o(1).
\]
Here $|\nabla T_m|\le CJ^{-m}$. Integration from a fixed large time gives
$\Theta(t)\le(1+o(1))t$.
Set $s(t)=\min(P,Q)/\max(P,Q)$. By \eqref{Tformula} and the mass
approximation,
\begin{equation}\label{refined-factorization}
 t^\beta e^{-\kappa t}S(e^{-t}\theta)
 =(1+o(1))\left(\frac{t}{\Theta(t)}\right)^\beta
 d_0^\beta\Phi_{a,\kappa}(s(t)),
\end{equation}
uniformly in $\theta$. Since $t/\Theta(t)\ge1-o(1)$ and
$d_0^\beta\Phi_{a,\kappa}\ge\CQ$, this proves \eqref{joint-lower}.

Since $\max(P,Q)\ge J/2$, the same expansion gives
\[
 t^\beta e^{-\kappa t}\max\{w_1,w_2\}(e^{-t}\theta)
 =(1+o(1))\left(\frac{t}{\Theta(t)}\right)^\beta
 d_0^\beta\kappa^{q\beta-1}I_m(s(t))^\beta.
\]
The function $I_m$ is decreasing and $I_m(1)=1/(m-1)$, with
$d_0^\beta\kappa^{q\beta-1}I_m(1)^\beta=C_0$.
This proves \eqref{joint-max-lower}.

For \eqref{joint-min-upper}, set
$N_i(x)=|x|^\kappa L(|x|)^\beta w_i(x)$ and suppose that
$M=\limsup_{x\to0}\min(N_1,N_2)<C_0$.
The preceding bound gives $\max(N_1,N_2)>\min(N_1,N_2)$ in a small
punctured ball. By continuity and connectedness, after interchanging
components, $w_1>w_2$ there. Thus $\limsup N_2=M$ and $w_1\ge S/2$.
The mass expansion then bounds $P/J$ away from zero and yields
\[
 e^{-\kappa t}w_1(e^{-t}\theta)
 =\kappa^{-1}P(t)^{1/a}(1+o(1))
\]
uniformly in $\theta$. For every $\varepsilon>0$, the bound on $N_2$
and the source inequality give, for a.e. sufficiently large $t$,
\[
 -P'(t)\le d_0^{-1}e^{-nt}\av_{\Sph}w_2(e^{-t}\theta)^q
 \le d_0^{-1}(M+\varepsilon)^q t^{-q\beta}.
\]
Integrating from infinity, using $P(\infty)=0$ and
$q\beta-1=a\beta$, and then letting $\varepsilon\to0$, we obtain
\[
 \limsup_{x\to0}N_1(x)
 \le\kappa^{-1}\left(\frac{M^q}{d_0a\beta}\right)^{1/a}
 =C_0\left(\frac{M}{C_0}\right)^m.
\]
If $M>0$, the right-hand side is strictly less than $M$, contrary to
$N_1>N_2$ and $\limsup N_2=M$. If $M=0$, both normalized components
tend to zero, contrary to \eqref{joint-max-lower}.
This proves \eqref{joint-min-upper}.
The diagonal example $w_1=w_2=C_0r^{-\kappa}L(r)^{-\beta}$
verified in the proof of Corollary~\ref{hessian-scalar} attains
both bounds on a sufficiently small punctured ball.
\end{proof}

\begin{proof}[The trichotomy and asymptotic radiality in
Theorem~\ref{mainjoint}]
For either operator, Lemma~\ref{fundamental-vector} gives
\eqref{component-fundamental-limits}. If $b_1+b_2>0$, this proves
alternative~(ii) and asymptotic radiality of $S=w_1+w_2$.

If $b_1=b_2=0$, the same lemma gives $|x|^\kappa S(x)\to0$.
By Lemma~\ref{radial-comparison}, either $S$ is bounded near $0$
or $S(x)\to\infty$ uniformly.
In the bounded case the densities are bounded and the extended measures
have no mass at $0$. The extensions are removable: for $-\Delta_p$
use removability across a point of zero $p$-capacity and local regularity
\cite{HKM,Tolksdorf}; for the Hessian operator use
\eqref{hessian-decomposition} and \eqref{hessian-oscillation}.
The strong minimum principle for the nonnegative $p$-superharmonic
components, or the strong maximum principle for the negative subharmonic
Hessian components, makes the extended $w_i$ strictly positive
at $0$. Continuity then gives the spherical ratio limit for $S$.
This is alternative~(i).

On the unbounded branch, both components tend to infinity. Indeed,
suppose that $\liminf_{x\to0}w_1(x)<\infty$, so the canonical extension
has $w_1(0)<\infty$. The lower estimate in \eqref{common-potential}
and monotonicity of the ball masses give, for some small $R_0>0$,
\[
 \int_0^{R_0}\mu_1(B_s)^{1/a}s^{-\kappa-1}\dd s<\infty,
 \qquad \mu_1(B_s)\le Cs^{a\kappa}\quad(0<s<R_0/2).
\]
For $x\in B_r$ and small $r$, apply the upper estimate in
\eqref{common-potential} with radius $2r$. Its infimum ball contains
$0$, and the potential uses only the restriction of $\mu_1$ to $B_{5r}$.
Thus \eqref{weak-wolff} and $\kappa q=n$ imply
\[
 |\{x\in B_r:w_1(x)>\lambda\}|
 \le C\mu_1(B_{5r})^{q/a}\lambda^{-q}
 \le Cr^n\lambda^{-q}\qquad(\lambda\ge\lambda_0),
\]
where $\lambda_0$ is independent of $r$.
On $A_r=B_r\setminus B_{r/2}$, Lemma~\ref{fundamental-vector}
gives $w_1\le Cr^{-\kappa}$. Integrating the distribution bound yields
\[
 \int_{A_r}w_1^q\dd x
 \le Cr^n\left(1+\int_{\lambda_0}^{Cr^{-\kappa}}
                  \frac{\dd\lambda}{\lambda}\right)
 \le Cr^n(1+L(r)).
\]
Since $\mu_2(\{0\})=0$ and $\mathcal Q[w_2]\le w_1^q$,
summing over dyadic annuli gives $\mu_2(B_r)\le Cr^n(1+L(r))$.
Consequently, using $n/a-\kappa=\ell/a>0$,
\[
 \int_0^{R_0}\mu_2(B_s)^{1/a}s^{-\kappa-1}\dd s
 \le C\int_0^{R_0}s^{\ell/a-1}(1+L(s))^{1/a}\dd s<\infty.
\]
Both centered potentials are finite, so the total-mass reconstruction
$z$ in \eqref{radial-z} is bounded. Lemma~\ref{radial-comparison}
then bounds $S$, a contradiction. Interchanging the components proves
$w_1,w_2\to\infty$.

The bounds proved above and Proposition~\ref{mass-approx} complete
alternative~(iii) and asymptotic radiality. The same fundamental-scale
and radial comparison arguments give the scalar classification in
Corollary~\ref{hessian-scalar}, whose sharp constant is proved below.
Section~\ref{sec:construction} proves sharpness of \eqref{joint-lower}.
\end{proof}

\begin{corollary}\label{finite-cluster-region}
Set
\begin{equation}\label{cluster-region}
 \mathcal D_{\mathcal Q}
 =\left\{(A_1,A_2)\in[0,\infty)^2\setminus\{(0,0)\}:
 d_0\kappa^aT_m(A_1^a,A_2^a)\le1\right\}.
\end{equation}
For every solution in alternative~\textup{(iii)} of
Theorem~\ref{mainjoint}, every finite cluster point of
\[
 (N_1,N_2)(t,\theta)
 =t^\beta e^{-\kappa t}(w_1,w_2)(e^{-t}\theta),
 \qquad t\to\infty,
\]
belongs to $\mathcal D_{\mathcal Q}$, even when the direction varies
with $t$. Conversely, there is a smooth radial solution in that alternative
whose set of finite cluster points is exactly $\mathcal D_{\mathcal Q}$.
\end{corollary}
\begin{proof}
Use $P,Q$ and $\Theta$ from the proof of Theorem~\ref{mainjoint},
and set $\widehat P=t^{a\beta}P$, $\widehat Q=t^{a\beta}Q$.
The mass approximation gives
\[
 (N_1,N_2)=\kappa^{-1}(\widehat P^{1/a},\widehat Q^{1/a})
       +o((\widehat P+\widehat Q)^{1/a}).
\]
After division by $(\widehat P+\widehat Q)^{1/a}$, the leading
vectors range over a compact subset of the nonnegative quadrant
that does not contain the origin. Continuity and positivity of $T_m$
on this set, followed by its homogeneity of degree $1-m$, give
\begin{equation}\label{normalized-time-constraint}
 d_0\kappa^aT_m(N_1^a,N_2^a)
 =(1+o(1))\frac{\Theta(t)}{t}\le1+o(1),
\end{equation}
uniformly in $\theta$. Here we used
$a\beta(m-1)=1$ and $am=q$. The lower bound
\eqref{joint-max-lower} excludes the zero vector as a cluster point,
so continuity proves the required inclusion.

In the half-plane $A_1\ge A_2$, equality in the constraint
\eqref{cluster-region} is parameterized by
\begin{equation}\label{cluster-boundary}
 (A_1(s),A_2(s))
 =\bigl(d_0\kappa^aI_m(s)\bigr)^\beta(1,s^{1/a}),
 \qquad 0\le s\le1.
\end{equation}
Every point of the region in this half-plane is
$\lambda^{-\beta}(A_1(s),A_2(s))$ for some $0<\lambda\le1$.
Section~\ref{sec:finite-realization} realizes each choice; interchanging
components gives the other half. In particular,
\[
 \CQ=\min_{(A_1,A_2)\in\mathcal D_{\mathcal Q}}(A_1+A_2).
\]
\end{proof}

\begin{proposition}\label{boundary-arc}
Let $(w_1,w_2)$ be in alternative~\textup{(iii)} of
Theorem~\ref{mainjoint}, and let $N=(N_1,N_2)$ be as in
Corollary~\ref{finite-cluster-region}. Suppose
$N(t_j,\theta_j)\to A$ for some $t_j\to\infty$, where
$c_*T_m(A_1^a,A_2^a)=1$ and $c_*=d_0\kappa^a$.
Then, locally uniformly for $0<\tau\le1$ and uniformly in $\theta$,
\begin{equation}\label{boundary-trajectory}
 \begin{gathered}
 N(\tau t_j,\theta)\longrightarrow
 \tau^\beta\bigl(U_1(\tau)^{1/a},U_2(\tau)^{1/a}\bigr),\\
 c_*U_1'=-U_2^m,\quad c_*U_2'=-U_1^m,\quad
 U(1)=(A_1^a,A_2^a).
 \end{gathered}
\end{equation}
In particular, the finite cluster set of this same solution contains
the entire arc of the equality curve in~\eqref{cluster-region}
joining $A$ to $(C_0,C_0)$.
\end{proposition}
\begin{proof}
Use $P,Q,\Theta$ from the proof of Theorem~\ref{mainjoint}.
By~\eqref{normalized-time-constraint}, $\Theta(t_j)/t_j\to1$.
Since $0\le\Theta'\le1+o(1)$ and $\Theta(t)\le(1+o(1))t$,
\[
 \frac{\Theta(\tau t_j)}{t_j}\longrightarrow\tau
 \quad\text{uniformly on }[\delta,1],\qquad \delta>0.
\]
Set $U_j(\tau)=\kappa^{-a}t_j^{a\beta}(P,Q)(\tau t_j)$.
Homogeneity gives $c_*T_m(U_j(\tau))=\Theta(\tau t_j)/t_j$.
Since $T_m(X,Y)$ is comparable to $(X+Y)^{1-m}$, the sums of
the components of $U_j$ are bounded above and away from zero on
$[\delta,1]$. Equation~\eqref{mass-control} gives equicontinuity.
Every subsequential limit $U$ therefore satisfies
\[
 \begin{gathered}
 c_*T_m(U(\tau))=\tau,\qquad U(1)=(A_1^a,A_2^a),\\
 0\le-c_*U_1'\le U_2^m,\qquad
 0\le-c_*U_2'\le U_1^m
 \quad\text{a.e.}
 \end{gathered}
\]
Both partial derivatives of $T_m$ are strictly negative, including
on the axes, by the formulas in Lemma~\ref{time-properties}.
Differentiating the identity and using~\eqref{time-pde} forces equality
in both differential inequalities. Uniqueness for this locally
Lipschitz system identifies every subsequential limit.
The mass approximation now proves~\eqref{boundary-trajectory} for
the original sequence.

The limiting normalized vector lies on the equality curve.
After interchanging components, assume $A_1\ge A_2$.
The invariant $U_1^{m+1}-U_2^{m+1}$ and
$c_*T_m(U)=\tau$ show that $U_2/U_1$ increases from
$A_2^a/A_1^a$ to $1$ as $\tau\to0$, unless $A$ is diagonal.
Thus the limiting trace tends to $(C_0,C_0)$ and covers
the stated arc except possibly this endpoint. Closedness of the finite
cluster set includes the endpoint as well.
\end{proof}

\begin{remark}\label{mass-ratio-refinement}
Suppose $t_j\to\infty$ and
$t_j^\beta e^{-\kappa t_j}S(e^{-t_j}\theta_j)\to\CQ$.
Both lower bounds in \eqref{refined-factorization} are then attained
asymptotically:
\[
 \operatorname{dist}\!\left(s(t_j),
 \operatorname*{argmin}_{[0,1]}\Phi_{a,\kappa}\right)\to0,
 \qquad \Theta(t_j)/t_j\to1.
\]
For $a\ge1$, Lemma~\ref{minimum-axis} gives $s(t_j)\to0$, and
\eqref{mass-expansion} yields, uniformly in $\theta$,
\begin{equation}\label{sharp-component-splitting}
 t_j^\beta e^{-\kappa t_j}
 \bigl(\max\{w_1,w_2\},\min\{w_1,w_2\}\bigr)(e^{-t_j}\theta)
 \to(\CQ,0).
\end{equation}
Indeed $J(t_j)\sim\max(P(t_j),Q(t_j))$, so the error remains negligible
on the larger mass scale; no relative estimate for the smaller component
is needed.
For $a<1$, Lemma~\ref{minimum-subquadratic} gives a unique minimizing
ratio $s_*>0$, equal to $1$ precisely when $a\le(3m+1)/(m+1)^2$.
Thus $s(t_j)\to s_*$, and both normalized components stay bounded
away from zero along these sequences, uniformly in $\theta$.
Finally, \eqref{joint-min-upper} gives
$\limsup_{x\to0}|x|^\kappa L(|x|)^\beta S(x)\ge2C_0$.
If $\CQ<2C_0$, the normalized sum of any solution attaining
\eqref{joint-lower} therefore cannot converge.
\end{remark}

\subsection{Evaluation of the sharp constant}

\begin{lemma}\label{minimum-axis}
For $a\ge1$, the function $(1+s^{1/a})I_m(s)^{1/[a(m-1)]}$ is strictly
increasing on $[0,1]$.
\end{lemma}
\begin{proof}
For $a=1$, write $N=X+Y$ and $D=X^{m+1}-Y^{m+1}\ge0$.
At fixed $N$, increasing $D$ increases the imbalance and, since $m>1$,
strictly increases $X^m+Y^m$. Since $D$ is constant along the flow,
the blow-up time from $N=1$, the integral of
$1/(X(N,D)^m+Y(N,D)^m)$ from $1$ to infinity, strictly decreases
with the initial imbalance.
The imbalance for $(X,Y)=(1,s)/(1+s)$ decreases with $s$. Hence
$(1+s)^{m-1}I_m(s)$ is strictly increasing.

For $a\ge1$, set $r=1/a\le1$. The identity
\[
 \big[(1+s^{1/a})I_m(s)^{1/[a(m-1)]}\big]^a
 =\frac{(1+s^r)^a}{1+s}
       \big[(1+s)I_m(s)^{1/(m-1)}\big]
\]
finishes the argument: the second factor is strictly increasing, and the
logarithmic derivative of the first is
$(s^{r-1}-1)/[(1+s^r)(1+s)]\ge0$.
\end{proof}
Consequently, when $a\ge1$ (all Hessian cases and $p\ge2$),
\begin{equation}\label{beta-constant}
 \CQ=\left[\frac{d_0\kappa^a}{m+1}
 B\left(\frac{m-1}{m+1},\frac1{m+1}\right)\right]^\beta<2C_0,
\end{equation}
where $B$ is the Euler beta function.
Its specialization to $(n,p)=(3,2)$ gives the constant in
Remark~\ref{three-dimensional-example}.

\begin{lemma}\label{minimum-subquadratic}
Let $0<a<1$ and $m>1$, and set
$a_c(m)=(3m+1)/(m+1)^2$. If $a\le a_c(m)$, the unique minimizer
of $\Phi_{a,\kappa}$ is $s=1$, and $\CQ=2C_0$.
If $a_c(m)<a<1$, there is a unique minimizer $s_*\in(0,1)$,
and $\CQ<2C_0$.
\end{lemma}
\begin{proof}
Let $\rho=1/a>1$ and, along \eqref{reverse-flow}, set
$R=(X^\rho+Y^\rho)^{1/\rho}$ and $s=Y/X\le1$. Then
\[
 \dot R=R^m h(s),\qquad
 h(s)=\frac{s^m+s^{\rho-1}}
 {(1+s^\rho)^{(m+\rho-1)/\rho}},\qquad \dot s>0\quad(s<1).
\]
The lifetime from $R=1$ is
$K(s)=(1+s^\rho)^{(m-1)/\rho}I_m(s)$, so differentiating the
remaining lifetime $R^{1-m}K(s)$ gives
\[
 K'(s)=\frac{(1+s^\rho)^{(m-1)/\rho}}{1-s^{m+1}}H(s),
 \qquad H(s)=(m-1)h(s)K(s)-1.
\]
To determine the shape of $h$, let $x=-\log s$,
$A=|m-\rho+1|$, and $B=m+\rho-1$. Then $A<B$ and
\[
 2\frac{\dd}{\dd x}\log h(e^{-x})
 =A\tanh(Ax/2)-B\tanh(\rho x/2).
\]
If $A\le\rho$, this is negative. If $A>\rho$, the ratio
$\tanh(Ax/2)/\tanh(\rho x/2)$ strictly decreases from $A/\rho$
to $1$, since $u/\sinh u$ strictly decreases for $u>0$.
As $A^2-B\rho=(m+1)^2-(3m+1)\rho$, it follows that $h$ is
strictly increasing on $(0,1)$ when
$\rho\ge(m+1)^2/(3m+1)$; otherwise it has one strict maximum.

Along the orbit starting from $R=1$ with ratio $s$, write $s(z)$
for the ratio at $R=z$. Separation gives
$K(s)=\int_1^\infty z^{-m}/h(s(z))\dd z$.
In the increasing case, $s(z)>s$ gives $H(s)<0$, so $K$ strictly
decreases and its unique minimum is at $1$.
In the one-maximum case, let $\sigma$ be the maximizer of $h$.
Then $H>0$ on $[\sigma,1)$, whereas $H(0)=-1$.
On $(0,\sigma)$, every zero satisfies $H'=(m-1)h'K>0$.
Thus $H$ has exactly one zero in $(0,1)$, and $K$ has
exactly one strict global minimum in $(0,1)$.
Since $\Phi_{a,\kappa}$ is a positive constant times
$K^{1/[a(m-1)]}$ and $\Phi_{a,\kappa}(1)=2d_0^{-\beta}C_0$,
the assertions follow.
\end{proof}

For the $p$-Laplacian, substituting $a=p-1$ and $m=n/(n-p)$
into $a=a_c(m)$ gives \eqref{subquadratic-threshold}.
Since $p-1-a_c(n/(n-p))$ is strictly increasing and changes sign
on $(1,2)$, there is a unique threshold $p_*(n)$, and the
classification stated in the introduction follows.

\begin{remark}\label{power-aggregation}
For $\rho>0$, write $H_\rho(x,y)=(x^\rho+y^\rho)^{1/\rho}$ and set
\[
 C_\rho=(d_0\kappa^a)^\beta\min_{0\le s\le1}
 (1+s^{\rho/a})^{1/\rho}I_m(s)^\beta,
 \qquad \rho_c=\frac{a(m+1)^2}{3m+1}.
\]
On branch~\textup{(iii)} of Theorem~\ref{mainjoint},
$\liminf_{x\to0}|x|^\kappa L(|x|)^\beta
H_\rho(w_1(x),w_2(x))\ge C_\rho$.
The unique minimizer is $s=0$ for $0<\rho\le a$, belongs to $(0,1)$
for $a<\rho<\rho_c$, and is $s=1$ for $\rho\ge\rho_c$.
Indeed, Corollary~\ref{finite-cluster-region} and coercivity of $H_\rho$
give the bound; raising the minimized expression to the power $\rho$
reduces the minimization to Lemmas~\ref{minimum-axis}--\ref{minimum-subquadratic}
with $a$ replaced by $a/\rho$.
The single radial pair realizing all of $\mathcal D_{\mathcal Q}$
attains these lower bounds simultaneously for every $\rho>0$.
\end{remark}

\subsection{The scalar inequality}

\begin{corollary}\label{hessian-scalar}
For either operator, a solution of $0\le\mathcal Q[w]\le\tau w^q$,
$\tau>0$, is asymptotically radial and satisfies the scalar version of
the trichotomy. Its zero-mass singular branch obeys the sharp bound
\begin{equation}\label{hessian-scalar-lower}
 \liminf_{x\to0}|x|^\kappa L(|x|)^\beta w(x)\ge\tau^{-\beta}C_0.
\end{equation}
For $\mathcal Q[w]=\tau w^q$, every
nonremovable singularity satisfies
\[
 |x|^\kappa L(|x|)^\beta w(x)\to\tau^{-\beta}C_0.
\]
\end{corollary}

The exact limits recover \cite{Aviles1987,MZII,MZKH}, while the
$p$-Laplace lower bound agrees with \eqref{scalar-lower-background}.

\begin{proof}[Proof of Corollary~\ref{hessian-scalar}]
We apply the spatial arguments in Sections~\ref{sec:local}--\ref{sec:mass}
to a single component. They give the scalar classification,
radiality, and, on the zero-mass singular branch, for
$P=\mu(B_{e^{-t}})/(\omega d_0)$,
\[
 e^{-\kappa t}w(e^{-t}\theta)
 =\kappa^{-1}P(t)^{1/a}(1+o(1)),\qquad
 0\le-P'\le\frac{\tau}{d_0\kappa^q}(1+o(1))P^m.
\]
Integrating $(P^{1-m})'$ gives \eqref{hessian-scalar-lower}.
For the exact equation the mass relation is an asymptotic equality,
so this is a limit. A positive mass at $0$ would force
$w^q\ge c|x|^{-n}$, contrary to local finiteness of the source measure.

For sharpness take $w=A r^{-\kappa}L(r)^{-\beta}$, where
$A=\tau^{-\beta}C_0$. The radial identity \eqref{common-radial-identity} gives, for large
$t=L(r)$,
\[
 \frac{\mathcal Q[w]}{w^q}
 =d_0a\beta A^{a-q}
 (\kappa-\beta/t)^{a-1}(\kappa-(\beta+1)/t)
 \le d_0a\beta\kappa^a A^{a-q}=\tau.
\]
For $a<1$, the last inequality follows by writing the ratio to
$\kappa^a$ as
$(1-\beta/(\kappa t))^a
 (1-1/(\kappa t-\beta))<1$; for $a\ge1$ it is immediate.
The expression is positive. In the Hessian case admissibility follows
from Lemma~\ref{radial-admissibility}. The limiting flux is zero,
and the required normalized limit equals $A$.
\end{proof}

\section{Radial constructions}\label{sec:construction}
Radial fluxes realize every state in Corollary~\ref{finite-cluster-region}.
Interleaving these stages with alternating flux reductions gives joint
sharpness and individual oscillations for the same pair.

\subsection{Admissible radial fluxes}

\begin{lemma}\label{flux-construction}
Let $P,Q>0$ be smooth, nonincreasing on $[1,\infty)$, tending to zero,
with $P(1)=Q(1)=1$ and
\begin{equation}\label{sufficient-control}
 0\le-P'\le d_0^{-1}\kappa^{-q}Q^m,\qquad
 0\le-Q'\le d_0^{-1}\kappa^{-q}P^m.
\end{equation}
Let $U,V$ solve
\begin{equation}\label{radial-filter}
 \begin{gathered}
 U'=P^{1/a}-\kappa U,\qquad V'=Q^{1/a}-\kappa V,\\[3pt]
 U(1)=V(1)=1/\kappa.
 \end{gathered}
\end{equation}
Then $(w_1,w_2)(r)=r^{-\kappa}(U,V)(L(r))$ is a smooth positive
radial solution of \eqref{upper-system} on $0<r<e^{-1}$,
with zero source masses at $0$.
\end{lemma}

\begin{proof}
We first show that $U,V$ are nonincreasing, tend to zero, and satisfy
\begin{equation}\label{filter-bounds}
 P(t)^{1/a}/\kappa\le U(t)\le1/\kappa,\qquad
 Q(t)^{1/a}/\kappa\le V(t)\le1/\kappa.
\end{equation}
The equation for $U-P^{1/a}/\kappa$ has a nonnegative source and zero
initial value, giving the lower bound; the same argument applies to $V$.
The filter equations give
$U',V'\le0$ and the upper bounds; exponential convolution gives
$U,V\to0$.

Direct differentiation and
\eqref{radial-flux}--\eqref{common-radial-identity} give
\[
 J_{w_1}(r)=P(L(r)),\qquad
 \mathcal Q[w_1]=-d_0r^{-n}P'(L(r)),
\]
and similarly for $w_2$. Since $w_2^q=r^{-n}V^q$,
\eqref{filter-bounds} and \eqref{sufficient-control} give the upper
inequalities. In the Hessian case, $(-w_i)'>0$ and the radial fluxes are
nondecreasing in $r$, so Lemma~\ref{radial-admissibility} gives admissibility,
including on constant-flux intervals. The limiting fluxes are zero, hence
both measures give zero mass to the origin. The physical functions
$e^{\kappa t}U(t)$ and $e^{\kappa t}V(t)$ are strictly increasing in $t$.
\end{proof}

\subsection{Realization of finite asymptotic states}\label{sec:finite-realization}

\begin{proof}[Realization in Corollary~\ref{finite-cluster-region}
and sharpness in Theorem~\ref{mainjoint}]
Fix $s\in[0,1]$, $\lambda\in(0,1]$, and $s_j\in(0,1]$ tending to $s$
(with $s_j\to0$ if $s=0$). Starting from $P(1)=Q(1)=1$,
we construct successive stages satisfying \eqref{sufficient-control}.
At each stage, reduce both fluxes to a common $B_j>0$ below half their
previous minimum, using flat smooth ramps with derivative magnitudes at
most $d_0^{-1}\kappa^{-q}B_j^m$. Keeping the fluxes above $B_j$ makes
this finite-time preparation admissible. Complete it before choosing $b$.

Choose $0<b<B_j/2$ and define
\[
 D_b=(1-s_j^{m+1})b^{m+1},\qquad
 C_b=(B_j^{m+1}-D_b)^{1/(m+1)}.
\]
Lower $Q$ from $B_j$ to $C_b$ over a unit interval while holding $P=B_j$.
For sufficiently small $b$, this is admissible because $B_j-C_b\to0$.
Starting from $(B_j,C_b)$, follow
\[
 P'=-\lambda d_0^{-1}\kappa^{-q}Q^m,\qquad
 Q'=-\lambda d_0^{-1}\kappa^{-q}P^m
\]
until reaching $(b,s_jb)$, preserving $P^{m+1}-Q^{m+1}=D_b$.
To smooth the joins, parameterize this orbit by its effective time as
$(\bar P(\sigma),\bar Q(\sigma))$, $0\le\sigma\le L_b$, with
\[
 L_b=\lambda^{-1}d_0\kappa^q
 \{T_m(b,s_jb)-T_m(B_j,C_b)\}.
\]
Choose a smooth nondecreasing function $\sigma=\sigma(t)$ on an interval
of length $L_b+O(1)$, constant near both endpoints, with
$0\le\sigma'\le1$, equal to $1$ except on two intervals of fixed length,
and satisfying $\sigma(\text{left})=0$ and
$\sigma(\text{right})=L_b$.  Set
$P(t)=\bar P(\sigma(t))$, $Q(t)=\bar Q(\sigma(t))$.  Then
\[
 -P'=\lambda\sigma'(t)d_0^{-1}\kappa^{-q}Q^m,
 \qquad -Q'=\lambda\sigma'(t)d_0^{-1}\kappa^{-q}P^m,
\]
Thus \eqref{sufficient-control} and the endpoint $(b,s_jb)$ are preserved,
and all joins are flat. The extra time is bounded independently of
$b$, so, with preceding stages fixed and $b\to0$, the terminal time satisfies
\begin{equation}\label{optimal-time}
 t_b=\lambda^{-1}d_0\kappa^q b^{1-m}I_m(s_j)(1+o(1)).
\end{equation}

At the endpoint, the filter convolution gives
\begin{equation}\label{tracking}
 U(t_b)=\frac{b^{1/a}}\kappa(1+o(1)),\qquad
 V(t_b)=\frac{(s_jb)^{1/a}}\kappa(1+o(1)).
\end{equation}
We now justify these estimates directly.  Since $L_b\asymp_j b^{1-m}$ while
$\log(1/b)=o(b^{1-m})$, the last $h_b=C_j\log(1/b)$ units of physical time
lie in the terminal part of the same orbit for $b$ small.  On that interval
the two fluxes are comparable, with constants depending on the fixed number
$s_j>0$, to their endpoint values; their relative variation is
$O_j(h_b b^{m-1})=o(1)$.  The bounded switch-off interval introduced above
only decreases that variation.  The contribution to the exponential
convolution from times earlier than $t_b-h_b$ is $O(e^{-\kappa h_b})$;
choosing $C_j$ sufficiently large makes this
$o((s_jb)^{1/a})$.  This proves \eqref{tracking}.  Here $s_j$ is fixed before $b$ is chosen sufficiently small.
Combining \eqref{optimal-time}--\eqref{tracking} gives
\[
 t_b^\beta(U(t_b),V(t_b))
 =\lambda^{-\beta}\bigl(d_0\kappa^aI_m(s_j)\bigr)^\beta
       (1,s_j^{1/a})+o(1).
\]
Choose $b=b_j$ so that the vector error is less than $1/j$, successive
endpoint times differ by more than $1$, and
$e^{\kappa t_b}U(t_b),e^{\kappa t_b}V(t_b)>j$.
This is possible since $t_b$ grows as $b^{1-m}$.
Repeated preparation gives $P,Q\to0$, and physical monotonicity gives
$w_1,w_2\to\infty$.
Lemma~\ref{flux-construction} gives the PDE assumptions. To obtain a
solution on $B_1\setminus\{0\}$, set
$\widetilde w_i(x)=\rho^\kappa w_i(\rho x)$ with $\rho=e^{-1}$.
The critical scaling preserves \eqref{upper-system}; its fixed translation
of $L(r)$ preserves finite cluster points and the normalized lower limit.
The endpoint vectors converge to
$\lambda^{-\beta}(A_1(s),A_2(s))$. Interchanging components realizes
the other half of $\mathcal D_{\mathcal Q}$ and proves the converse in
Corollary~\ref{finite-cluster-region}. Taking $\lambda=1$ and $s$ to
be a minimizer of $\Phi_{a,\kappa}$ gives equality in the normalized
lower limit in \eqref{joint-lower}.

For simultaneous realization, choose a sequence in
$\mathcal D_{\mathcal Q}\cap(0,\infty)^2$ having every point of
$\mathcal D_{\mathcal Q}$ as a subsequential limit.
At stage $j$, fix its parameters $s_j,\lambda_j$ and orientation
before taking $b_j$ small enough to make the vector error less than $1/j$.
The preparation and tracking estimates apply after any preceding stages.
The construction therefore retains vanishing masses and $w_i\to\infty$,
and its finite cluster set is exactly $\mathcal D_{\mathcal Q}$ by the
proved inclusion. This same pair also satisfies all the conclusions of
Theorem~\ref{individual-failure}.
\end{proof}

\subsection{Alternating small values of the individual components}

\begin{proof}[Proof of Theorem~\ref{individual-failure}]
Using Lemma~\ref{flux-construction}, we alternate oscillation stages
with stages from Section~\ref{sec:finite-realization} attaining the
joint constant.
Let $\theta\in C^\infty(\R)$ be nondecreasing, equal to $0$ on a neighborhood
of $(-\infty,0]$, equal to $1$ on a neighborhood of $[1,\infty)$, and satisfy
$0\le\theta'\le M$ for some fixed $M>0$.

At time $S$, let the selected flux be $X=A>0$ and the other $Y=B>0$,
selecting $X=P$ at odd cycles and $X=Q$ at even cycles.
Keep $Y$ constant, set $T=1+2Md_0\kappa^q A B^{-q/a}$, and choose
$0<\varepsilon<\min\{A,B,1\}/2$ below. Define
\begin{equation}\label{eq:ramp}
 X(t)=A-(A-\varepsilon)\theta((t-S)/T),\qquad Y(t)=B,
 \quad S\le t\le S+T.
\end{equation}
Then
\[
 0\le-X'(t)\le \frac{MA}T\le \frac12d_0^{-1}\kappa^{-q}B^{q/a},
 \qquad -Y'=0.
\]
Thus \eqref{sufficient-control} holds independently of $\varepsilon$.

After time $S+T$, hold both fluxes constant for an additional time
\begin{equation}\label{eq:waiting}
 L_\varepsilon=\frac{2}{a\kappa}\log(1/\varepsilon),
\end{equation}
with $L_\varepsilon\ge1$, and set $t_\varepsilon=S+T+L_\varepsilon$.
Let $Z,\widehat Z$ be the scaled components corresponding to $X,Y$.
The constant-flux interval gives
\[
 Z(t_\varepsilon)=\frac{\varepsilon^{1/a}}\kappa
       +\left(Z(S+T)-\frac{\varepsilon^{1/a}}\kappa\right)e^{-\kappa L_\varepsilon}.
\]
The bounds of Lemma~\ref{flux-construction} therefore give
\begin{equation}\label{eq:endpoint-estimates}
 \frac{\varepsilon^{1/a}}\kappa\le Z(t_\varepsilon)
 \le \frac{2\varepsilon^{1/a}}\kappa,
 \qquad \widehat Z(t_\varepsilon)\ge\frac{B^{1/a}}\kappa.
\end{equation}
Here we used $e^{-\kappa L_\varepsilon}=\varepsilon^{2/a}$.
Since $S,T,A,B$ are fixed before choosing $\varepsilon$, we have
\begin{equation}\label{eq:time-limits}
 t_\varepsilon^\beta Z(t_\varepsilon)\to0,
 \qquad t_\varepsilon^\beta\widehat Z(t_\varepsilon)\to\infty
 \quad(\varepsilon\to0).
\end{equation}
At the same time,
\[
 e^{\kappa t_\varepsilon}Z(t_\varepsilon)
 \ge\frac{e^{\kappa(S+T)}}\kappa\varepsilon^{-1/a}\to\infty,
 \qquad
 e^{\kappa t_\varepsilon}\widehat Z(t_\varepsilon)\to\infty.
\]
At cycle $j$, we can therefore choose $\varepsilon$ so that
\begin{equation}\label{eq:stage-choice}
 \begin{gathered}
 t_\varepsilon^\beta Z(t_\varepsilon)<\frac1j,\qquad
 t_\varepsilon^\beta\widehat Z(t_\varepsilon)>j,\\[3pt]
 e^{\kappa t_\varepsilon}U(t_\varepsilon)>j,\qquad
 e^{\kappa t_\varepsilon}V(t_\varepsilon)>j.
 \end{gathered}
\end{equation}
After this oscillation stage, insert one complete realization stage
from Section~\ref{sec:finite-realization}, with $\lambda=1$ and
$s_j\in(0,1]$ tending to a minimizer of $\Phi_{a,\kappa}$. Choose its
terminal parameter $b_j$ after all preceding stages have been fixed,
so that at its endpoint $\tau_j>t_\varepsilon+1$,
\[
 \left|\tau_j^\beta\bigl(U(\tau_j)+V(\tau_j)\bigr)-\CQ\right|
 <\frac1j.
\]
Choose $s_j$ sufficiently close to the minimizing set, then use
\eqref{optimal-time}--\eqref{tracking} and continuity of $\Phi_{a,\kappa}$.
Preparation starts from the current flux values, without resetting the
filters; the terminal tracking estimate makes their earlier contribution
exponentially small. At $\tau_j$, start the next cycle with the selected
component exchanged.

Start at $S=1$ with $P=Q=1$. Since $L_\varepsilon\ge1$, the cycles
have no finite accumulation time. Flat ramps and reparameterizations
make $P,Q$ smooth; preparation below half the preceding minimum gives
$P,Q\to0$. Lemma~\ref{flux-construction} proves the system and
$r^\kappa w_1,r^\kappa w_2\to0$. Physical monotonicity and the last
two inequalities in \eqref{eq:stage-choice} give $w_1,w_2\to\infty$.
Since the selected component alternates, the first two inequalities in
\eqref{eq:stage-choice} give both limits in \eqref{individual-oscillation}
for each component.
The values at $\tau_j$, together with \eqref{joint-lower}, give
\eqref{simultaneous-joint-attainment}.
\end{proof}

The sources may vanish on intervals of logarithmic radii. Positive lower
coupling, treated next, rules out this behavior.

\section{Comparison across the isolated point}\label{sec:coupling}

\subsection{Comparison tools}

Comparison across the puncture uses an inner cutoff for $-\Delta_p$
and a fundamental-solution perturbation for the Hessian operator.
It gives comparability under lower coupling and a bounded difference
for the exact symmetric system.

\begin{lemma}\label{inner-boundary}
For the $p$-Laplace parameters, let $w_1,w_2\ge0$ be locally bounded
weak supersolutions, $f_i=-\Delta_pw_i$, with
$0\le f_1+f_2\le D(w_1+w_2)^q$ and
$|x|^\kappa(w_1+w_2)\to0$. Then
\begin{equation}\label{small-gradient}
 |x|^{\kappa+1}(|\nabla w_1|+|\nabla w_2|)\to0.
\end{equation}
\end{lemma}

Consequently, any cutoff $0\le\eta_\varepsilon\le1$ equal to $0$ on
$B_\varepsilon$ and $1$ outside $B_{2\varepsilon}$, with
$|\nabla\eta_\varepsilon|\le C/\varepsilon$, satisfies
\begin{equation}\label{cutoff-error}
 \int_{\varepsilon<|x|<2\varepsilon}
 (|\nabla w_1|^a+|\nabla w_2|^a)
 |\nabla\eta_\varepsilon|\dd x\to0.
\end{equation}
\begin{proof}
We apply \eqref{gradient-input} on balls with radius comparable to $|x|$.
Their suprema are $o(r^{-\kappa})$ and their densities are $o(r^{-n})$.
Both terms in that estimate are $o(r^{-\kappa-1})$, since
$(n-1)/a=\kappa+1$. The integral in \eqref{cutoff-error} is therefore
$o(r^{n-1-a(\kappa+1)})=o(1)$.
\end{proof}

\begin{lemma}\label{hessian-bounded-difference}
For the Hessian operator, suppose $b_1=b_2=0$ and
\[
 \mathcal Q[w_1]=A(x)w_2^q,\qquad
 \mathcal Q[w_2]=B(x)w_1^q,\qquad 0<c\le A,B\le C<\infty.
\]
If $\lambda>0$ and $A\le\lambda^{k+q}B$ on $B_R\setminus\{0\}$, then
\begin{equation}\label{hessian-weighted-comparison}
 w_1\le\lambda w_2+M_{\lambda,R},\qquad
 M_{\lambda,R}=\max_{\partial B_R}(w_1-\lambda w_2)_+.
\end{equation}
\end{lemma}
\begin{proof}
We write $u=-w_1$, $v=-w_2$, and $M=M_{\lambda,R}$.
If necessary, replace $A,B$ by the continuous ratios
$\sigma_k(D^2u)/(-v)^q$ and $\sigma_k(D^2v)/(-u)^q$;
the stated bounds and comparison condition are preserved.
By Lemma~\ref{fundamental-vector}, $r^\kappa(-u-v)\to0$.
For $\varepsilon>0$, the function
$v_\varepsilon=\lambda v-\varepsilon r^{-\kappa}$ is $k$-admissible and
\[
 \sigma_k(D^2v_\varepsilon)\ge\lambda^k\sigma_k(D^2v).
\]
Indeed, concavity and homogeneity of $\sigma_k^{1/k}$ on the convex
admissible cone give superadditivity~\cite[Section~2]{TW2}, and
$-r^{-\kappa}$ is $k$-harmonic away from $0$.
The function $v_\varepsilon-u-M$ is nonpositive on $\partial B_R$
and every sufficiently small inner sphere. At a positive interior maximum,
$w_1>\lambda w_2$ and $D^2v_\varepsilon\le D^2u$, so ellipticity gives
\[
 \sigma_k(D^2u)\ge\sigma_k(D^2v_\varepsilon)
 \ge\lambda^kB w_1^q
 >\lambda^{k+q}B w_2^q
 \ge A w_2^q=\sigma_k(D^2u),
\]
a contradiction. Letting $\varepsilon\to0$ yields
\eqref{hessian-weighted-comparison}, completing the proof.
\end{proof}

\subsection{Comparison and growth under lower coupling}

\begin{proposition}\label{maintwosided}
Suppose
\begin{equation}\label{two-sided}
 cw_2^q\le\mathcal Q[w_1]\le Cw_2^q,\qquad
 cw_1^q\le\mathcal Q[w_2]\le Cw_1^q,
 \qquad 0<c\le C<\infty.
\end{equation}
For $K=(C/c)^{1/(a+q)}$ and every sufficiently small $R>0$,
\begin{equation}\label{comparison-main}
 w_i\le Kw_{3-i}+M_i
 \quad\text{in }B_R\setminus\{0\},\qquad i=1,2.
\end{equation}
Moreover, either both singularities are removable or
\begin{equation}\label{component-orders}
 w_i(x)\asymp |x|^{-\kappa}L(|x|)^{-\beta},
 \qquad i=1,2,\quad x\to0.
\end{equation}
In the exact symmetric case $c=C=1$, for every $B_R\Subset B_1$,
\begin{equation}\label{bounded-difference}
 |w_1-w_2|\le\max_{\partial B_R}|w_1-w_2|
 \quad\text{in }B_R\setminus\{0\}.
\end{equation}
\end{proposition}

Here $M_i=\max_{\partial B_R}(w_i-Kw_{3-i})_+$.
The constants in \eqref{component-orders} depend only on the operator
parameters and $c,C$; the radius of validity may depend on the solution.

\begin{proof}
Write $f_i=\mathcal Q[w_i]$. Lemma~\ref{vanishing} gives $b_1=b_2=0$ and
$|x|^\kappa(w_1+w_2)\to0$. The source bound and
Lemma~\ref{radial-comparison} make the sum either bounded or uniformly
divergent. In the bounded case, bounded densities give removability
by Section~\ref{sec:local}.

For $-\Delta_p$, let $\Aop(\xi)=|\xi|^{p-2}\xi$. We fix $M>M_1$
and define $w=Kw_2$ and
$h=(w_1-w-M)_+$. Then $h=0$ near $\partial B_R$, and on $\{h>0\}$,
\begin{equation}\label{comparison-sign}
 f_1-K^af_2\le Cw_2^q-cK^aw_1^q
 \le(C-cK^{a+q})w_2^q=0.
\end{equation}
For $N>0$, define $T_N(s)=\min\{s,N\}$ for $s\ge0$. Subtract the
weak equations for $w_1$ and $w$, and test with
$\eta_\varepsilon T_N(h)$, where $\eta_\varepsilon$ is as in
Lemma~\ref{inner-boundary}. The right-hand side is nonpositive, while
\begin{align}\label{comparison-integral}
 &\int \eta_\varepsilon\boldsymbol1_{\{0<h<N\}}
  (\Aop(\nabla w_1)-\Aop(\nabla w))
      \cdot(\nabla w_1-\nabla w)\dd x\notag\\
 &\qquad\le
 N\int_{\varepsilon<|x|<2\varepsilon}
   (|\nabla w_1|^a+K^a|\nabla w_2|^a)
       |\nabla\eta_\varepsilon|\dd x=o(1).
\end{align}
The integrand on the first line is nonnegative. Fatou's lemma and strict
monotonicity show that $\nabla T_N(h)=0$ almost everywhere on the
punctured ball. This truncation is locally in $W^{1,p}$, and the punctured
ball is connected; hence it is constant. It vanishes near the outer
boundary, so $h=0$. Let $M\to M_1$ from above and interchange the components.

For the Hessian operator, write
$A=\mathcal Q[w_1]/w_2^q$ and $B=\mathcal Q[w_2]/w_1^q$.
Then $c\le A,B\le C$, and Lemma~\ref{hessian-bounded-difference} gives
\eqref{comparison-main} with $\lambda=K$ in both directions.

For $c=C=1$, take $K=1$ in both directions on any $B_R\Subset B_1$
to obtain \eqref{bounded-difference}.

If the sum diverges, \eqref{comparison-main} forces both components
to infinity and makes them comparable. Set $S=w_1+w_2$, $r=e^{-t}$, and
\[
 P(t)=\frac{\mu_1(B_r)}{\omega d_0},\qquad
 Q(t)=\frac{\mu_2(B_r)}{\omega d_0},\qquad J=P+Q.
\]
The mass expansion \eqref{mass-expansion} gives, uniformly in direction,
\[
 r^\kappa S(r\theta)
 =\kappa^{-1}(P^{1/a}+Q^{1/a})+o(J^{1/a})\asymp J^{1/a},
\]
since $P^{1/a}+Q^{1/a}\asymp(P+Q)^{1/a}$. Also,
$w_1^q+w_2^q\asymp S^q$ for every $q>0$. Thus \eqref{two-sided} and
$\kappa q=n$ yield, almost everywhere,
\begin{equation}\label{two-sided-mass-ode}
 -J'(t)=\frac{r^n}{d_0}\av_{\Sph}(f_1+f_2)(r\theta)
 \asymp J(t)^m.
\end{equation}
Integration gives $J(t)\asymp t^{-1/(m-1)}=t^{-a\beta}$, hence
$S(r\theta)\asymp r^{-\kappa}t^{-\beta}$. Finally,
\eqref{comparison-main} and $S\to\infty$ give $w_i\asymp S$, proving
\eqref{component-orders}. All comparison constants depend only on the
operator parameters and $c,C$ after a solution-dependent increase of $t$.
\end{proof}

\section{Exact asymptotics and synchronization}\label{sec:exact}

Comparison and the mass equations give the common leading term.
Smooth annular convergence and a cylinder estimate then yield derivative
limits and angular decay. The averaged equations determine the common
second-order profile and synchronization coefficient, followed by radial
realization and an extension to averaged coupling.

\subsection{The common leading term}\label{sec:exact-leading}

\begin{proof}[Proof of the dichotomy and \eqref{exact-common-limit} in
Theorem~\ref{mainexact}]
For either operator in~\eqref{parameters}, let $w_1,w_2$ be positive functions. Lemma~\ref{vanishing} gives $b_1=b_2=0$.
Proposition~\ref{maintwosided}, with $c=C=1$, gives the
removable-or-divergent alternative and \eqref{bounded-difference}.
On the singular branch, $w_1,w_2\to\infty$, so this bound implies
$w_1/w_2\to1$.
Define
\[
 P(t)=\frac{\mu_1(B_{e^{-t}})}{\omega d_0},\qquad
 Q(t)=\frac{\mu_2(B_{e^{-t}})}{\omega d_0},\qquad J=P+Q.
\]
Proposition~\ref{mass-approx} gives $P/J,Q/J\to1/2$: for any convergent
subsequence of $P/J\in[0,1]$, divide \eqref{mass-expansion} by $J^{1/a}$.
Then $w_1/w_2\to1$ forces its limit to be $1/2$, without dividing by
a possibly smaller component.

The exact equations and $\kappa q=n$ now give, almost everywhere,
\begin{equation}\label{exact-total-mass}
 \begin{aligned}
 -J'(t)&=(1+o(1))(d_0\kappa^q)^{-1}
       (P(t)^m+Q(t)^m)\\[3pt]
       &=(1+o(1))(d_0\kappa^q)^{-1}2^{1-m}J(t)^m.
 \end{aligned}
\end{equation}
Since $J\to0$, integration yields
\begin{equation}\label{exact-mass-rate}
 J(t)^{1-m}\sim(m-1)(d_0\kappa^q)^{-1}2^{1-m}t.
\end{equation}
Substitution in \eqref{mass-expansion} gives, uniformly in $\theta$,
\[
 \begin{aligned}
 \lim_{t\to\infty}t^\beta e^{-\kappa t}w_i(e^{-t}\theta)
 &=\kappa^{-1}2^{-1/a}
   [(m-1)(d_0\kappa^q)^{-1}2^{1-m}]^{-\beta}\\
 &=\kappa^{-1}[(m-1)(d_0\kappa^q)^{-1}]^{-\beta}
 =(d_0a\beta\kappa^a)^\beta=C_0.
 \end{aligned}
\]
Here $(m-1)\beta=1/a$ cancels the powers of $2$, and
$q\beta-1=a\beta$ gives the last equality. This proves
\eqref{exact-common-limit} for both operators.
\end{proof}

\subsection{Perturbative Hessian regularity}

\begin{lemma}\label{hessian-annular-smooth}
Let $u=-w_1$, $v=-w_2$ for a nonremovable solution of the exact
Hessian system. Then, as $r\to0$,
\[
 r^\kappa L(r)^\beta(u(r\,\cdot),v(r\,\cdot))
 \to-C_0|\cdot|^{-\kappa}(1,1)
 \quad\text{in }C^\infty_{\loc}(\R^n\setminus\{0\}).
\]
\end{lemma}

\begin{proof}
Set $t=L(r)$ and
\[
 u_r(y)=r^\kappa t^\beta u(ry),\qquad
 v_r(y)=r^\kappa t^\beta v(ry),\qquad
 \phi(y)=-C_0|y|^{-\kappa}.
\]
\noindent\emph{Uniform H\"older control.}
The leading asymptotic gives local uniform convergence, and scaling gives
\[
 \sigma_k(D^2u_r)=t^{-1}(-v_r)^q,\qquad
 \sigma_k(D^2v_r)=t^{-1}(-u_r)^q.
\]
On fixed annuli the absolute values of the components are uniformly
bounded above and away from zero. Estimate~\eqref{hessian-oscillation}
gives a uniform $C^{0,\alpha_0}$ bound for some $\alpha_0\in(0,1)$.
Thus, for $0<\alpha<\alpha_0$, both right-hand sides have
$C^{0,\alpha}$ norm $O(t^{-1})$ on smaller annuli. This uses the
equations and potential estimate, not merely uniform convergence.

\smallskip
\noindent\emph{Ellipticity near the limiting function.}
Fan's inhomogeneous small-perturbation theorem~\cite[Theorem~1.7]{Fan} applies
even though $\sigma_k(D^2\phi)=0$. Set
$A=D^2\phi(y_0)$ at a point $y_0$ of the annulus and
$b=C_0\kappa |y_0|^{-\kappa-2}$.
In a frame whose first vector is radial,
\[
 A=b\,\operatorname{diag}(-\kappa-1,1,\ldots,1).
\]
The eigenvalues of $T_{k-1}(A)=D\sigma_k(A)$ are
\[
 \binom{n-1}{k-1}b^{k-1},\qquad
 \frac{\kappa+1}{n-1}\binom{n-1}{k-1}b^{k-1}
 \quad\text{(with multiplicity $n-1$)}.
\]
Moreover, $\sigma_j(A)>0$ for $1\le j<k$, while $\sigma_k(A)=0$.
Thus the linearization is positive definite, uniformly for $y_0$ in
a fixed compact annulus.

Define a globally degenerate elliptic extension as follows. For $s\ge0$ and $M\in\operatorname{Sym}(n)$, let
$F_s(M)$ be the unique real number $\zeta$ such that
\[
 M-\zeta I\in\overline\Gamma_k,\qquad
 \sigma_k(M-\zeta I)=s,
\]
where for $s=0$ the matrix lies on $\partial\Gamma_k$, and for $s>0$
it lies in $\Gamma_k$. Equivalently,
\[
 F_s(M)=\sup\{\zeta:M-\zeta I\in\overline\Gamma_k,
                         \ \sigma_k(M-\zeta I)\ge s\}.
\]
Existence and uniqueness follow by moving along $M-\zeta I$ from
$\zeta=-\infty$ to its first intersection with
$\partial\Gamma_k$: $\sigma_k$ decreases strictly from infinity to zero, since
$\partial_\zeta\sigma_k(M-\zeta I)=-(n-k+1)\sigma_{k-1}(M-\zeta I)<0$
in $\Gamma_k$.
Stability of the cone under addition of nonnegative matrices and
monotonicity of $\sigma_k$ imply global degenerate ellipticity of $F_s$.
Together with $F_s(M+cI)=F_s(M)+c$, this gives global Lipschitz
continuity in $M$. Since $\operatorname{tr}T_{k-1}(A)>0$, the implicit function theorem
gives a smooth extension near $(A,0)$ in $(M,s)$, including small
negative $s$. In that neighborhood,
\[
 \begin{aligned}
 D_MF_s(M)&=
 \frac{T_{k-1}(M-F_s(M)I)}
      {\operatorname{tr}T_{k-1}(M-F_s(M)I)},\\[3pt]
 \partial_sF_s(M)&=
 -\frac1{\operatorname{tr}T_{k-1}(M-F_s(M)I)}.
 \end{aligned}
\]
The first matrix is uniformly positive definite near $(A,0)$.
For $s\ge0$, positivity of the lower elementary symmetric functions
and $\sigma_k=s$ place the implicit branch in the admissible cone,
where it agrees with the global definition.

\smallskip
\noindent\emph{Small perturbations and smooth convergence.}
Let $Q$ be the second-order Taylor polynomial of $\phi$ at $y_0$.
Choose a small fixed $\rho>0$ and set
\[
 h_r(z)=\rho^{-2}\bigl(u_r(y_0+\rho z)-Q(y_0+\rho z)\bigr),
 \qquad s_r(z)=t^{-1}(-v_r(y_0+\rho z))^q.
\]
Define
\[
 G_r(N,z)=F_{s_r(z)}(A+N)-F_{s_r(z)}(A),\qquad
 g_r(z)=-F_{s_r(z)}(A).
\]
Then $G_r(D^2h_r,z)=g_r(z)$ and $G_r(0,z)=0$.
Indeed, the original admissible equation is equivalent to
$F_{s_r(z)}(A+D^2h_r)=0$. The operators $G_r$ are globally degenerate elliptic and uniformly
elliptic in a fixed neighborhood of $N=0$. Choose a fixed small matrix
radius so that $A+N$ stays there uniformly for $y_0$ in the compact annulus. The derivative formulas, and
$\|s_r\|_{C^{0,\alpha}}=O(t^{-1})$, give the uniform estimates
\[
 \|D_NG_r(N,z)-D_NG_r(N',z')\|
 \le C\bigl(\|N-N'\|+t^{-1}|z-z'|^\alpha\bigr),
\]
and
\[
 |G_r(N,z)-G_r(N,z')|\le Ct^{-1}|z-z'|^\alpha,
 \qquad \|g_r\|_{C^{0,\alpha}(B_1)}\le Ct^{-1}.
\]
Thus the matrix derivatives share a local modulus of continuity,
and the coefficient oscillation and right-hand side are H\"older small.
For each fixed $r$, $s_r>0$ and $s_r\in C^1(B_1)$ by the original
$C^2$ regularity. Since the defining matrix lies in $\Gamma_k$ with
$\operatorname{tr}T_{k-1}>0$, the implicit function theorem gives
smoothness of $(M,s)\mapsto F_s(M)$ for $s>0$. Thus $G_r$ is jointly
$C^1$ in $(N,z)$. Fan's matrix-variable hypothesis holds with the
common local modulus of $D_NG_r$ established above; no bound on
$D_zs_r$ uniform in $r$ is needed. Since $G_r$ has no first-order or zeroth-order arguments,
the remaining structure conditions are automatic.

Finally,
\[
 \|h_r\|_{L^\infty(B_1)}
 \le \rho^{-2}\|u_r-\phi\|_{L^\infty(B_\rho(y_0))}+C\rho.
\]
Choose $\rho$ small, then $r$ small. Fan's theorem gives a uniform
$C^{2,\alpha}$ bound on $B_{1/2}$, also for $v_r$.
A finite covering and compactness upgrade the known uniform
convergence to $C^{2,\alpha'}$ convergence for every
$0<\alpha'<\alpha$. Thus the equations are uniformly elliptic
along the actual normalized solutions on smaller annuli.
Schauder estimates, using difference quotients and then successive
differentiation, give uniform bounds of every fixed order. Compactness completes the proof.
\end{proof}
\subsection{Quantitative angular decay}

We first prove exponential decay of the nonconstant spherical part.
For a cylinder function write
\[
 \langle h\rangle(t)=\frac1\omega\int_{\Sph}h(t,\theta)\dd\theta,
 \qquad \Pi h=h-\langle h\rangle.
\]
Cylinder derivatives refer to $t$ and the round metric on $\Sph$.

\begin{lemma}\label{angular-decay}
Let $(w_1,w_2)$ be a nonremovable solution of \eqref{exact-symmetric}
for either operator, and write $\bar w_i(r)=\av_{\Sph}w_i(r\theta)$.
For some $\alpha\in(0,1)$ and every $0<\eta<1$,
\begin{equation}\label{absolute-angular-decay}
 \sum_{i=1}^2\|w_i(e^{-\cdot}\,\cdot)-\bar w_i(e^{-\cdot})\|
 _{C^{2,\alpha}([t,t+1]\times\Sph)}=O(e^{-\eta t}).
\end{equation}
\end{lemma}

\begin{proof}
For either operator, set
\[
 U(t,\theta)=e^{-\kappa t}w_1(e^{-t}\theta),\qquad
 V(t,\theta)=e^{-\kappa t}w_2(e^{-t}\theta),\qquad
 W=t^\beta U,\quad Z=t^\beta V.
\]
The translates of $W,Z$ converge smoothly to $C_0$ on compact
cylinder strips. For the Hessian operator this is
Lemma~\ref{hessian-annular-smooth}. For $-\Delta_p$, the normalized functions
$w_{i,r}(\xi)=r^\kappa L(r)^\beta w_i(r\xi)$ converge uniformly on compact annuli to $C_0|\xi|^{-\kappa}$ and satisfy
$-\Delta_pw_{i,r}=L(r)^{-1}w_{3-i,r}^q$.
Interior $C^{1,\alpha}$ estimates~\cite{Tolksdorf} give $C^1$
convergence on smaller annuli. The nonvanishing limiting gradient and positivity permit Schauder
estimates and differentiation to every fixed order, also for $p<2$
or $q<1$. Compactness gives smooth convergence.

\smallskip
\noindent\emph{The two cylinder operators.}
Set $b(t)=\kappa-\beta/t$. In the $p$-Laplace case write
\[
 D_th=h_t+bh,\qquad
 F_t(h)=\big((D_th)^2+|\nabla_\theta h|^2\big)^{(p-2)/2},
\]
and
\[
 \mathcal N_t(h)=\partial_t(F_t(h)D_th)
 +\operatorname{div}_\theta(F_t(h)\nabla_\theta h)
 -\frac{a\beta}{t}F_t(h)D_th.
\]
The exact cylinder equations are
\begin{equation}\label{normalized-cylinder-equations}
 \mathcal N_t(W)+t^{-1}Z^q=0,\qquad
 \mathcal N_t(Z)+t^{-1}W^q=0.
\end{equation}
The limiting differential at $C_0$ is
\[
 a(\kappa C_0)^{p-2}
 \bigl(h_{tt}+\kappa h_t+a^{-1}\Delta_{\Sph}h\bigr),
 \qquad a^{-1}=\frac{\kappa+1}{n-1}.
\]

For the Hessian operator, define
\[
 \mathsf H_\kappa[h]=
 \begin{pmatrix}
 -h_{tt}-(2\kappa+1)h_t-\kappa(\kappa+1)h
   & (\nabla_\theta h_t+(\kappa+1)\nabla_\theta h)^T\\
 \nabla_\theta h_t+(\kappa+1)\nabla_\theta h
   & (h_t+\kappa h)I-\nabla_\theta^2h
 \end{pmatrix}.
\]
Then $D^2(-w_1)=r^{-\kappa-2}\mathsf H_\kappa[U]$ in the
orthonormal radial--tangential frame. Define $\mathsf H_t[h]$ by
replacing its radial entry with
$-h_{tt}-(2b+1)h_t-(b'+b(b+1))h$, its mixed entry with
$\nabla_\theta h_t+(b+1)\nabla_\theta h$, and its tangential block
with $(h_t+bh)I-\nabla_\theta^2h$. The equations become
\begin{equation}\label{hessian-cylinder-exact}
 \sigma_k(\mathsf H_t[W])=t^{-1}Z^q,\qquad
 \sigma_k(\mathsf H_t[Z])=t^{-1}W^q.
\end{equation}
At $C_0$, the matrix $\mathsf H_\kappa[C_0]$ has eigenvalues
$-(\kappa+1)\lambda,\lambda,\ldots,\lambda$, where $\lambda=\kappa C_0$.
Its Newton tensor has radial eigenvalue
$T_r=\binom{n-1}{k-1}\lambda^{k-1}$ and tangential eigenvalue
$T_\theta=T_r(\kappa+1)/(n-1)$. Hence
\[
 D\bigl[\sigma_k(\mathsf H_\kappa[\cdot])\bigr]_{C_0}h
 =-T_r\left(h_{tt}+\kappa h_t+
              \frac{\kappa+1}{n-1}\Delta_{\Sph}h\right).
\]

\smallskip
\noindent\emph{Projection and decay.}
In either system, subtract the expressions at the spherical averages
of $W,Z$ and apply $\Pi$. These radial expressions are annihilated
by $\Pi$; they need not solve the equations.
Integrating the differential along the intervening segments and dividing
by the corresponding limiting radial coefficient gives
\begin{equation}\label{projected-angular-equation}
 \begin{gathered}
 \mathcal L_0\mathbf w+\Pi\mathcal R\mathbf w=0,
 \qquad \mathbf w=(\Pi W,\Pi Z),\\[3pt]
 \mathcal L_0=\partial_t^2+\kappa\partial_t+
                  \frac{\kappa+1}{n-1}\Delta_{\Sph}.
 \end{gathered}
\end{equation}
For the coefficient bound needed by the cylinder inverse, write $I_t=[t,t+1]\times\Sph$,
$I_t^*=[t-1,t+2]\times\Sph$, and set
\[
 \delta(t)=t^{-1}+\|W-C_0\|_{C^{2,\alpha}(I_t^*)}
                    +\|Z-C_0\|_{C^{2,\alpha}(I_t^*)}\to0.
\]
For $s\in[0,1]$, let $W_s=\langle W\rangle+s\Pi W$ and
$Z_s=\langle Z\rangle+s\Pi Z$. Boundedness of spherical averaging on these H\"older spaces gives
the same control along both segments.
Write $\lambda=\kappa C_0$.

In the $p$-Laplace case define
\[
 \xi_s=(D_tW_s,\nabla_\theta W_s),\qquad
 \mathsf A(\xi)=|\xi|^{p-2}I+(p-2)|\xi|^{p-4}\xi\otimes\xi.
\]
The principal tensor of the differential of $\mathcal N_t$ at $W_s$
is $\mathsf A(\xi_s)$. Since $\xi_s$ stays uniformly away from $0$,
\[
 \sup_{0\le s\le1}\left\|
 \frac{\mathsf A(\xi_s)}{a\lambda^{p-2}}
       -\operatorname{diag}(1,a^{-1}I_{n-1})
 \right\|_{C^{1,\alpha}(I_t)}\le C\delta(t).
\]
This also controls the coefficients obtained by expanding the divergence. In the Hessian case, with
$c_\theta=(\kappa+1)/(n-1)$, the corresponding estimate is
\[
 \sup_{0\le s\le1}\left\|
 \frac{T_{k-1}(\mathsf H_t[W_s])}{T_r}
       -\operatorname{diag}(1,c_\theta I_{n-1})
 \right\|_{C^{0,\alpha}(I_t)}\le C\delta(t).
\]
Indeed the differential is $T_{k-1}(\mathsf H_t[W_s]):\mathsf H_t[h]$,
and the Newton tensor is polynomial. Division by $-T_r$ gives the
principal part of $\mathcal L_0$ with the stated small coefficient errors. Both estimates hold with $Z_s$ in place of $W_s$.
The terms involving $b-\kappa$ and $b'$ are controlled by $Ct^{-1}$.
Also $W_s,Z_s\ge C_0/2$ on sufficiently late strips, so the coupling
coefficient $t^{-1}q\int_0^1Z_s^{q-1}\dd s$, and its counterpart with
$W_s$, have $C^{0,\alpha}$ norm at most $Ct^{-1}$, even when $q<1$.

It follows that every coefficient of $\mathcal R$, written in covariant
nondivergence form, has $C^{0,\alpha}(I_t)$ norm at most $C\delta(t)$.
The product estimate and boundedness of $\Pi$ therefore give, for every
$\mathbf h\in C^{2,\alpha}(I_t;\R^2)$,
\begin{equation}\label{angular-perturbation-norm}
 \|\Pi\mathcal R\mathbf h\|_{C^{0,\alpha}(I_t)}
 \le\varepsilon(t)\|\mathbf h\|_{C^{2,\alpha}(I_t)},
 \qquad \varepsilon(t)=C\delta(t)\to0.
\end{equation}
Thus $\sup_{t\ge T}\varepsilon(t)\to0$ as $T\to\infty$, as required
by Lemma~\ref{cylinder-inverse}, without a quantitative convergence rate.

The first nonconstant spherical eigenvalue is $n-1$, so the
corresponding characteristic roots are $1$ and $-(\kappa+1)$.
Lemma~\ref{cylinder-inverse}, applied to the bounded projected pair,
gives for every $0<\sigma<\kappa+1$,
\begin{equation}\label{normalized-angular-exponential}
 \|\Pi W\|_{C^{2,\alpha}([t,t+1]\times\Sph)}
 +\|\Pi Z\|_{C^{2,\alpha}([t,t+1]\times\Sph)}
 =O(e^{-\sigma t}).
\end{equation}
This inverse estimate does not assume exponential decay beforehand.
We choose $\kappa+\eta<\sigma<\kappa+1$ and use
\[
 w_1(e^{-t}\theta)-\bar w_1(e^{-t})
 =e^{\kappa t}t^{-\beta}\Pi W(t,\theta),
\]
with the analogous identity for $w_2$. These identities and their
cylinder derivatives prove \eqref{absolute-angular-decay}.
\end{proof}

The same smooth annular convergence gives, for both operators and $i=1,2$,
\begin{equation}\label{exact-derivative-limits}
 \begin{aligned}
 r^{\kappa+1}L(r)^\beta\nabla w_i(r\theta)
   &\longrightarrow-\kappa C_0\theta,\\
 r^{\kappa+2}L(r)^\beta D^2w_i(r\theta)
   &\longrightarrow\kappa C_0
       \bigl((\kappa+2)\theta\otimes\theta-I\bigr),
 \end{aligned}
\end{equation}
uniformly in $\theta$, by differentiating the rescaled limit
$C_0|y|^{-\kappa}$ on a fixed annulus.

\subsection{The common profile and the difference expansion}

The spherical means and ball masses satisfy the following system,
whose average and difference determine the two expansions. Separating
exponential forcing from errors proportional to the difference preserves
the angular decay rate when the synchronization coefficient vanishes.

\begin{lemma}\label{averaged-difference-lemma}
Let $a,\kappa>0$, $q>a$, and set
\[
 \beta=(q-a)^{-1},\qquad \gamma=q\beta,\qquad
 C_{\mathrm{av}}=(a\beta\kappa^a)^\beta.
\]
Suppose $f,g,P,Q>0$ are $C^1$ on $[T,\infty)$, $T>1$, and satisfy
\begin{equation}\label{averaged-four-ode}
 \begin{gathered}
 \begin{alignedat}{2}
 f'&=P^{1/a}-\kappa f+E_1,&\qquad g'&=Q^{1/a}-\kappa g+E_2,\\[3pt]
 P'&=-g^q+E_3,& Q'&=-f^q+E_4,
 \end{alignedat}\\[4pt]
 E_i=O(e^{-(\kappa+c_0)t})
 \end{gathered}
\end{equation}
for some $c_0>0$. If $f,g\sim C_{\mathrm{av}}t^{-\beta}$ and
$H(t)=e^{\kappa t}(f-g)=O(1)$, then there is a common constant
$b_*\in\R$ such that
\begin{equation}\label{averaged-common-phase}
 \begin{aligned}
 f(t),g(t)=C_{\mathrm{av}}t^{-\beta}\left[
 1-\frac{\beta\gamma}{\kappa}\frac{\log t}{t}
 +\frac{b_*}{t}
 +O\left(\frac{(\log t)^2}{t^2}\right)\right].
 \end{aligned}
\end{equation}
Moreover, for some $d\in\R$,
\begin{equation}\label{averaged-difference-asymptotic}
 H(t)=dt^{-\gamma}-\frac{d\gamma^2}{\kappa}t^{-\gamma-1}\log t
       +O(t^{-\gamma-1}).
\end{equation}
If $d=0$, then $H(t)=O(e^{-\eta t})$ for every $0<\eta<c_0$.
\end{lemma}

\begin{proof}
Since $q\beta>1$, $P,Q$ have finite nonnegative limits, which vanish
by the filter equations and $f,g\to0$. Integration from infinity yields
\[
 P,Q\sim\frac{C_{\mathrm{av}}^q}{q\beta-1}t^{1-q\beta}
       =(\kappa C_{\mathrm{av}})^a t^{-a\beta}.
\]
\noindent\emph{Step 1: accuracy of the common term.}
Set
\[
 S=\frac{f+g}{2},\qquad R=\frac{P+Q}{2},\qquad
 \delta=\frac{P-Q}{2}.
\]
The assumed bound on $H$ gives $f-g=O(e^{-\kappa t})$.
From \eqref{averaged-four-ode},
\[
 \delta'=\tfrac12(f^q-g^q)+\tfrac12(E_3-E_4),\qquad
 \delta(\infty)=0.
\]
Consequently $\delta=O(e^{-\kappa t}(1+t)^M)$ for some fixed $M$.
Taylor expansion of the averages in \eqref{averaged-four-ode} therefore
gives, for some $c>0$,
\begin{equation}\label{average-scalar-system}
 S'=R^{1/a}-\kappa S+O(e^{-ct}),\qquad
 R'=-S^q+O(e^{-ct}).
\end{equation}
The constants $c>0$ in exponential remainders may decrease in subsequent estimates.

Set $X=S^{q-a}$ and $\rho=R/(\kappa^aS^a)$.
Then $X\sim a\beta\kappa^a/t$ and $\rho\to1$.
Direct differentiation of \eqref{average-scalar-system} gives
\begin{align}
 \rho'&=-\kappa^{-a}X-a\kappa\rho(\rho^{1/a}-1)+O(e^{-ct}),
 \label{rho-eq}\\
 X'/X&=(\kappa/\beta)(\rho^{1/a}-1)+O(e^{-ct}).
 \label{X-eq}
\end{align}
Writing $Y=1-\rho=o(1)$, we obtain
\begin{equation}\label{YX-bootstrap}
 \begin{aligned}
 Y'+\kappa Y&=\kappa^{-a}X+O(Y^2)+O(e^{-ct}),\\[3pt]
 X'/X&=-\frac\kappa{a\beta}Y+O(Y^2)+O(e^{-ct}).
 \end{aligned}
\end{equation}
For large $t$, write the $O(Y^2)$ term as an $o(1)$ coefficient times
$Y$ and absorb it into $\kappa$. Variation of constants and $X=O(t^{-1})$
give $Y=O(t^{-1})$, then the second equation gives $X'=O(t^{-2})$.
Applying the first equation to $Y-\kappa^{-a-1}X$ yields
\[
 Y=\kappa^{-a-1}X+O(t^{-2}),\qquad
 X'/X=-\frac{X}{a\beta\kappa^a}+O(t^{-2}).
\]
Since $X\asymp t^{-1}$, integration gives
$X=(a\beta\kappa^a/t)(1+O(t^{-1}\log t))$.
To determine the logarithmic coefficient, set
\[
 \alpha_*=\kappa^{-a-1},\qquad
 d_*=(a\beta\kappa^a)^{-1},\qquad
 e_*=\frac{q(q-a)}{a^2\kappa^{2a+1}}.
\]
Expanding \eqref{rho-eq}--\eqref{X-eq} one order further gives
\[
 \begin{aligned}
 Y&=\alpha_*X+\frac{2q-a+1}{2a}\alpha_*^2X^2+O(t^{-3}),\\
 X'&=-d_*X^2-e_*X^3+O(t^{-4}).
 \end{aligned}
\]
Indeed, subtracting the first two terms on the right from $Y$ gives
$Z$ satisfying $Z'+\kappa Z=O(t^{-1}|Z|)+O(t^{-3})$.
Variation of constants gives $Z=O(t^{-3})$, and substitution in
\eqref{X-eq} gives the equation for $X$; no derivatives of the errors
are needed. In particular, $(1/X)'=d_*+e_*X+O(t^{-2})$.
Using $X=(d_*t)^{-1}+O(t^{-2}\log t)$, we obtain
\[
 \left(\frac1X-d_*t-\frac{e_*}{d_*}\log t\right)'
 =e_*\left(X-\frac1{d_*t}\right)+O(t^{-2})
 =O\left(\frac{\log t}{t^2}\right).
\]
The last bound is integrable. Hence, for some $h_\infty\in\R$,
\[
 \frac1X=d_*t+\frac{e_*}{d_*}\log t+h_\infty
             +O\left(\frac{\log t}{t}\right).
\]
Since $e_*/d_*^2=\gamma/\kappa$ and $S=X^\beta$, inversion and the
$\beta$th power give \eqref{averaged-common-phase}, with
$b_*=-\beta h_\infty/d_*$. The exponentially small difference $f-g$
does not affect this expansion. We also record the coarser estimates
\begin{equation}\label{common-profile-second}
 \begin{aligned}
 X&=\frac{a\beta\kappa^a}{t}
       \left(1-\frac\gamma\kappa\frac{\log t}{t}+O(t^{-1})\right),\\
 Y&=\frac{a\beta}{\kappa t}
       +O\left(\frac{\log t}{t^2}\right).
 \end{aligned}
\end{equation}

\smallskip
\noindent\emph{Step 2: the difference equation with a controlled remainder.}
For positive arguments define the positive divided differences
\[
 A_a(x,y)=\int_0^1\frac1a\big(y+s(x-y)\big)^{1/a-1}\dd s,
 \qquad
 B_q(x,y)=\int_0^1q\big(y+s(x-y)\big)^{q-1}\dd s.
\]
Integrating the equation for $\delta$ from $t$ to infinity gives
\[
 \delta(t)=-\frac12\int_t^\infty
   B_q(f(s),g(s))e^{-\kappa s}H(s)\dd s
   +O(e^{-(\kappa+c_0)t}).
\]
Also $H'=2e^{\kappa t}A_a(P,Q)\delta+e^{\kappa t}(E_1-E_2)$.
Fix $0<\eta<c_0$. Since $A_a(P,Q)$ grows at most polynomially,
we obtain
\begin{equation}\label{nonradial-H-volterra}
 H'(t)=-e^{\kappa t}A_a(P(t),Q(t))
       \int_t^\infty B_q(f(s),g(s))e^{-\kappa s}H(s)\dd s
       +O(e^{-\eta t}).
\end{equation}
Constants may depend on $\eta$. The independent remainder comes only
from the $E_i$ in \eqref{averaged-four-ode}; multiplication by
$e^{\kappa t}$ and polynomial factors preserves every decay exponent
strictly below $c_0$.

Since $P-Q$ and $f-g$ are exponentially small relative to their averages,
\eqref{common-profile-second} gives
\begin{equation}\label{divided-difference-product}
 \begin{split}
 \frac1\kappa A_a(P,Q)B_q(f,g)
 &=\frac{q}{a\kappa}R^{1/a-1}S^{q-1}(1+O(e^{-ct}))\\
 &=\frac{q}{a\kappa^a}X\rho^{1/a-1}(1+O(e^{-ct}))\\
 &=\frac\gamma t-\frac{\gamma^2}{\kappa}\frac{\log t}{t^2}
   +O(t^{-2}).
 \end{split}
\end{equation}
Let $M_H(t)=\sup_{s\ge t}|H(s)|$. By
\eqref{averaged-four-ode}--\eqref{common-profile-second}, $A_a(P,Q)$ and
$B_q(f,g)$ have logarithmic derivatives $O(t^{-1})$; they and their
reciprocals grow at most polynomially. Equation
\eqref{nonradial-H-volterra} gives
\[
 |H'(t)|\le Ct^{-1}M_H(t)+Ce^{-\eta t}.
\]
For $0\le y\le t/2$, this implies
$|H(t+y)-H(t)|\le Cy t^{-1}M_H(t)+Ce^{-\eta t}y$; also
\[
 B_q(f(t+y),g(t+y))
 =B_q(f(t),g(t))(1+O(y/t)).
\]
Insert these estimates in \eqref{nonradial-H-volterra} with $s=t+y$.
The integrals of $e^{-\kappa y}$ and $ye^{-\kappa y}$ are finite.
The part $y>t/2$, bounded by a polynomial in $t$ times
$e^{-\kappa t/2}M_H(t)$, is absorbed into $Ct^{-2}M_H(t)$.
It is not independent forcing: all coefficient and localization errors
multiply $H$ or $M_H$.
Together with \eqref{divided-difference-product}, this proves
\begin{equation}\label{nonradial-H-refined-ode}
 \left|H'+\left(\frac\gamma t-
          \frac{\gamma^2}{\kappa}\frac{\log t}{t^2}\right)H\right|
 \le Ct^{-2}M_H(t)+Ce^{-\eta t}.
\end{equation}
In particular,
\begin{equation}\label{nonradial-H-first-order}
 \left|H'+\frac\gamma tH\right|
 \le C\frac{\log t}{t^2}M_H(t)+Ce^{-\eta t}.
\end{equation}

\smallskip
\noindent\emph{Step 3: existence of the coefficient and the uniform error.}
For fixed large $T$, variation of constants in
\eqref{nonradial-H-first-order} gives
\begin{equation}\label{difference-bootstrap-integral}
 |H(t)|\le C_Tt^{-\gamma}
 +Ct^{-\gamma}\int_T^t s^{\gamma-2}(\log s)M_H(s)\dd s.
\end{equation}
The exponentially decaying forcing is included in $C_T$, since
$\int_T^\infty s^\gamma e^{-\eta s}\dd s<\infty$.
If $M_H(t)=O(t^{-\alpha}(\log t)^N)$, the integral gives
$H(t)=O(t^{-\alpha-1}(\log t)^{N+1})$ when $\alpha+1<\gamma$,
$H(t)=O(t^{-\gamma}(\log t)^{N+2})$ when $\alpha+1=\gamma$,
and $H(t)=O(t^{-\gamma})$ when $\alpha+1>\gamma$.
These bounds are eventually decreasing, so their suprema on $[t,\infty)$
give the same bounds for $M_H(t)$.
Starting from $M_H=O(1)$, finite iteration gives $M_H(t)=O(t^{-\gamma+\varepsilon})$ for any
fixed $0<\varepsilon<1$.
The exponential forcing contributes only $O(t^{-\gamma})$ at each step.
Multiplying \eqref{nonradial-H-first-order} by $t^\gamma$ gives
\[
 |(t^\gamma H)'|
 \le Ct^{-2+\varepsilon}\log t+Ct^\gamma e^{-\eta t}\in L^1(T,\infty).
\]
Hence $t^\gamma H\to d$ for some $d\in\R$, and in particular
$M_H=O(t^{-\gamma})$.  Substitution back into
\eqref{nonradial-H-first-order} and integration from $t$ to infinity give
\[
 |t^\gamma H(t)-d|
 \le C\int_t^\infty s^{-2}\log s\dd s+Ct^\gamma e^{-\eta t}
 =O(t^{-1}\log t).
\]
If $d=0$, set $K=t^\gamma H$ and
$M_K(t)=\sup_{s\ge t}|K(s)|$. Since $M_K$ is bounded and $K(\infty)=0$,
\eqref{nonradial-H-first-order} gives
\[
 \begin{aligned}
 |K'(t)|&\le C\frac{\log t}{t^2}M_K(t)+Ct^\gamma e^{-\eta t},\\
 |K(t)|&\le C\frac{1+\log t}{t}M_K(t)+Ct^\gamma e^{-\eta t}.
 \end{aligned}
\]
Taking the supremum over $[t,\infty)$ and absorbing the first term
for large $t$ gives $M_K(t)\le Ct^\gamma e^{-\eta t}$.
Thus $H=O(e^{-\eta t})$.

If $d\ne0$, then $H\sim dt^{-\gamma}$ and $M_H=O(|H|)$.
Dividing \eqref{nonradial-H-refined-ode} by $|H|$ yields
\[
 \frac{H'}H=-\frac\gamma t+
       \frac{\gamma^2}{\kappa}\frac{\log t}{t^2}+O(t^{-2}).
\]
Integration from $t$ to infinity, using $t^\gamma H\to d$, gives
\[
 H(t)=dt^{-\gamma}
       \left(1-\frac{\gamma^2}{\kappa}\frac{\log t}{t}
                    +O(t^{-1})\right).
\]
Together with the case $d=0$, this proves
\eqref{averaged-difference-asymptotic} and the asserted exponential decay.
\end{proof}

\begin{proof}[Proof of \eqref{exact-common-phase} and
\eqref{nonradial-second-order-limit} in Theorem~\ref{mainexact}]
To apply Lemma~\ref{averaged-difference-lemma}, for either operator set $U=e^{-\kappa t}w_1(e^{-t}\theta)$,
$V=e^{-\kappa t}w_2(e^{-t}\theta)$ and
\[
 f=\langle U\rangle,\quad g=\langle V\rangle,\qquad
 P=\frac{\mu_1(B_{e^{-t}})}{\omega d_0},\quad
 Q=\frac{\mu_2(B_{e^{-t}})}{\omega d_0}.
\]
Both point masses vanish. For the $p$-Laplace operator, $d_0=1$ and
\begin{gather}
 P=\left\langle
 \big((\kappa U+U_t)^2+|\nabla_\theta U|^2\big)^{(p-2)/2}
 (\kappa U+U_t)\right\rangle,\label{nonradial-flux-identity}\\[3pt]
 P'=-\langle V^q\rangle,\qquad
 Q'=-\langle U^q\rangle.\label{nonradial-mass-identities}
\end{gather}
For the Hessian operator, write $d=d_0$ and $u=-w_1$, $v=-w_2$.
The Newton identity
$k\sigma_k(D^2u)=\operatorname{div}(T_{k-1}(D^2u)\nabla u)$ gives
\begin{equation}\label{hessian-averaged-flux}
 P=\frac1{kd}\left\langle
 \left[T_{k-1}(\mathsf H_\kappa[U])
 \begin{pmatrix}\kappa U+U_t\\-\nabla_\theta U\end{pmatrix}
 \right]_r\right\rangle,
 \qquad P'=-d^{-1}\langle V^q\rangle,
\end{equation}
and the analogous identities for $Q$. The powers of $r$ cancel since $k(\kappa+2)=n$. Apply the divergence
identity on annuli: Lemma~\ref{hessian-annular-smooth} gives inner flux
$O(t^{-k\beta})\to0$, so the zero-mass extension yields the ball-mass
identity.
For radial $U=f$, its Newton expression is $(\kappa f+f')^k$.

The leading limits and smooth annular convergence give, in both cases,
\[
 f,g\sim C_0t^{-\beta},\qquad
 \kappa f+f',\ \kappa g+g'\sim\kappa C_0t^{-\beta},\qquad
 P,Q\sim(\kappa C_0)^a t^{-a\beta}.
\]
Fix $0<c_0<\eta_0<1$. Lemma~\ref{angular-decay} gives errors of order
$O(e^{-(\kappa+\eta_0)t})$ for $U-f,V-g$ and their first two
derivatives. The flux formulas, the mass identities, and inversion of
the positive flux powers consequently yield
\begin{equation}\label{hessian-averaged-four-ode}
 \begin{gathered}
 \begin{alignedat}{2}
 f'&=P^{1/a}-\kappa f+E_1,&\qquad g'&=Q^{1/a}-\kappa g+E_2,\\[3pt]
 P'&=-d_0^{-1}g^q+E_3,& Q'&=-d_0^{-1}f^q+E_4,
 \end{alignedat}\\[4pt]
 |E_i|\le Ce^{-(\kappa+c_0)t}
 \end{gathered}
\end{equation}
with the prescribed $c_0$. The Newton tensor is polynomial in the
Hessian, and all derivatives of the other powers and the inverse flux
power are bounded by fixed powers of $t$. The strict inequality
$c_0<\eta_0$ absorbs these factors; this also handles $p<2$ or $q<1$.

We replace $(f,g,P,Q)$ by
$(d_0^{-\beta}f,d_0^{-\beta}g,d_0^{-a\beta}P,d_0^{-a\beta}Q)$
and omit the tildes. Since $\beta(q-a)=1$, the rescaled functions satisfy
\eqref{averaged-four-ode}, with the same exponential gain $c_0>0$.
Since $d_0^{-\beta}C_0=(a\beta\kappa^a)^\beta$, $f,g$ have the required
asymptotics. The flux and mass identities ensure positivity and $C^1$
regularity of the averaged functions and masses on a sufficiently late tail. The bounded-difference
estimates imply $H=e^{\kappa t}(f-g)=O(1)$.
Lemma~\ref{averaged-difference-lemma}, with the factor $d_0^\beta$ and
angular error restored, gives \eqref{exact-common-phase} with the same
$b_*$ for both components. It also gives
\[
 H(t)=\widetilde c\,t^{-\gamma}
       -\frac{\widetilde c\gamma^2}{\kappa}t^{-\gamma-1}\log t
       +O(t^{-\gamma-1})
\]
for some $\widetilde c\in\R$. If $\widetilde c=0$, then
$H=O(e^{-\eta t})$ for every $0<\eta<c_0$. Finally
Lemma~\ref{angular-decay} implies, uniformly in $\theta$,
\[
 w_1(e^{-t}\theta)-w_2(e^{-t}\theta)
 =d_0^\beta H(t)+O(e^{-\eta_0t}).
\]
Taking $c=d_0^\beta\widetilde c$ proves
\eqref{nonradial-second-order-limit}, with $\gamma=n/p$ or
$\gamma=n/(2k)$ according to the operator. In particular the
coefficient is a single real number, independent of direction.
Given any $0<\delta<1$, choose
$\delta<c_0<\eta_0<1$. If $c=0$, the two exponential bounds then give
$w_1(r\theta)-w_2(r\theta)=O(r^\delta)$, uniformly in $\theta$.
\end{proof}

\begin{remark}\label{difference-gradient}
Under the nonremovable alternative of Theorem~\ref{mainexact}, let
$h=w_1-w_2$ and let $c$ be the coefficient in
\eqref{nonradial-second-order-limit}. Then, uniformly in $\theta$,
\[
 rL(r)^{\gamma+1}\nabla h(r\theta)\longrightarrow c\gamma\theta.
\]
With $h(0)=0$, the difference is Dini continuous at $0$ and belongs to
$W^{1,n}_{\loc}(B_1)$. If $c\ne0$, it belongs to no $W^{1,s}(B_R)$
with $s>n$ and $0<R<1$.

Equation~\eqref{nonradial-H-first-order} and
$M_H(t)=O(t^{-\gamma})$ give
$H'(t)=-\gamma\widetilde c\,t^{-\gamma-1}
+O(t^{-\gamma-2}\log t)$ in the preceding proof.
Restoring the factor $d_0^\beta$ and using the cylinder derivative
estimates in \eqref{absolute-angular-decay} proves the gradient limit
and the bound $|\nabla h|=O(r^{-1}L^{-\gamma-1})$.
The regularity assertions follow by radial integration, since
$h=O(L^{-\gamma})$ and $\gamma>1$; when $c\ne0$, the gradient limit
also gives the matching lower bound.
\end{remark}

\subsection{Radial realization and rigidity}

\begin{proposition}\label{nonzero-radial-coefficient}
For either operator, \eqref{exact-symmetric} has a smooth radial
singular solution with $c\ne0$ in \eqref{nonradial-second-order-limit}.
Among radial singular solutions, $c=0$ if and only if $w_1=w_2$.
\end{proposition}

\begin{proof}[Proof of Proposition~\ref{nonzero-radial-coefficient}]
We first construct a nondiagonal positive solution on a tail of the
normalized exact system
\begin{equation}\label{exact-radial-four-ode}
 f'=P^{1/a}-\kappa f,\quad g'=Q^{1/a}-\kappa g,
 \qquad P'=-g^q,\quad Q'=-f^q.
\end{equation}
Only $a,\kappa>0$ and $q>a$ are used in this construction. Introduce
\[
 x=f^{q-a},\qquad h=\log(g/f),\qquad
 r=\frac{P}{\kappa^af^a},\qquad
 s=\frac{Q}{\kappa^ag^a}.
\]
The resulting autonomous system is
\begin{equation}\label{radial-smooth-system}
 \begin{aligned}
 x'&=(\kappa/\beta)x(r^{1/a}-1),&
 h'&=\kappa(s^{1/a}-r^{1/a}),\\
 r'&=-\kappa^{-a}e^{qh}x-a\kappa r(r^{1/a}-1),&
 s'&=-\kappa^{-a}e^{-ah}x-a\kappa s(s^{1/a}-1).
 \end{aligned}
\end{equation}
The vector field is smooth near $(x,h,r,s)=(0,0,1,1)$, allowing
either sign of $x$ in this local extension.

\smallskip
\noindent\emph{A diagonal slow trajectory.}
The subspace $h=0$, $r=s$ is invariant. Set
$b_0=\kappa^{-a-1}$ and choose $0<b_-<b_0<b_+$.
For sufficiently small $\varepsilon_0>0$, the region
\[
 0<x\le\varepsilon_0,\qquad 1-b_+x\le r\le1-b_-x
\]
is positively invariant for the diagonal two-dimensional system.
Indeed, on $r=1-bx$,
\[
 (r-1+bx)'=(\kappa b-\kappa^{-a})x+O(x^2),
\]
which points inward on both boundary curves; in the region,
$-C_2x^2\le x'\le-C_1x^2$ for some $C_2>C_1>0$.
Starting inside this region at a time $T$ gives a trajectory
$(x_0,0,r_0,r_0)$ on $[T,\infty)$ with
\begin{equation}\label{diagonal-polynomial-bounds}
 \frac{x_0(T)}{1+C_2x_0(T)(t-T)}
 \le x_0(t)\le
 \frac{x_0(T)}{1+C_1x_0(T)(t-T)}.
\end{equation}
Expanding the diagonal equation gives
\[
 (r_0-1+b_0x_0)'=-\kappa(r_0-1+b_0x_0)+O(x_0^2).
\]
Variation of constants and \eqref{diagonal-polynomial-bounds} imply
$r_0-1+b_0x_0=O(t^{-2})$. Hence
\[
 x_0'=-\frac{x_0^2}{a\beta\kappa^a}+O(t^{-3}),\qquad
 x_0(t)=\frac{a\beta\kappa^a}{t}
       \left(1+O\left(\frac{\log t}{t}\right)\right).
\]
In particular $f_0=x_0^\beta\sim(a\beta\kappa^a)^\beta t^{-\beta}$.

\smallskip
\noindent\emph{An exact nondiagonal trajectory.}
Write $z=(x,h,r-1,s-1)$ and denote the vector field in
\eqref{radial-smooth-system} by $\mathcal F$.
Its derivative at $z=0$ is
\[
 A(x,h,u,v)=
 \left(0,\frac\kappa a(v-u),-\kappa^{-a}x-\kappa u,
                         -\kappa^{-a}x-\kappa v\right).
\]
It has the complementary eigenspaces
\[
 E_c=\{(x,h,-b_0x,-b_0x):x,h\in\mathbb R\},\qquad
 E_s=\{(0,(u-v)/a,u,v):u,v\in\mathbb R\},
\]
with eigenvalues $0$ and $-\kappa$, respectively. Let $\Pi_c,\Pi_s$
be their projections and set
$z_0=(x_0,0,r_0-1,r_0-1)$.
For a difference $v=z-z_0$, write
\[
 N(t,v)=\mathcal F(z_0(t)+v)-\mathcal F(z_0(t))-Av.
\]
Take $x_0(T)$ and the allowed radius for $v$ small. Smoothness gives,
uniformly for $t\ge T$,
\[
 N(t,0)=0,\qquad |N(t,v)-N(t,w)|\le\varepsilon|v-w|,
\]
where $\varepsilon$ can be arbitrarily small.
Fix $\lambda=\kappa/2$. Given $\xi\in E_s$, solve the mixed
Lyapunov--Perron equations
\begin{equation}\label{radial-stable-fiber-integrals}
 \begin{aligned}
 \Pi_sv(t)&=e^{-\kappa(t-T)}\xi+
       \int_T^t e^{-\kappa(t-\tau)}\Pi_sN(\tau,v(\tau))\dd\tau,\\
 \Pi_cv(t)&=-\int_t^\infty\Pi_cN(\tau,v(\tau))\dd\tau.
 \end{aligned}
\end{equation}
In the norm $\|v\|_\lambda=\sup_{t\ge T}e^{\lambda(t-T)}|v(t)|$,
the nonlinear part has Lipschitz constant at most
\[
 C_A\varepsilon\left((\kappa-\lambda)^{-1}+\lambda^{-1}\right),
\]
where $C_A$ depends only on the two fixed projections. Choose $\varepsilon$ so that this is less than $1/2$. For small $|\xi|$,
the map in \eqref{radial-stable-fiber-integrals} is a contraction of
the ball $\|v\|_\lambda\le2C_A|\xi|$ into itself. Its fixed point is
differentiable and solves $v'=Av+N(t,v)$, so $z_0+v$ is an exact
solution. It satisfies
\begin{equation}\label{radial-stable-fiber-bound}
 |v(t)|\le C|\xi|e^{-\lambda(t-T)},\qquad
 v(T)=\xi+O(\varepsilon|\xi|).
\end{equation}
Choose $\xi=\alpha(0,1,a/2,-a/2)$, with $\alpha\ne0$ sufficiently
small. The $h$-coordinate at $T$ is $\alpha+O(\varepsilon|\alpha|)$,
and is nonzero after decreasing $\varepsilon$.
The lower bound in \eqref{diagonal-polynomial-bounds}, together with
\eqref{radial-stable-fiber-bound}, shows that $x(t)>0$ for every
$t\ge T$ when $|\alpha|\le c_1 x_0(T)$, with a sufficiently small
fixed $c_1>0$: the ratio of the perturbation bound to $x_0(t)$ is at
most
\[
 C\frac{|\alpha|}{x_0(T)}
       (1+C_2x_0(T)(t-T))e^{-\lambda(t-T)}.
\]
The last factor is uniformly bounded since $x_0(T)\le\varepsilon_0$.
Also $r,s>0$ throughout the tail. Thus
\[
 f=x^\beta,\quad g=e^hf,\quad P=\kappa^arf^a,\quad Q=\kappa^asg^a
\]
are positive and solve \eqref{exact-radial-four-ode}. The exponential
tracking and \eqref{diagonal-polynomial-bounds} give
$f,g\sim(a\beta\kappa^a)^\beta t^{-\beta}$; the pair is nondiagonal
because $h(T)\ne0$.

\smallskip
\noindent\emph{Return to the PDE.}
Restore the normalization by setting
$\bar f=d_0^\beta f$, $\bar g=d_0^\beta g$ and
$\bar P=d_0^{a\beta}P$, $\bar Q=d_0^{a\beta}Q$.
Then $\bar P'=-d_0^{-1}\bar g^q$ and
$\bar Q'=-d_0^{-1}\bar f^q$, while the two filter equations retain
their form. Consequently
\[
 w_1(r)=r^{-\kappa}\bar f(L(r)),\qquad
 w_2(r)=r^{-\kappa}\bar g(L(r))
\]
solve \eqref{exact-symmetric} on $0<r<e^{-T}$ and have the common
leading coefficient $C_0$. The fluxes tend to zero at the origin.
In the Hessian case, the positive radial derivatives and sources give
strict $k$-admissibility by Lemma~\ref{radial-admissibility}. The critical dilation
$\widehat w_i(x)=R^\kappa w_i(Rx)$, $R=e^{-T}$, yields a pair on
$B_1\setminus\{0\}$ and preserves exactness, singularity and
nondiagonality.

\smallskip
\noindent\emph{Vanishing of the radial coefficient is rigid.}
For any exact radial singular pair, all the $E_i$ in
\eqref{averaged-four-ode} vanish. In the normalized variables,
$H=e^{\kappa t}(f-g)$ therefore satisfies the exact homogeneous identity
\[
 H'(t)=-e^{\kappa t}A_a(P(t),Q(t))
       \int_t^\infty B_q(f(\tau),g(\tau))e^{-\kappa\tau}H(\tau)\dd\tau.
\]
The coefficient estimates and localization used to obtain
\eqref{nonradial-H-first-order} now give the homogeneous estimate
\[
 \left|H'+\frac\gamma tH\right|
 \le C\frac{\log t}{t^2}M_H(t).
\]
There is no forcing: every localization error, including the far tail,
is bounded by a coefficient times $M_H(t)$.
By Theorem~\ref{mainexact}, $K(t)=t^\gamma H(t)$ has a finite limit.
If $c=0$, this limit is zero. With
$M_K(t)=\sup_{\tau\ge t}|K(\tau)|<\infty$, we have
$M_H(t)\le t^{-\gamma}M_K(t)$ and hence
\[
 |K'(t)|\le C\frac{\log t}{t^2}M_K(t),\qquad
 |K(t)|\le C\frac{1+\log t}{t}M_K(t).
\]
Since $(1+\log \tau)/\tau$ decreases for $\tau\ge1$, taking the supremum
over $[t,\infty)$ gives a coefficient less than one for large $t$,
so $M_K(t)=0$. Therefore $f=g$ and
$P=Q$ on a tail. Uniqueness for the positive smooth ODE propagates
this identity throughout the radial interval. The converse is immediate.
The nondiagonal pair constructed above consequently has $c\ne0$,
which realizes the leading exponent and the nonzero logarithmic correction
in~\eqref{nonradial-second-order-limit}.
\end{proof}

\subsection{Coupling coefficients with averaged limits}

For bounded nonnegative measurable coefficients $A,B$ on a punctured
neighborhood of $0$, write, for $t\ge T_0$,
\[
 \bar A(t)=\av_{\Sph}A(e^{-t}\theta),\qquad
 \bar B(t)=\av_{\Sph}B(e^{-t}\theta).
\]

\begin{corollary}\label{mainvariable}
Suppose
\begin{equation}\label{coefficient-cesaro-limits}
 \frac1T\int_{T_0}^T\bar A(t)\dd t\to A_0>0,\qquad
 \frac1T\int_{T_0}^T\bar B(t)\dd t\to B_0>0.
\end{equation}
Every positive solution of
\begin{equation}\label{exact-variable}
 \mathcal Q[w_1]=A(x)w_2^q,\qquad
 \mathcal Q[w_2]=B(x)w_1^q
\end{equation}
either has two removable singularities or satisfies, uniformly in $\theta$,
\begin{equation}\label{exact-limits}
 r^\kappa L(r)^\beta(w_1,w_2)(r\theta)
 \to C_0\bigl(A_0^{-a/E}B_0^{-q/E},A_0^{-q/E}B_0^{-a/E}\bigr),
 \qquad E=q^2-a^2.
\end{equation}
\end{corollary}

Positive pointwise limits imply these assumptions; for coefficients
depending only on direction, $A_0,B_0$ are their spherical means. Only leading terms,
not \eqref{nonradial-second-order-limit}, are asserted under these
averaged hypotheses. For example,
\[
 A(r)=1+\sin L(r),\qquad B(r)=1-\sin L(r),
\]
have no pointwise limits or positive lower bounds near the origin,
but have logarithmic Ces\`aro means $1$. Thus every singular solution
with these coefficients has leading coefficients $C_0$.

\begin{proof}
Write $f_1=Aw_2^q$ and $f_2=Bw_1^q$. Boundedness of $A,B$ gives
$f_1+f_2\le C(w_1+w_2)^q$, so Lemma~\ref{fundamental-vector} applies.
If $\mu_1(\{0\})>0$, its uniform coefficient limit gives
$w_1(r\theta)\ge c_1r^{-\kappa}$ for all small $r$. Hence
\[
 \int_{B_R} f_2\dd x
 \ge c_2\int_0^R\av_{\Sph}B(r\theta)\frac{\dd r}{r}
 =\infty.
\]
The divergence follows from the positive Ces\`aro limit of
$\bar B$ and contradicts local finiteness of the source measure.
Interchanging the components excludes the other point mass.
Lemma~\ref{radial-comparison} now gives the bounded-or-divergent
alternative for the sum. If bounded, both sources are bounded and both
singularities are removable as in the proof of Theorem~\ref{mainjoint}.

Assume the sum diverges. Set
\begin{equation}\label{coefficient-normalization}
 c_1=A_0^{a/E}B_0^{q/E},\qquad
 c_2=A_0^{q/E}B_0^{a/E},\qquad E=q^2-a^2.
\end{equation}
By homogeneity, $(c_1w_1,c_2w_2)$ solves \eqref{exact-variable} with
coefficients $A/A_0,B/B_0$. We may therefore assume $A_0=B_0=1$.
Let $P,Q$ be the normalized ball masses and $J=P+Q$. The uniform
expansion \eqref{mass-expansion}, followed by spherical integration of
the sources, gives, for a.e. large $t$,
\begin{equation}\label{averaged-coefficient-mass}
 \begin{split}
 -P'&=c_*\bar A(t)\bigl(Q^m+\varepsilon_1(t)J^m\bigr),\\
 -Q'&=c_*\bar B(t)\bigl(P^m+\varepsilon_2(t)J^m\bigr),
 \end{split}
 \qquad c_*=(d_0\kappa^q)^{-1}.
\end{equation}
Here $\varepsilon_i(t)\to0$. The mass expansion is uniform in direction,
even when one mass is small. Nonnegativity preserves the spherical
average as a factor in each error; if it vanishes, set $\varepsilon_i=0$.

We claim that
\begin{equation}\label{averaged-coefficient-mass-rate}
 J(t)\asymp t^{-\rho},\qquad \rho=(m-1)^{-1}=a\beta.
\end{equation}
Here $J>0$ on a tail, since the radial reconstruction is unbounded.
Boundedness of the coefficients gives $0\le-J'\le CJ^m$, and hence
$J(t)\ge ct^{-\rho}$. For the reverse bound,
\eqref{coefficient-cesaro-limits} gives, for all sufficiently large $t$,
\[
 \int_t^{2t}\bar A(s)\dd s\ge ct,\qquad
 \int_t^{2t}\bar B(s)\dd s\ge ct.
\]
Fix $0<\vartheta<2^{-\rho}$. If $J(2t)\ge\vartheta J(t)$, then
one of $P(2t),Q(2t)$ is at least $\vartheta J(t)/2$.
Suppose it is $P(2t)$. Monotonicity gives
$P(s)\ge P(2t)\ge\vartheta J(t)/2$ and $J(s)\le J(t)$ on $[t,2t]$,
so the second error in \eqref{averaged-coefficient-mass} can be
absorbed for large $t$. Thus
\[
 J(t)\ge Q(t)-Q(2t)
 \ge\frac{c_*}{2}P(2t)^m\int_t^{2t}\bar B(s)\dd s
 \ge c tJ(t)^m.
\]
The other case is identical, and both give $J(t)\le Ct^{-\rho}$.
If $J(2t)<\vartheta J(t)$ instead, we already have a fixed contraction.
Consequently, for every large $t$,
\[
 J(2t)\le\vartheta J(t)+Ct^{-\rho}.
\]
Since $2^\rho\vartheta<1$, iteration on dyadic times and monotonicity
between them give the upper bound in \eqref{averaged-coefficient-mass-rate}.

To compare the normalized masses, define
\[
 D(t)=P(t)^{m+1}-Q(t)^{m+1},\qquad
 W(t)=P(t)^mQ(t)^m.
\]
The masses are nonnegative and nonincreasing, so $W$ is locally
absolutely continuous and nonincreasing. Equation
\eqref{averaged-coefficient-mass} gives
\begin{equation}\label{averaged-coefficient-difference}
 D'=(m+1)c_*\bigl(\bar B-\bar A\bigr)W+o(J^{2m}).
\end{equation}
By \eqref{coefficient-cesaro-limits},
$G(t)=\int_{T_0}^t(\bar B-\bar A)\dd s=o(t)$.
Also $W(t)=O(t^{-2m\rho})$ and $2m\rho>1$.
Integration by parts is therefore valid on the tail and yields
\[
 \int_t^\infty(\bar B-\bar A)W\dd s
 =-G(t)W(t)-\int_t^\infty G(s)W'(s)\dd s
 =o(t^{1-2m\rho}).
\]
Indeed, let
$\varepsilon_t=\sup_{s\ge t}|G(s)|/s\to0$. The absolute value is at
most
\[
 \varepsilon_t\left(tW(t)+\int_t^\infty s(-W'(s))\dd s\right)
 =\varepsilon_t\left(2tW(t)+\int_t^\infty W(s)\dd s\right)
 =o(t^{1-2m\rho}).
\]
The tail integral of $o(J^{2m})$ has the same order, by
\eqref{averaged-coefficient-mass-rate}. Since $D(\infty)=0$ and
$1-2m\rho=-(m+1)\rho$, integrating
\eqref{averaged-coefficient-difference} gives
$D=o(J^{m+1})$. Consequently
\begin{equation}\label{averaged-coefficient-ratio}
 P/J\to\tfrac12,\qquad Q/J\to\tfrac12.
\end{equation}
This uses cancellation against a monotone weight, without assuming
pointwise coefficient limits.

We finally sum \eqref{averaged-coefficient-mass} and use
\eqref{averaged-coefficient-ratio} to obtain
\begin{equation}\label{variable-total-mass}
 (J^{1-m})'=(m-1)c_*2^{-m}
                  (\bar A+\bar B)+o(1).
\end{equation}
Integrating and using the two normalized Ces\`aro limits yields
$J^{1-m}\sim(m-1)c_*2^{1-m}t$, exactly as in
\eqref{exact-mass-rate}. Substitution in \eqref{mass-expansion} gives
$t^\beta e^{-\kappa t}w_i(e^{-t}\theta)\to C_0$ uniformly in
$\theta$ for the normalized pair. Dividing by the constants in
\eqref{coefficient-normalization} proves \eqref{exact-limits}.
\end{proof}

\begin{corollary}[A common coefficient]\label{common-coefficient-clock}
Suppose $0\le A\le C$ and
$\mathcal T(t)=\int_{T_0}^t\bar A(s)\dd s\to\infty$.
Every positive solution of \eqref{exact-variable} with $B=A$ either has
two removable singularities or satisfies, uniformly in $\theta$,
\[
 r^\kappa\mathcal T(L(r))^\beta(w_1,w_2)(r\theta)
 \longrightarrow C_0(1,1).
\]
\end{corollary}
\begin{proof}
The point-mass exclusion and bounded-or-divergent alternative in the
proof of Corollary~\ref{mainvariable} use only
$\int_{T_0}^\infty\bar A=\infty$ and remain valid.
On the divergent branch, use the same masses $P,Q,J$ and set
$D=P^{m+1}-Q^{m+1}$. Equation~\eqref{averaged-coefficient-mass} gives
$-J'\asymp\bar A J^m$ and $D'=\bar A\,o(J^{2m})$.
Thus $|D'|\le\varepsilon(t)(-J')J^m$, with $\varepsilon(t)\downarrow0$;
integration from $t$ to infinity yields $D=o(J^{m+1})$.
Consequently $P/J,Q/J\to1/2$, and summing the mass equations gives
\[
 (J^{1-m})'=(m-1)(d_0\kappa^q)^{-1}2^{1-m}
              \bar A(t)(1+o(1)).
\]
Integration against the divergent clock $\mathcal T$, followed by
\eqref{mass-expansion}, proves the claim.
\end{proof}
For example, $A(r)=L(r)^{-\sigma}$ with $0<\sigma<1$ gives
$w_i\sim C_0(1-\sigma)^\beta r^{-\kappa}L(r)^{-(1-\sigma)\beta}$,
whereas $A(r)=L(r)^{-1}$ gives
$w_i\sim C_0r^{-\kappa}(\log L(r))^{-\beta}$ on the singular branch.

\section{The nonsymmetric Serrin critical curve}\label{sec:nonsymmetric}

Here $q_1,q_2,\eta,D$ satisfy
\eqref{nonsymmetric-parameters}; in particular,
$0<\eta<\kappa$ and $D=\ell(a+q_1)/\kappa>0$.
We use the identities
\begin{equation}\label{nonsymmetric-identities}
 q_1\eta=n,\qquad
 q_2\kappa=\ell+a\eta<n,\qquad
 \frac{\eta}{\kappa}=\frac{a+q_2}{a+q_1}.
\end{equation}
The canonical measures from Section~\ref{sec:local} have densities
$w_2^{q_1}$ and $w_1^{q_2}$ away from $0$, both locally integrable across
the origin. A concave-power comparison reduces the first component to a
scalar critical inequality; its mass asymptotics determine the second
through a radial source problem.

\subsection{Natural scales and vanishing masses}

\begin{lemma}\label{nonsymmetric-upper}
Every solution of \eqref{nonsymmetric-system} satisfies
\begin{equation}\label{nonsymmetric-smallness}
 \sup_{\theta\in\Sph}
 \bigl(r^\kappa w_1(r\theta)+r^\eta w_2(r\theta)\bigr)
 \to0\qquad(r\to0).
\end{equation}
\end{lemma}
\begin{proof}
Set $M=w_1^{1/\kappa}+w_2^{1/\eta}$. Suppose
$|x|M(x)$ is unbounded near zero. As in the proof of
Lemma~\ref{fundamental-vector}, apply the doubling lemma
\cite[Lemma~5.1]{PQS} on shrinking punctured balls. We obtain
$x_j\to0$ and $\lambda_j=M(x_j)^{-1}=o(|x_j|)$ such that, on
balls $B_j$ avoiding the puncture, the functions
\[
 U_j(y)=\lambda_j^\kappa w_1(x_j+\lambda_j y),\qquad
 V_j(y)=\lambda_j^\eta w_2(x_j+\lambda_j y)
\]
are uniformly bounded and
$U_j(0)^{1/\kappa}+V_j(0)^{1/\eta}=1$.
By \eqref{nonsymmetric-identities}, the scaled equations are
\[
 \mathcal Q[U_j]=V_j^{q_1},\qquad
 \mathcal Q[V_j]=U_j^{q_2}.
\]
For every fixed $R>0$,
\[
 \int_{B_R}V_j^{q_1}\dd y
 =\int_{x_j+\lambda_jB_R}w_2^{q_1}\dd x
 \to0.
\]
Both scaled densities are locally bounded. By Section~\ref{sec:local},
extraction gives limits $U,V$ locally in $C^1$ for the $p$-Laplacian
and locally uniformly with weak measure convergence for the Hessian
operator.
The integral gives $V=0$, and the second equation forces
$U^{q_2}=0$, contradicting the normalization.
Hence $|x|M(x)$ is bounded, giving
$w_1=O(r^{-\kappa})$ and $w_2=O(r^{-\eta})$.

For any decreasing sequence $r_j\to0$, rescale on fixed annuli by
\[
 U_j(y)=r_j^\kappa w_1(r_jy),\qquad
 V_j(y)=r_j^\eta w_2(r_jy).
\]
The scale bounds and local regularity give subsequential convergence
on compact annuli. Since
\[
 \int_K V_j^{q_1}\dd y
 =\int_{r_jK}w_2^{q_1}\dd x\to0,
\]
the limit of $V_j$ is zero, and the second equation gives the same for
$U_j$. This proves \eqref{nonsymmetric-smallness} for every sequence.
A positive mass at $0$ would force the corresponding component to be at
least $c|x|^{-\kappa}$ by \eqref{common-potential}, a contradiction.
\end{proof}

\subsection{Concave-power comparison across the puncture}

\begin{lemma}
\label{nonsymmetric-crosspoint}
Let $s,t>0$ be supersolutions for either operator in the prescribed
solution class, with $s+t=O(|x|^{-\eta})$ for some $0<\eta<\kappa$.
If $M_0\ge0$, $t\le s+M_0$ on an outer sphere, and
\[
 \mathcal Q[s]\ge\mathcal Q[t]
 \quad\text{on }\{t>s+M_0\},
\]
then $t\le s+M_0$ inside that sphere.
\end{lemma}
\begin{proof}
For $\mathcal Q=-\Delta_p$, a nonnegative weak
supersolution $h=O(r^{-\eta})$ obeys the annular Caccioppoli bound
\begin{equation}\label{nonsymmetric-caccioppoli}
 \int_{\{r<|x|<2r\}}|\nabla h|^p\dd x
 \le Cr^{n-p-p\eta}.
\end{equation}
Indeed, on a larger concentric annulus take
$H=\sup h$ and test the supersolution inequality with
$(H-h)\chi^p$, where $\chi$ is an annular cutoff. By H\"older's
inequality, for an inner cutoff $\chi_r$ with
$|\nabla\chi_r|\le C/r$,
\begin{equation}\label{nonsymmetric-cutoff}
 \int_{\{r<|x|<2r\}}
 (|\nabla s|^{p-1}+|\nabla t|^{p-1})
 |\nabla\chi_r|\dd x
 \le Cr^{n-p-(p-1)\eta}
 =Cr^{a(\kappa-\eta)}\to0.
\end{equation}
Subtract the weak inequalities on $\{t>s+M_0\}$ and test with
an inner cutoff times
$T_b((t-s-M_0)_+)$, where $T_b(h)=\min\{h,b\}$.
The cutoff error is $O(b)$ times the left-hand side of
\eqref{nonsymmetric-cutoff}. As the inner radius tends to zero,
strict monotonicity of $|\zeta|^{p-2}\zeta$ gives
$(t-s-M_0)_+=0$.

For the Hessian operator, $-r^{-\kappa}$ is $k$-convex and has
zero $k$-Hessian away from zero. G\aa rding's inequality and
cone monotonicity imply that $-(s+\varepsilon r^{-\kappa})$
is $k$-convex and
\[
 \mathcal Q[s+\varepsilon r^{-\kappa}]\ge\mathcal Q[s].
\]
For fixed $\varepsilon>0$, the $r^{-\kappa}$ term dominates $t$ on a
small inner sphere. Compare $t$ with $s+M_0+\varepsilon r^{-\kappa}$
on the annulus: where $t$ is larger, its source is no greater than that
of the comparison function. The $k$-Hessian comparison principle applies. Letting the inner radius tend to zero and then $\varepsilon\to0$,
we obtain $t\le s+M_0$ throughout the punctured ball, completing the proof.
\end{proof}

The semilinear antecedent of our concave-power comparison is
\cite[Theorem~1.2(i)]{BV2000}. Set
\begin{equation}\label{nonsymmetric-concave-power}
 \xi=\frac\eta\kappa=\frac{a+q_2}{a+q_1}\in(0,1),
 \qquad K=\xi^{-a/(a+q_1)},\qquad P=K w_1^\xi.
\end{equation}
The $p$-Laplace chain rule and, respectively, the Hessian identity
\[
 D^2(-P)=K\xi w_1^{\xi-1}D^2(-w_1)
 +K\xi(1-\xi)w_1^{\xi-2}\nabla w_1\otimes\nabla w_1
\]
give
\begin{equation}\label{nonsymmetric-concave-ineq}
 \mathcal Q[P]\ge(K\xi)^a
 w_1^{a(\xi-1)}\mathcal Q[w_1]
 =(K\xi)^a w_1^{a(\xi-1)}w_2^{q_1}.
\end{equation}
Thus $P$ is a positive supersolution, and
Lemma~\ref{nonsymmetric-upper} gives
$P+w_2=O(r^{-\eta})$. On a small outer sphere $\partial B_R$, choose
\[
 M_0>\max\{0,\max_{\partial B_R}(w_2-P)\}.
\]
On $\{w_2>P+M_0\}$,
\[
 \mathcal Q[P]\ge K^{a+q_1}\xi^a
 w_1^{a(\xi-1)+q_1\xi}
 =w_1^{q_2}=\mathcal Q[w_2].
\]
Lemma~\ref{nonsymmetric-crosspoint} therefore proves
\begin{equation}\label{nonsymmetric-power-bound}
 w_2\le K w_1^\xi+M_0.
\end{equation}
Since $q_1\xi=q_S$, translation invariance and
\eqref{nonsymmetric-power-bound} give, with $h=w_1+1$,
\begin{equation}\label{nonsymmetric-scalar-ineq}
 0\le\mathcal Q[h]=w_2^{q_1}\le C h^{q_S},
 \qquad r^\kappa h(r\theta)\to0
 \quad\text{uniformly in }\theta.
\end{equation}
The scalar instances of Lemmas~\ref{radial-comparison}
and~\ref{radial-ode} and Proposition~\ref{mass-approx} apply to $h$.
If $h$ is bounded, \eqref{nonsymmetric-power-bound} bounds $w_2$,
and both singularities are removable. Otherwise the reconstruction
$\zeta$ from $\mu_1$ in \eqref{radial-z} satisfies
\begin{equation}\label{nonsymmetric-scalar-mass}
 \begin{gathered}
 \frac{w_1(r\theta)}{\zeta(r)}\to1,
 \qquad\zeta(r)\to\infty,
 \qquad-\frac{r\zeta'(r)}{\zeta(r)}\to\kappa,\\
 Z(t):=r^\kappa\zeta(r)
 \sim\frac1\kappa J(t)^{1/a},\qquad
 J(t):=\frac{\mu_1(B_{e^{-t}})}{\omega d_0}\to0,
 \quad t=L(r),
 \end{gathered}
\end{equation}
uniformly in direction.

\subsection{Recovery of the second component and exact constants}

On the unbounded branch of \eqref{nonsymmetric-scalar-mass},
define a positive decreasing radial profile $V$ by
\begin{equation}\label{nonsymmetric-radial-source}
 \begin{split}
 J_V(r)&=\frac1{d_0}
   \int_0^r s^{n-1}\zeta(s)^{q_2}\dd s,\\
 V(r)&=K_1+\int_r^RJ_V(s)^{1/a}
   s^{-\kappa-1}\dd s,
 \end{split}
\end{equation}
where $K_1>0$. Since $q_2\kappa<n$, the first integral is
finite. By \eqref{common-radial-identity},
$\mathcal Q[V]=\zeta^{q_2}$ and $J_V(0)=0$.
For the Hessian operator, $-V$ is $k$-admissible by
Lemma~\ref{radial-admissibility}.
Write
\[
 h_0=n-\kappa q_2=a(\kappa-\eta)>0.
\]
The derivative limit in \eqref{nonsymmetric-scalar-mass}
and two applications of l'H\^opital's rule in
\eqref{nonsymmetric-radial-source} give
\begin{equation}\label{nonsymmetric-V-asymptotic}
 \begin{split}
 J_V(r)&\sim
 \frac{r^n\zeta(r)^{q_2}}{d_0h_0},\\
 V(r)&\sim C_2^{-1/a}
 r^{\ell/a}\zeta(r)^{q_2/a}.
 \end{split}
\end{equation}
In particular $V\to\infty$ and
$V(r)=O(r^{-\eta})$, since $r^\kappa\zeta(r)\to0$.

For every $\varepsilon\in(0,1)$, \eqref{nonsymmetric-scalar-mass} gives
\[
 (1-\varepsilon)\zeta(r)^{q_2}
 \le\mathcal Q[w_2](r\theta)
 \le(1+\varepsilon)\zeta(r)^{q_2}
\]
for small $r$. Compare with $(1\pm\varepsilon)^{1/a}V$ using
Lemma~\ref{nonsymmetric-crosspoint}, with constants chosen on an outer
sphere, to obtain
\[
 (1-\varepsilon)^{1/a}V(r)-M_\varepsilon
 \le w_2(r\theta)\le
 (1+\varepsilon)^{1/a}V(r)+M_\varepsilon.
\]
Since $V\to\infty$, sending $r\to0$ and then
$\varepsilon\to0$ gives, uniformly in direction,
\begin{equation}\label{nonsymmetric-component-two}
 w_2(r\theta)\sim C_2^{-1/a}
 r^{\ell/a}\zeta(r)^{q_2/a}.
\end{equation}

Differentiating the ball mass gives, for almost every $t$,
\[
 -J'(t)=\frac{r^n}{\omega d_0}
 \int_{\Sph}w_2(r\theta)^{q_1}\dd\theta,
 \qquad r=e^{-t}.
\]
Using \eqref{nonsymmetric-component-two},
$q_1\eta=n$, and \eqref{nonsymmetric-scalar-mass}, we get
\begin{equation}\label{nonsymmetric-mass-ode}
 -J'(t)\sim c_*J(t)^{m_*},\qquad
 m_*=\frac{q_1q_2}{a^2}>1,\qquad
 c_*=\frac{C_2^{-q_1/a}}{d_0\kappa^{q_1q_2/a}}.
\end{equation}
Because $J(t)>0$ and $J(t)\to0$,
integration of $(J^{1-m_*})'\sim(m_*-1)c_*$ yields
\begin{equation}\label{nonsymmetric-mass-limit}
 J(t)\sim\bigl((m_*-1)c_*t\bigr)^{-1/(m_*-1)}.
\end{equation}
Equations \eqref{nonsymmetric-scalar-mass} and
\eqref{nonsymmetric-component-two} now imply
\eqref{nonsymmetric-limits}. Their positive coefficients obey
\[
 B^{q_1}=\frac{d_0a^2\kappa^a}{D}A^a=C_1A^a,
 \qquad A^{q_2}=d_0a\eta^a(\kappa-\eta)B^a=C_2B^a.
\]
Solving these identities gives
\eqref{nonsymmetric-constants} and proves the classification.

\subsection{Radial realization}

To realize these asymptotics, use a radial dynamical formulation
(cf.\ \cite{BG}). With $t=-\log r$, let $X,Y$ denote the logarithmic
derivatives of the components and $Z,H$ the negatives of those of
their radial fluxes. We construct a positive trajectory of
\begin{equation}\label{nonsymmetric-radial-system}
 \begin{aligned}
 \dot X&=X(\kappa-X-Z/a),&
 \dot Y&=Y(\kappa-Y-H/a),\\
 \dot Z&=Z(-n+q_1Y+Z),&
 \dot H&=H(-n+q_2X+H).
 \end{aligned}
\end{equation}
Its equilibrium
\[
 E=(\kappa,\eta,0,h_0),\qquad
 h_0=a(\kappa-\eta)>0
\]
has eigenvalues $-\kappa,-\eta,0,h_0$. The simple zero eigenvalue
has eigenvector
\[
 (-1/a,-q_2/a^2,1,q_2/a).
\]
The center-manifold theorem \cite{Carr} gives a $C^3$ invariant curve
parameterized by $Z$, with
\[
 \begin{aligned}
 X&=\kappa-Z/a+O(Z^2),\\
 Y&=\eta-q_2Z/a^2+O(Z^2),\\
 H&=h_0+q_2Z/a+O(Z^2).
 \end{aligned}
\]
The induced equation is
\begin{equation}\label{nonsymmetric-center-eqn}
 \dot Z=-\frac{D}{a^2}Z^2+O(Z^3).
\end{equation}
For small $Z(0)>0$, the other coordinates remain positive and
\eqref{nonsymmetric-center-eqn} gives a forward solution satisfying
\begin{equation}\label{nonsymmetric-center-rate}
 Z(t)>0,\qquad Z(t)\sim\frac{a^2}{Dt}.
\end{equation}
Reconstruct positive radial functions by
\begin{equation}\label{nonsymmetric-reconstruction}
 \begin{split}
 w_1(e^{-t})&=e^{\kappa t}
 \bigl(d_0^{a+q_1}Z^aH^{q_1}
 X^{a^2}Y^{aq_1}\bigr)^{1/D},\\
 w_2(e^{-t})&=e^{\eta t}
 \bigl(d_0^{q_2+a}Z^{q_2}H^a
 X^{aq_2}Y^{a^2}\bigr)^{1/D}.
 \end{split}
\end{equation}
Differentiating these expressions using
\eqref{nonsymmetric-radial-system} gives
\[
 \frac{\dd}{\dd t}\log w_1(e^{-t})=X,
 \qquad
 \frac{\dd}{\dd t}\log w_2(e^{-t})=Y.
\]
The defining algebra also gives
\[
 d_0ZX^aw_1^a=r^\ell w_2^{q_1},
 \qquad d_0HY^aw_2^a=r^\ell w_1^{q_2}.
\]
Consequently, the radial fluxes satisfy
\[
 \begin{alignedat}{2}
 J_{w_1}&=r^{a\kappa}(Xw_1)^a,
 &\qquad\dot{J}_{w_1}&=-ZJ_{w_1},\\[3pt]
 J_{w_2}&=r^{a\kappa}(Yw_2)^a,
 &\qquad\dot{J}_{w_2}&=-HJ_{w_2}.
 \end{alignedat}
\]
\eqref{common-radial-identity} proves the exact radial equations.
The rates \eqref{nonsymmetric-center-rate} and reconstruction
\eqref{nonsymmetric-reconstruction} give \eqref{nonsymmetric-limits}
and \eqref{nonsymmetric-constants}; both fluxes vanish at zero.

In the Hessian case, both signed components have positive radial
derivatives and strictly positive sources. Lemma~\ref{radial-admissibility} gives strict $k$-admissibility.
This completes the proof of Theorem~\ref{main-nonsymmetric}.

\subsection{Source masses and gradient asymptotics}

The source masses separate on the two power scales, with the
normalization of Section~\ref{sec:mass}.

\begin{corollary}\label{nonsymmetric-masses-gradients}
On the nonremovable branch of Theorem~\ref{main-nonsymmetric},
\begin{equation}\label{nonsymmetric-source-masses}
 \begin{split}
 \mu_1(B_r)&\sim\omega d_0(\kappa A)^aL(r)^{-a^2/D},\\
 \mu_2(B_r)&\sim\omega d_0(\eta B)^a
 r^{a(\kappa-\eta)}L(r)^{-aq_2/D}.
 \end{split}
\end{equation}
In particular,
\begin{equation}\label{nonsymmetric-mass-ratio}
 \frac{\mu_2(B_r)}{\mu_1(B_r)}
 \sim\left(\frac{\eta B}{\kappa A}\right)^a
 r^{a(\kappa-\eta)}L(r)^{a(a-q_2)/D}\to0.
\end{equation}
For $\mathcal Q=-\Delta_p$, one also has, uniformly in $\theta$,
\begin{equation}\label{nonsymmetric-gradient-limits}
 \begin{split}
 r^{\kappa+1}L(r)^{a/D}\nabla w_1(r\theta)&\to-\kappa A\theta,\\
 r^{\eta+1}L(r)^{q_2/D}\nabla w_2(r\theta)&\to-\eta B\theta.
 \end{split}
\end{equation}
\end{corollary}
\begin{proof}
Both atoms vanish by Lemma~\ref{nonsymmetric-upper}. Integrating the
uniform asymptotics \eqref{nonsymmetric-limits} in polar coordinates gives
\[
 \begin{split}
 \mu_1(B_r)&\sim\omega B^{q_1}
 \int_0^r s^{-1}L(s)^{-q_1q_2/D}\dd s
 =\frac{\omega D}{a^2}B^{q_1}L(r)^{-a^2/D},\\
 \mu_2(B_r)&\sim\omega A^{q_2}
 \int_0^r s^{h_0-1}L(s)^{-aq_2/D}\dd s
 \sim\frac{\omega A^{q_2}}{h_0}
 r^{h_0}L(r)^{-aq_2/D},
 \end{split}
\]
where $h_0=a(\kappa-\eta)>0$. The coefficient identities following
\eqref{nonsymmetric-mass-limit} prove
\eqref{nonsymmetric-source-masses}, and division gives
\eqref{nonsymmetric-mass-ratio}.

For the gradient limits, set $t=L(r)$ and
\[
 W_{1,r}(y)=r^\kappa t^{a/D}w_1(ry),\qquad
 W_{2,r}(y)=r^\eta t^{q_2/D}w_2(ry).
\]
On each fixed compact annulus these functions converge uniformly to
$A|y|^{-\kappa}$ and $B|y|^{-\eta}$, respectively. Their equations are
\[
 -\Delta_pW_{1,r}=t^{-1}W_{2,r}^{q_1},\qquad
 -\Delta_pW_{2,r}=W_{1,r}^{q_2}.
\]
Both right-hand sides are uniformly bounded on compact annuli.
Local $C^1$ compactness (Section~\ref{sec:local}) upgrades both
convergences to $C^1_{\loc}$. Evaluating their gradients on the
unit sphere gives \eqref{nonsymmetric-gradient-limits}.
\end{proof}

\begin{remark}\label{individual-mass-error-obstruction}
The error scale in Proposition~\ref{mass-approx} cannot in general
be replaced by the individual mass scale. Indeed, the nonremovable
solutions of Theorem~\ref{main-nonsymmetric} satisfy
\eqref{total-source}: the concave-power bound
\eqref{nonsymmetric-power-bound} controls $w_2^{q_1}$ by
$Cw_1^{q_S}$ near zero, while $q_2<q_S$ controls the other source.
With $P_2(t)=\mu_2(B_{e^{-t}})/(\omega d_0)$, the second formulas in
\eqref{nonsymmetric-limits} and \eqref{nonsymmetric-source-masses} give
\[
 \frac{\kappa^{-1}P_2(t)^{1/a}}
 {e^{-\kappa t}w_2(e^{-t}\theta)}\to\frac\eta\kappa\ne1
 \quad\text{uniformly in }\theta.
\]
Thus an $o(P_2^{1/a})$ error would give an incorrect leading coefficient.
\end{remark}

\appendix

\section{A bounded inverse for the nonconstant spherical part}
\label{app:cylinder}

Set
\[
 \mathcal L_0=\partial_t^2+\kappa\partial_t+c\Delta_{\Sph},
 \qquad c=\frac{\kappa+1}{n-1}.
\]
Here $c=1/(p-1)$ for $-\Delta_p$. The estimate used in
Lemma~\ref{angular-decay} excludes the constant mode and its zero root.

\begin{lemma}\label{cylinder-inverse}
Suppose $\mathbf w$ has zero spherical mean, uniformly bounded
$C^{2,\alpha}$ norms on unit strips of a half-cylinder, and solves
\[
 \mathcal L_0\mathbf w+\Pi\mathcal R\mathbf w=0,
\]
where $\mathcal R$ is a linear differential operator of order at most two
with matrix coefficients tending to zero in $C^{0,\alpha}$ on unit strips.
Then, for every $0<\sigma<\kappa+1$,
\[
 \|\mathbf w\|_{C^{2,\alpha}([t,t+1]\times\Sph)}=O(e^{-\sigma t}).
\]
\end{lemma}

\begin{proof}
Fix a starting point $T$ and define
\[
 \|h\|_{k,\alpha,\sigma;T}
 =\sup_{s\ge T}e^{\sigma(s-T)}
       \|h\|_{C^{k,\alpha}([s,s+1]\times\Sph)}.
\]
All spaces below have zero spherical mean. We construct the bounded
Dirichlet inverse $G_T$ of $\mathcal L_0$ at $t=T$ and prove

\begin{equation}\label{weighted-cylinder-inverse}
 \|G_Tf\|_{2,\alpha,\sigma;T}
 \le C_\sigma\|f\|_{0,\alpha,\sigma;T},
 \qquad 0\le\sigma<\kappa+1.
\end{equation}
The constants are independent of $T$.

Let $Y_{\ell,j}$ be an orthonormal spherical harmonic basis, with
$-\Delta_{\Sph}Y_{\ell,j}=\lambda_\ell Y_{\ell,j}$ and
$\lambda_\ell=\ell(\ell+n-2)$.  The zero-mean condition restricts us to
$\ell\ge1$.  In this component the equation is
\[
 h''+\kappa h'-(c\lambda_\ell)h=f_{\ell,j},
\]
whose roots are
\[
 \rho_\ell^\pm=
 \frac{-\kappa\pm\sqrt{\kappa^2+4c\lambda_\ell}}2.
\]
Since $c\lambda_1=\kappa+1$,
\begin{equation}\label{cylinder-root-gap}
 \rho_\ell^+\ge1,\qquad \rho_\ell^-\le-(\kappa+1),
 \qquad \rho_1^+=1,\quad \rho_1^-=-(\kappa+1).
\end{equation}
The scalar Dirichlet Green kernel on $[T,\infty)$ is
\[
 K_{\ell,T}(t,s)=-\frac1{\rho_\ell^+-\rho_\ell^-}
 \begin{cases}
 e^{\rho_\ell^-(t-s)}
  -e^{\rho_\ell^-(t-T)+\rho_\ell^+(T-s)},&t\ge s,\\
 e^{\rho_\ell^+(t-s)}
  -e^{\rho_\ell^-(t-T)+\rho_\ell^+(T-s)},&t<s.
 \end{cases}
\]
This kernel vanishes at $T$, is continuous at $t=s$, and has unit
derivative jump. The correction is bounded by the first term, so the
spherical spectral multiplier satisfies

\begin{equation}\label{green-kernel-bound}
 \|K_T(t,s)\|_{L^2(\Sph)\to L^2(\Sph)}\le C
 \begin{cases}
 e^{-(\kappa+1)(t-s)},&t\ge s,\\
 e^{-(s-t)},&t<s.
 \end{cases}
\end{equation}
Define $G_Tf(t)=\int_T^\infty K_T(t,s)f(s)\dd s$.
Minkowski's inequality and \eqref{green-kernel-bound} give the weighted
$L^2(\Sph)$ bound, with constants
$(\kappa+1-\sigma)^{-1}$ and $(1+\sigma)^{-1}$, and hence the weighted
$L^2$ bound on unit strips. For
$I_t=[t,t+1]\times\Sph$ and
$I_t^*=[t-1,t+2]\times\Sph$, the local estimate is
\[
 \|h\|_{C^{2,\alpha}(I_t)}\le
 C\bigl(\|\mathcal L_0h\|_{C^{0,\alpha}(I_t^*)}
             +\|h\|_{L^2(I_t^*)}\bigr),\qquad t\ge T+1.
\]
On the first strip use $[T,T+2]\times\Sph$ and add the
$C^{2,\alpha}$ norm of the boundary data, zero for $G_Tf$.
These Schauder estimates with local interpolation~\cite[Chapter~6]{GT}
have uniform constants; comparable weights on overlapping strips give
\eqref{weighted-cylinder-inverse}. Smooth approximation by the spherical
heat semigroup and in $t$, with uniformly bounded H\"older norms,
justifies regularity of the $L^2$ limit. For uniqueness, boundedness
excludes each growing mode, and the zero boundary value excludes the
remaining mode.

Similarly, for zero-mean boundary data $\varphi\in C^{2,\alpha}(\Sph)$,
the homogeneous extension
\[
 H_T\varphi(t)=\sum_{\ell\ge1,j}
 e^{\rho_\ell^-(t-T)}\varphi_{\ell,j}Y_{\ell,j}
\]
satisfies
\begin{equation}\label{weighted-cylinder-boundary}
 \|H_T\varphi\|_{2,\alpha,\sigma;T}
 \le C_\sigma\|\varphi\|_{C^{2,\alpha}(\Sph)}.
\end{equation}
This follows from the same $L^2$ bound and boundary and interior estimates.
All these arguments apply componentwise to vectors.

Let $\varepsilon_T$ be the supremum of the $C^{0,\alpha}$ norms of the
coefficients of $\mathcal R$ on strips beyond $T$.  The projection $\Pi$
is bounded on the cylinder H\"older spaces, and multiplication of H\"older
functions gives
\[
 \|\Pi\mathcal R h\|_{0,\alpha,\sigma;T}
 \le C\varepsilon_T\|h\|_{2,\alpha,\sigma;T}.
\]
We choose $T$ so large that the norm of $G_T\Pi\mathcal R$ is less than
$1/2$ in both the unweighted space and the space of weight $\sigma$.
For $\varphi=\mathbf w(T,\cdot)$, the equation
\begin{equation}\label{cylinder-contraction}
 h=H_T\varphi-G_T\Pi\mathcal R h
\end{equation}
has a unique weighted solution by a Neumann series, and that solution
is bounded. Uniqueness of the bounded Dirichlet inverse shows that
$\mathbf w$ satisfies the same equation in the unweighted space.
The unweighted contraction identifies the two solutions, proving
$\|\mathbf w\|_{2,\alpha,\sigma;T}<\infty$ without presupposing this bound.
\end{proof}

\begingroup
\small\linespread{1.0}\selectfont

\endgroup
\bigskip

\noindent Shiguang Ma: School of Mathematical Sciences and LPMC, Nankai University\\ 
Tianjin 300071, China \\
\textit{e-mail}: msgdyx8741@nankai.edu.cn \\

\medskip

\noindent Shengyang Zang: School of Mathematical Sciences, Nankai University\\ 
Tianjin 300071, China\\
\textit{e-mail}: 1120230032@mail.nankai.edu.cn
\end{document}